\pdfoutput=1
\documentclass{amsart}
\usepackage[utf8]{inputenc}

\usepackage[margin=1in]{geometry}
\usepackage[expansion=false]{microtype}

\usepackage{amsmath, amsfonts, amssymb, amsthm, amsopn}
\usepackage{eucal}

\usepackage{tikz}
\usetikzlibrary{diagrams}

\usepackage[pdfusetitle,unicode,hidelinks,bookmarksopen]{hyperref}

\newcommand{\resp}{{\sfcode`\.1000 resp.}}
\newcommand{\ie}{{\sfcode`\.1000 i.e.}}
\newcommand{\eg}{{\sfcode`\.1000 e.g.}}

\numberwithin{equation}{section}

\theoremstyle{plain}
\newtheorem*{theorem*}{Theorem}
\newtheorem{theorem}[equation]{Theorem}
\newtheorem{proposition}[equation]{Proposition}
\newtheorem{lemma}[equation]{Lemma}
\newtheorem{corollary}[equation]{Corollary}

\theoremstyle{definition}
\newtheorem{definition}[equation]{Definition}
\newtheorem{example}[equation]{Example}

\theoremstyle{remark}
\newtheorem{remark}[equation]{Remark}

\let\scr=\mathcal
\let\bb=\mathbb
\def\N{\bb N}
\def\Z{\bb Z}
\def\Q{\bb Q}

\def\A{\bb A}
\def\P{\bb P}
\def\V{\bb V}
\def\1{\mathbf 1}

\def\G{\mathbb G}
\def\ph{\mathord-}
\def\pt{{\scriptscriptstyle\bullet}}

\let\from=\leftarrow
\let\into=\hookrightarrow
\let\onto=\twoheadrightarrow
\let\tens=\otimes
\def\suchthat{\:\vert\:}

\DeclareMathOperator{\Sym}{Sym}
\def\id{\mathrm{id}}
\DeclareMathOperator{\Corr}{Corr}
\def\all{\mathrm{all}}
\def\cpt{\mathrm{comp}}
\def\open{\mathrm{open}}
\def\prop{\mathrm{prop}}
\DeclareMathOperator{\Hom}{Hom}

\DeclareMathOperator{\Map}{Map{}}

\DeclareMathOperator{\Spec}{Spec}
\DeclareMathOperator{\Proj}{Proj}
\def\Th{\mathrm{Th}}

\def\Nis{\mathrm{Nis}}

\def\mot{\mathrm{mot}}
\def\et{\mathrm{\acute et}}

\def\QCoh{\mathrm{QCoh}{}}
\def\Ex{\mathit{Ex}}

\def\Sph{\mathrm{Sph}{}}

\def\H{\mathrm H}
\def\SH{\mathrm S\mathrm H}
\def\Mod{\mathrm{M}\mathrm{od}{}}
\def\Sch{\mathrm{S}\mathrm{ch}{}}
\def\Sm{\mathrm{S}\mathrm{m}}
\def\Et{\mathrm{E}\mathrm{t}}
\def\s{\mathrm{S}}
\def\op{\mathrm{op}}
\def\Pr{\mathcal{P}\mathrm{r}}
\def\CAlg{\mathrm{CAlg}}

\def\minus{\smallsetminus}
\def\longinto{\lhook\joinrel\longrightarrow}

\def\Cat{\mathcal{C}\mathrm{at}{}}
\def\PSh{\mathcal{P}}

\def\Aff{\mathrm{Aff}}
\def\GL{\mathrm{GL}}
\def\htp{\mathrm{htp}}
\def\Tw{\mathit{Tw}{}}
\def\KH{K\mskip -1mu H}

\DeclareMathOperator{\Stab}{Stab}

\let\segre=\sigma

\let\lim=\relax
\DeclareMathOperator*{\lim}{lim}
\DeclareMathOperator*{\colim}{colim}

\title{The six operations in equivariant motivic homotopy theory}
\author{Marc Hoyois}
\date{\today}
\address{Department of Mathematics, Massachusetts Institute of Technology, Cambridge, MA, USA}
\email{hoyois@mit.edu}
\urladdr{\url{http://math.mit.edu/~hoyois/}}

\begin{document}

\begin{abstract}
We introduce and study the homotopy theory of motivic spaces and spectra parametrized by quotient stacks $[X/G]$, where $G$ is a linearly reductive linear algebraic group. We extend to this equivariant setting the main foundational results of motivic homotopy theory: the (unstable) purity and gluing theorems of Morel–Voevodsky and the (stable) ambidexterity theorem of Ayoub.
Our proof of the latter is different than Ayoub's and is of interest even when $G$ is trivial.
Using these results, we construct a formalism of six operations for equivariant motivic spectra, and we deduce that any cohomology theory for $G$-schemes that is represented by an absolute motivic spectrum satisfies descent for the cdh topology.
\end{abstract}

\maketitle
\vspace{-1em}
\tableofcontents
\vspace{-3em}
\newpage
\section{Introduction}

\begingroup
\def\nocontent#1#2#3{}

The goal of this paper is to develop the formalism of six operations in stable equivariant motivic homotopy theory. An equivariant version of motivic homotopy theory was first considered by Voevodsky in \cite{DeligneNote} and played a small part in his proof of the Bloch–Kato conjecture \cite{Voevodsky:2008}, more precisely in the construction of symmetric power functors on the $\A^1$-homotopy category \cite[\S2.1]{MEMS}. A variety of definitions of equivariant motivic homotopy theory were later proposed by several authors: by Hu, Kriz, and Ormsby \cite{HuKrOr}, by Herrmann \cite{HerrmannThesis,Herrmann}, by Heller, Krishna, and Østvær \cite{HeKrOs}, and by Carlsson and Joshua \cite{CarlssonJoshua}.
In these approaches, equivariant motivic homotopy theory is a setting in which to study invariants of smooth $G$-schemes with some specific properties, for $G$ an algebraic group. 
We will discuss in \S\ref{sub:comparisons} below how these approaches relate to the one developed in this paper.
In any case, our starting point is somewhat different: we view equivariant homotopy theory as a natural extension of parametrized homotopy theory, and the formalism of six operations serves as a guiding principle in our definitions.

The ``yoga'' of six operations is due to Grothendieck and was first developed in \cite{SGA4-3}, in collaboration with Artin and Deligne, for the étale cohomology of schemes.
The operations in question are the pushforward along a morphism and its left adjoint, the compactly supported pushforward and its right adjoint, the tensor product, and the corresponding internal homomorphism object. These six functors are related by several identities, such as those appearing in Theorem~\ref{thm:intro} below. Despite their mundane appearance, these identities subsume and unify several nontrivial theorems, such as Poincaré duality and the Lefschetz trace formula.
The formalism of six operations was later developed for many other coefficient systems in various geometric contexts, and in particular by Ayoub for the stable motivic homotopy theory of schemes \cite{Ayoub}.

\subsection{Equivariant homotopy theory}
\label{sub:equivariant}

For $S$ a scheme, Voevodsky constructed the $\infty$-category $\SH(S)$ of motivic spectra over $S$ \cite{Voevodsky:1998}. It is in many ways an algebraic analog to the $\infty$-category $\SH^\mathrm{top}(S)$ of sheaves of spectra over a topological space $S$. In particular, both constructions support a formalism of six operations as the base $S$ varies. 

In topology, a basic observation that leads us from parametrized to equivariant homotopy theory is that the objects being parametrized (\eg{}, CW complexes or smooth manifolds) have nontrivial automorphisms and hence can vary in families parametrized by topological \emph{stacks} rather than just spaces. We can therefore expect that the functor $S\mapsto \SH^\mathrm{top}(S)$ admits an interesting extension
\begin{tikzmath}
	\diagram{
	\{\text{topological spaces}\} & \{\text{$\infty$-categories}\}\rlap. \\
	\{\text{topological stacks}\} & \\
	};
	\arrows (11-) edge node[above]{$\SH^\mathrm{top}$} (-12) (11) edge[c->] (21) (21) edge[dashed] (12);
\end{tikzmath}
Similarly, smooth schemes can vary in families parametrized by algebraic stacks, and hence we might expect an interesting extension of $S\mapsto\SH(S)$ to algebraic stacks:
\begin{tikzmath}
	\diagram{
	\{\text{schemes}\} & \{\text{$\infty$-categories}\}\rlap. \\
	\{\text{algebraic stacks}\} & \\
	};
	\arrows (11-) edge node[above]{$\SH$} (-12) (11) edge[c->] (21) (21) edge[dashed] (12);
\end{tikzmath}
Of course, it is not hard to construct such extensions, and there are even several sensible possibilities. If a topological stack is presented by a simplicial topological space $X_\bullet$, one can consider the limit
\[
\SH^\mathrm{top}(X_\bullet)=\lim_{n\in\Delta} \SH^\mathrm{top}(X_n).
\]
For example, if $G$ is a discrete group and $B_\bullet G$ is the usual bar construction on $G$, then $\SH^\mathrm{top}(B_\bullet G)$ is the $\infty$-category of spectra equipped with a homotopy coherent action of $G$.
It turns out that $\SH^\mathrm{top}(X_\bullet)$ depends only on the topological stack presented by $X_\bullet$.
In the same way one can plug into $\SH(\ph)$ arbitrary simplicial schemes (and indeed arbitrary small diagrams of schemes); this construction was studied in an axiomatic setting by Ayoub in \cite[\S2.4]{Ayoub}. 
The problem with
this naive extension of $\SH(\ph)$ is that, unlike in topology, it depends on presentations of stacks rather than on stacks themselves, because $\SH(\ph)$ does not satisfy étale descent.
Another extension of $\SH^\mathrm{top}(\ph)$ to topological stacks with a more geometric flavor is the so-called ``naive'' stable equivariant homotopy theory.
This coefficient system includes finer cohomological invariants of topological stacks, like topological $K$-theory. Its main defect is that it does not satisfy Atiyah duality, which is an important feature of a formalism of six operations.
This is rectified by passing to
``genuine'' stable equivariant homotopy theory.
Although one could give a very general definition,
this theory only works well for nice enough stacks, such as quotients of topological spaces by actions of compact Lie groups.
The principal example is the stack $\mathbf BG=[*/G]$ for $G$ a compact Lie group, in which case $\SH^\mathrm{top}(\mathbf BG)$ is the usual $\infty$-category of genuine $G$-spectra. A partial formalism of six operations in this context, encompassing only locally constant coefficients, was developed by Hu \cite{HuEquiv} and May--Sigurdsson \cite{MaySigurdsson}.

Our goal in this paper is to extend the assignment $S\mapsto \SH(S)$, together with its formalism of six operations, from schemes to a suitable class of algebraic stacks. In other words, our goal is to develop a good theory of ``genuine'' equivariant motivic spectra over varying base stacks.
We will be able to achieve this,
essentially, for
 stacks with affine stabilizers whose derived category is compactly generated by vector bundles. 
 These are somewhat analogous to quotients of spaces by compact Lie groups, in that they have a well-behaved representation theory compared to more general stacks.

At this point we need to come clean about the fact that the word ``stack'' will not much appear beyond this introduction. Indeed, the above-mentioned restriction on stacks implies that we will not lose much generality by considering only global quotient stacks $[X/G]$ for a fixed algebraic group $G$ (e.g., $\GL_n$ for large $n$). It will therefore be much simpler to work directly with $G$-schemes rather than the corresponding quotient stacks.
For a $G$-scheme $X$, we will write $\SH^G(X)$ for $\SH([X/G])$.
 The reader will rightly object that the category of $G$-schemes for fixed $G$ only accounts for morphisms of stacks that are schematic (i.e., representable by schemes). However, this is not a significant drawback because the interesting features of our formalism of six operations only exist for schematic morphisms.
 For example, we do not expect an exceptional adjunction $f_!\dashv f^!$ or a left adjoint $f_\sharp$ to $f^*$ unless $f$ is schematic.
These restrictions already exist in topology: if $G$ is a nontrivial compact Lie group, pullback along the morphism $f\colon \mathbf BG\to *$ does not have a left adjoint (in other words, there are no ``genuine $G$-orbits''). Enforcing the existence of such left adjoints naturally leads to a parametrized version of \emph{global} homotopy theory in the sense of Schwede \cite{Schwede},
which should also have a motivic analog, 
but we will not discuss it further here.

It is worth noting that all the intricacies of equivariant homotopy theory disappear in the \emph{étale} version of the theory. Indeed, $\SH^\et(\ph)$ satisfies descent for the smooth topology and hence, if $\mathfrak X$ is Artin stack (or even an Artin $\infty$-stack) presented by a simplicial scheme $X_\bullet$ with smooth face maps, it is perfectly adequate to define $\SH^\et(\mathfrak X)$ as the limit $\lim_{n\in\Delta}\SH^\et(X_n)$. The result is a theory that is already ``genuine'' and ``global''.
 This is why, modulo some serious technicalities, it is straightforward to extend the formalism of six operations in $\ell$-adic cohomology from schemes to Artin $\infty$-stacks (see \cite{LiuZheng2}).
Our work is thus motivated by invariants such as algebraic $K$-theory, Chow groups, and algebraic cobordism, which do not satisfy étale descent.

\subsection{The formalism of six operations}

Fix a quasi-compact quasi-separated base scheme $B$ and a flat finitely presented group scheme $G$ over $B$.
Throughout this introduction, we will assume for simplicity that $B$ has the $G$-resolution property, \ie, that every finitely generated quasi-coherent $G$-module over $B$ is the quotient of a locally free one.

As explained in \S\ref{sub:equivariant}, our objective is to construct a functor $S\mapsto \SH^G(S)=\SH([S/G])$ from $G$-schemes to $\infty$-categories,
together with a formalism of six operations.
We are only able to develop a good theory under the assumption that $G$ is \emph{tame} (see Definition~\ref{def:tame}).
The following are the essential examples of tame group schemes:
\begin{itemize}
	\item $G$ is finite locally free of order invertible on $B$;
	\item $G$ is of multiplicative type;
	\item $G$ is reductive and $B$ has characteristic $0$ (\ie, there exists $B\to \Spec\Q$).
\end{itemize}
Moreover, for the same reason as in \cite[\S1.3.5]{Ayoub}, we will restrict our attention to $G$-schemes that are \emph{$G$-quasi-projective}, \ie, that admit a $G$-equivariant immersion into the projectivization of a $G$-equivariant vector bundle over $B$.
If $G$ is finite locally free, any quasi-projective $G$-scheme is $G$-quasi-projective, and if $G$ is reductive, any \emph{normal} quasi-projective $G$-scheme is $G$-quasi-projective, by Sumihiro's theorem \cite{Sumihiro}.
In the nonequivariant setting, it was shown in \cite[\S2]{CD} how to extend the formalism of six operations from quasi-projective to more general schemes using Chow's lemma. Similar tricks will work in our setting, but unfortunately the reach of the equivariant versions of Chow's lemma is much more limited. For this reason, we will only discuss this generalization when $G$ is discrete (see Remark~\ref{rmk:nagatachow}).

For $G$ a tame group and $S$ a finitely presented $G$-quasi-projective scheme (or an arbitrary qcqs $G$-scheme if $G$ is discrete), we will construct a closed symmetric monoidal stable $\infty$-category $\SH^G(S)$. In particular, $\SH^G(S)$ is equipped with a tensor product $\tens$ and an internal mapping object $\Hom$.
For every $G$-equivariant morphism $f\colon T\to S$, we have a pullback–pushforward adjunction
\[
f^*:\SH^G(S) \rightleftarrows \SH^G(T) : f_*,
\]
where $f^*$ is symmetric monoidal. If $f$ is smooth, then $f^*$ also admits a left adjoint $f_\sharp$.
If $f$ is separated and of finite type, we further have the ``exceptional'' adjunction
\[
f_! :\SH^G(T) \rightleftarrows \SH^G(S) : f^!,
\]
where $f_!$ is \emph{pushforward with compact support}: for any factorization $f=p\circ j$ where $j$ is an open immersion and $p$ is proper, we have $f_!\simeq p_*\circ j_\sharp$.
Our formalism of six operations is summarized in the following theorem:

\begin{theorem}[Theorem~\ref{thm:main} and Proposition~\ref{prop:constructible}]
\label{thm:intro}
Let $B$ be a qcqs scheme and $G$ a tame group scheme over $B$. 
If $G$ is not finite, we assume that $B$ has the $G$-resolution property.
Then the six operations 
\[(\ph)^*,\; (\ph)_*,\; (\ph)_!,\; (\ph)^!,\; \tens,\; \Hom\]
satisfy the following properties on finitely presented $G$-quasi-projective $B$-schemes
(or on all qcqs $G$-schemes if $G$ is discrete), whenever the exceptional functors are defined.
	\begin{enumerate}
		\item \textnormal{(Proper pushforward)} If $f$ is a proper $G$-morphism, there is an equivalence
		\[
		f_!\simeq f_*.
		\]
		\item \textnormal{(Smooth pullback)} If $f$ is a smooth $G$-morphism, there is a self-equivalence $\Tw_f$ and an equivalence
		\begin{gather*}
		\Tw_f\circ f^!\simeq f^*.
		\end{gather*}
		\item \textnormal{(Base change)} If
		\begin{tikzmath}
				\diagram{\bullet & \bullet \\ \bullet & \bullet \\};
				\arrows (11-) edge node[above]{$g$} (-12) (11) edge node[left]{$q$} (21) (21-) edge node[below]{$f$} (-22) (12) edge node[right]{$p$} (22);
			\end{tikzmath}
			is a cartesian square of $G$-schemes, there are equivalences
			\begin{gather*}
			f^*p_!\simeq q_!g^*,\\
			f^!p_*\simeq q_*g^!.
			\end{gather*}
		\item \textnormal{(Gluing)} If $i$ is a closed $G$-immersion with complementary open $G$-immersion $j$, there are cofiber sequences
		\begin{gather*}
		j_!j^! \to \id \to i_*i^*,\\
		i_!i^! \to \id \to j_*j^*.
		\end{gather*}
		\item \textnormal{(Immersive pushforward)} If $i$ is a $G$-immersion, the functors $i_*$ and $i_!$ are fully faithful.
		\item \textnormal{(Monoidality)} If $f$ is any $G$-morphism, there is an equivalence
		\[
		f^*(\ph\tens\ph)\simeq f^*(\ph)\tens f^*(\ph).
		\]
		\item \textnormal{(Projection formulas)} If $f$ is any $G$-morphism, there are equivalences
		\begin{gather*}
		f_!(\ph\tens f^*(\ph))\simeq f_!(\ph)\tens \ph,\\
		\Hom(f_!(\ph),\ph)\simeq f_*\Hom(\ph, f^!(\ph)),\\
		f^!\Hom(\ph,\ph)\simeq \Hom(f^*(\ph), f^!(\ph)).
		\end{gather*}
		\item \textnormal{(Homotopy invariance)} If $f$ is a $G$-affine bundle, the functors $f^*$ and $f^!$ are fully faithful.
		\item \textnormal{(Constructible separation)} If $\{f_i\}$ is a cover for the $G$-equivariant constructible topology, the families of functors $\{f_i^*\}$ and $\{f_i^!\}$ are conservative.
	\end{enumerate}
\end{theorem}

If $G$ is finite locally free, standard model-categorical constructions provide a presentation of $\SH^G(S)$ by a symmetric monoidal simplicial model category that is combinatorial and left proper. For more general $G$, however, we do not know an explicit such presentation, although there exists one by the main result of \cite{NikolausSagave}. The language of $\infty$-categories is also useful to make sense of the following statement, which is an easy corollary of Theorem~\ref{thm:intro}:

\begin{theorem}[Proposition~\ref{prop:cdh}]
	Under the assumptions of Theorem~\ref{thm:intro}, the assignment $S\mapsto \SH^G(S)$, $f\mapsto f^*$, is a sheaf for the $G$-equivariant cdh topology.
\end{theorem}

Let us give some brief commentary on the proof of Theorem~\ref{thm:intro}. In the nonequivariant case, Ayoub \cite{Ayoub} and Cisinski–Déglise \cite{CD} have shown how to derive such a formalism from just a few key results. 
Similarly, Theorem~\ref{thm:intro} is reduced by abstract nonsense to two nontrivial results: \emph{gluing} for complementary open–closed pairs and \emph{ambidexterity} for smooth proper morphisms.
A third important result is \emph{purity} for smooth closed pairs: although purity does not directly enter the proof of Theorem~\ref{thm:intro}, it plays a significant role in our proof of ambidexterity; it also provides an explicit description of the twisting equivalence $\Tw_f$, showing that it depends only on the relative tangent bundle of $f$.
We will discuss these three key results in more details in \S\ref{sub:comparisons} below, along with the actual definition of $\SH^G(S)$.

In characteristic zero, the restriction to tame group schemes is not very significant.
Indeed, a theorem of Gross \cite[Theorem A]{Gross} states that any qcqs Artin stack with affine stabilizers and the resolution property is the quotient of a quasi-affine scheme by an action of $\GL_n$. Thus, if $B$ has characteristic zero, our formalism of six operations includes all finitely presented Artin stacks over $B$ with affine stabilizers and the resolution property. In arbitrary characteristic, it includes tame Deligne–Mumford stacks that are quotients of quasi-projective schemes by finite étale groups, but also some Artin stacks that are not Deligne–Mumford (\eg, quotients by tori). Unfortunately, we do not know how to set up a satisfying theory for more general stacks in positive characteristic.

The formalism of six operations described above is not the most complete possible as it does not mention dualizing objects, constructibility, and absolute purity. In ordinary stable motivic homotopy theory, a good theory of dualizing and constructible objects requires some form of resolutions of singularities \cite[\S2.3.10]{Ayoub}. It is likely that the arguments in \emph{loc.\ cit.} can be applied in the present context assuming the existence of suitable equivariant resolutions of singularities, which are known to exist over fields of characteristic zero.
Absolute purity, on the other hand, seems out of reach since it is not even known for $\SH(\ph)$.

\subsection{Summary of the construction}
\label{sub:comparisons}

Let $G$ be a tame group scheme over $B$ and let $S$ be a finitely presented $G$-quasi-projective $B$-scheme.
As in the nonequivariant case, our construction of $\SH^G(S)$ proceeds in several steps. First, we define the unstable equivariant motivic homotopy $\infty$-category $\H^G(S)$ as a localization of the $\infty$-category of presheaves on the category $\Sm_S^G$ of smooth $G$-quasi-projective $S$-schemes. 
An object in $\H^G(S)$ is thus a presheaf $F$ on $\Sm_S^G$, and it is subject to the following two conditions:
\begin{itemize}
	\item (Homotopy invariance) If $Y\to X$ is a $G$-equivariant torsor under a $G$-vector bundle, then the restriction map $F(X)\to F(Y)$ is an equivalence.
	\item (Nisnevich excision) $F(\emptyset)\simeq *$, and if $V\to X$ is an étale map in $\Sm_S^G$ that is an isomorphism over a $G$-invariant finitely presented closed subscheme $Z\subset X$, then the following square is cartesian:
	\begin{tikzmath}
		\diagram{F(X) & F(X\minus Z) \\ F(V) & F(V\minus Z)\rlap. \\};
		\arrows (11-) edge (-12) (11) edge (21)
		(21-) edge (-22) (12) edge (22);
	\end{tikzmath}
\end{itemize}
In the context of equivariant algebraic $K$-theory, what we call homotopy invariance is sometimes called \emph{strong} homotopy invariance to distinguish it from the weaker property of $\A^1$-homotopy invariance. There are several reasons for using the stronger condition: perhaps the most important one is that it plays a crucial role in our proof of ambidexterity; it also implies that $\H^G(S)$ is generated by the affine schemes in $\Sm_S^G$, a fact that is used in the proofs of all the main results.
The Nisnevich excision condition was first considered by Voevodsky for finite locally free groups in \cite{DeligneNote}, where it was shown to be equivalent to a topological descent condition.
It was further studied in \cite{Herrmann,HVO,HeKrOs}.

When $B$ is noetherian of finite Krull dimension and $G$ is a finite discrete group, the unstable equivariant motivic homotopy $\infty$-category $\H^G(B)$ is equivalent to the underlying $\infty$-category of the model category constructed by Heller, Krishna, and Østvær in \cite{HeKrOs}; see Remark~\ref{rmk:comparison}. For more general tame groups $G$, however, we do not know if they agree. In any case, our unstable category does not appear to be equivalent to the constructions in \cite{HuKrOr} and \cite{Herrmann}, where a finer version of the Nisnevich topology is used.

The following are our two main results about $\H^G(S)$. They are equivariant generalizations of the purity and gluing theorems of Morel and Voevodsky \cite[\S3, Theorems 2.23 and 2.21]{MV}.

\begin{theorem}[Theorem~\ref{thm:wexcisive}]
	\label{thm:intro_purity}
	Let $S$ be a $G$-scheme and let $Z\into X$ be a $G$-equivariant closed immersion in $\Sm_S^G$. Deformation to the normal bundle induces a canonical equivalence
	\[
	\frac{X}{X\minus Z} \simeq \frac{N_ZX}{N_ZX\minus Z}
	\]
	in $\H^G(S)$, where $N_ZX$ is the normal bundle of $Z$ in $X$ with the induced action of $G$.
\end{theorem}

\begin{theorem}[Theorem~\ref{thm:gluing}]
	\label{thm:intro_gluing}
	Let $i\colon Z\into S$ be a $G$-equivariant closed immersion with open complement $j\colon U\into S$. Then, for every $F\in\H^G(S)$, there is a cocartesian square
	\begin{tikzmath}
		\def\colsep{2em}
		\diagram{j_\sharp j^* F & F \\ U & i_*i^* F\rlap. \\};
		\arrows (11-) edge (-12) (11) edge(21) (21-) edge (-22) (12) edge (22);
	\end{tikzmath}
\end{theorem}

Theorem~\ref{thm:intro_purity} was also proved in \cite[Theorem 7.6]{HeKrOs} assuming that $S$ is the spectrum of a perfect field, that $G$ is a finite discrete abelian group acting trivially on $S$, and that $S$ contains a primitive $e$th root of unity, where $e$ is the exponent of $G$.

Our proofs of Theorems \ref{thm:intro_purity} and~\ref{thm:intro_gluing} follow the same ideas as in \cite[\S3.2]{MV}. The main obstacle, which is the source of the restrictions in the statement of the purity theorem in \cite{HeKrOs}, is that a smooth $G$-equivariant morphism is not locally the composition of an étale map and a vector bundle, even for $G$ finite discrete. For gluing, an additional complication comes from the fact that, if $G$ is not discrete, the inclusion of the subcategory of affine $G$-schemes in $\Sm_S^G$ is not necessarily cocontinuous for the Nisnevich topology. We remark that these two theorems already fully exploit the definition of $\H^G(S)$, in the sense that no obvious weakening of the conditions of homotopy invariance and Nisnevich excision would make them work.

Let $\H_\pt^G(S)$ denote the undercategory $\H^G(S)_{S/}$, \ie, the $\infty$-category of \emph{pointed} presheaves on $\Sm_S^G$ that are homotopy invariant and Nisnevich excisive. It is a symmetric monoidal $\infty$-category under the ``smash product'', which we denote by $\tens$.
 For every locally free $G$-module of finite rank $\scr E$ on $S$, there is a corresponding ``representation sphere'' $\s^{\scr E}\in \H^G_\pt(S)$ and a suspension functor $\Sigma^{\scr E}$ defined by 
\[
\s^{\scr E}=\frac{\V(\scr E)}{\V(\scr E)\minus S}\quad\text{and}\quad \Sigma^{\scr E}=\s^{\scr E}\tens(\ph),
\]
where $\V(\scr E)=\Spec(\Sym(\scr E))$.
The stable equivariant motivic homotopy $\infty$-category $\SH^G(S)$ is then defined by formally adjoining a $\tens$-inverse $\s^{-\scr E}$ for each representation sphere $\s^{\scr E}$:
\[
\SH^G(S) = \H_\pt^G(S)[\s^{-\scr E} \suchthat \text{$\scr E$ is a locally free $G$-module of finite rank on $S$}].
\]
We denote by $\Sigma^\infty\colon \H^G_\pt(S)\to\SH^G(S)$ the canonical functor.
The precise definition of $\SH^G(S)$ is a universal construction in the $\infty$-category of presentably symmetric monoidal $\infty$-categories, and we do not spell it out here. We only note that the symmetric monoidal $\infty$-category of genuine $G$-spectra can be obtained in the same way from that of pointed $G$-spaces. In the equivariant motivic setting, this stabilization procedure recovers the ones considered in \cite{HuKrOr} and \cite{Herrmann} (for finite discrete groups), but it has the advantage that the resulting $\infty$-category $\SH^G(S)$ comes equipped with a homotopy coherent symmetric monoidal structure.

Our final main result is what we call \emph{ambidexterity}, since it identifies, up to a twist, the left and right adjoints of the base change functor:

\begin{theorem}[Theorem~\ref{thm:stableduality}]
	\label{thm:intro_duality}
	Let $f\colon X\to S$ be a smooth proper $G$-morphism. Then there is a canonical equivalence $f_*\simeq f_\sharp\Sigma^{-\Omega_f}$, where $f_\sharp$ and $f_*$ are left and right adjoint to $f^*\colon \SH^G(S)\to \SH^G(X)$.
\end{theorem}

Here, $\Omega_f$ is the sheaf of differentials of $X$ over $S$, which is a locally free $G$-module of finite rank on $X$.
An easy consequence is \emph{Atiyah duality}, which states that $\Sigma^\infty X_+$ is strongly dual to the Thom spectrum of the stable normal bundle of $X$ over $S$ in the symmetric monoidal $\infty$-category $\SH^G(S)$.

In the nonequivariant case, this ambidexterity theorem was proved by Ayoub \cite[Théorème 1.7.17]{Ayoub}. His argument implicitly relies on the fact that any vector bundle is Zariski-locally a sum of line bundles. The same argument can be used to prove Theorem~\ref{thm:intro_duality} if $X$ can be embedded into a projective bundle $\P(\scr E)$ where $\scr E$ is a sum of $G$-line bundles over $S$. This suffices to prove Atiyah duality in the stable equivariant motivic homotopy category of a diagonalizable group over an algebraically closed field, but it is far from sufficient in general.
Our strategy to prove Theorem~\ref{thm:intro_duality} is to explicitly write down a unit and a counit for an adjunction between $f^*$ and $f_\sharp\Sigma^{-\Omega_f}$, and to verify the triangle identities. This is rather interesting even in the case of a trivial group, as Ayoub's proof does not provide a geometric description of this adjunction. The central construction is an algebro-geometric analog of the Pontryagin–Thom collapse map, which is essentially due to Voevodsky \cite[Theorem 2.11]{VV}, although he only used it to prove a pale motivic-cohomological shadow of Atiyah duality.

In fact, we will prove a finer \emph{unstable} version of Theorem~\ref{thm:intro_duality} (Theorem~\ref{thm:duality}), which shows that the passage from $\H_\pt^G(\ph)$ to $\SH^G(\ph)$ is exactly what is needed to enforce ambidexterity for smooth projective morphisms. For instance, if $f\colon \P^n_S\to S$ is the projection, we will show that the unit and counit for the adjunction $f^*\dashv f_\sharp\Sigma^{-\Omega_f}$, as well as homotopies witnessing the triangle identities, are in the image of $\Sigma^\infty$ after tensoring with $\s^{\scr O}$ at least $2n^2+4n+1$ times.

\subsection{Some applications}

\subsubsection{Motivic Wirthmüller and Adams isomorphisms}

As was observed by Hu \cite[\S6]{HuEquiv}, the Wirthmüller and Adams isomorphisms in stable equivariant homotopy theory are merely instances of the ambidexterity isomorphism $f_*\simeq f_\sharp\Sigma^{-\Omega_f}$. One can also consider these instances in the motivic context. However, they are not isomorphisms in general, because quotients of tame groups are not always proper.

Let $G$ be a tame group scheme over $B$ and let $H\subset G$ be a subgroup.
Suppose that the quotient $B$-scheme $G/H$ exists and is smooth and $G$-quasi-projective over $B$.
Let $p\colon G/H\to B$ be the structural $G$-morphism, presenting the morphism of stacks $\mathbf BH\to \mathbf BG$. The functor $p^*\colon \SH^G(B) \to \SH^G(G/H)\simeq\SH^H(B)$ is thus the ``forgetful'' functor from $G$-spectra to $H$-spectra, and its left and right adjoints $p_\sharp$ and $p_*$ are the induction and coinduction functors, respectively.
The Wirthmüller morphism is then the composition
\[p_\sharp\Sigma^{-\Omega_{G/H}} \simeq p_!\to p_*;\]
it is an equivalence if $G/H$ is proper. See \cite[Theorem 1.1]{MayWirthmuller} for a statement of the classical Wirthmüller isomorphism in this form.

The Adams morphism is more subtle as it involves a nonschematic pushforward (see Remark~\ref{rmk:stacks}).
Let $H\subset G$ be as above and let $N$ be a smooth normal subgroup of $G$ such that $G/N$ is tame and $N\cap H$ is trivial. Assume further that $G/NH$ is a $G/N$-quasi-projective $B$-scheme. 
We then have smooth morphisms of stacks
\[
\mathbf BH \stackrel p\to \mathbf BG \stackrel q\to \mathbf B(G/N),
\]
with $p$ and $qp$ schematic, and a canonical map $(qp)_! \to q_*p_!$. 
Let $\mathfrak n$ be the Lie algebra of $N$ with conjugation action of $G$, viewed as a vector bundle over $\mathbf BG$.
The canonical fiber sequence of cotangent complexes
over $\mathbf BH$ induces an equivalence $\Sigma^{-\Omega_p}\Sigma^{\Omega_{qp}}\simeq \Sigma^{-p^*(\mathfrak n)}$.
The Adams morphism for $H$-induced $G$-spectra is then the composition
\[(qp)_\sharp\simeq (qp)_!\Sigma^{\Omega_{qp}}\to q_*p_!\Sigma^{\Omega_{qp}}\simeq q_*p_\sharp \Sigma^{-\Omega_{p}}\Sigma^{\Omega_{qp}}\simeq q_* \Sigma^{-\mathfrak n}p_\sharp;\]
 it is an equivalence if $G/H$ is proper.
In classical notation, given an $H$-spectrum $X$, the above morphism reads $G/N\ltimes_H X\to (\Sigma^{-\mathfrak n}(G\ltimes_HX))^N$.
As in the classical case, one can do slightly better and construct an Adams morphism for suitably defined ``$N$-free $G$-spectra'', but this requires some variants of our results.
 
\subsubsection{An equivariant Lefschetz trace formula}

The formalism of six operations developed here implies an abstract version of the Lefschetz trace formula, as in \cite[Theorem 1.3]{HoyoisGLV}. Let $S$ be a $G$-scheme, let $p\colon X\to S$ be a smooth $G$-projective morphism, let $f\colon X\to X$ be a $G$-equivariant endomorphism over $S$, and let $i\colon X^f\into X$ be the inclusion of the scheme of fixed points of $f$. Suppose that $X^f$ is smooth over $S$ and that the endomorphism $\id - i^*(df)$ of $i^*(\Omega_{X})$ restricts to an automorphism $\phi$ of the conormal sheaf $\scr N_i$. 
Then the trace of $\Sigma^\infty f_+$ in $\SH^G(S)$ is equal to the trace of the automorphism of $\Sigma^\infty X^f_+$ induced by $\phi$ via the $J$-homomorphism.

When $B$ is the spectrum of a field and $G$ is finite and discrete, the group of endomorphisms of the $G$-equivariant motivic sphere spectrum over $B$ was recently computed by Gepner and Heller \cite{HellerOW} in terms of a motivic tom Dieck splitting.
It would be interesting to express the right-hand side of the trace formula in terms of their computation.

\subsubsection{Cdh descent for the homotopy $K$-theory of tame stacks}

The question
of cdh descent for the homotopy $K$-theory of tame Deligne–Mumford stacks
was raised in \cite[\S1]{KO} (it was answered affirmatively for schemes in \cite{Cisinski}).
Our formalism of six operations yields a positive answer to this question for nice enough stacks. We only sketch the proof here, the details of which appear in \cite{HoyoisKH}.
Suppose that $G$ is a tame finite étale group scheme over $B$. For $S$ a quasi-projective $G$-scheme, let $\mathbb{K}^G(S)$ denote the nonconnective $K$-theory of perfect complexes on the quotient stack $[S/G]$.
Following Weibel \cite{Weibel}, define the \emph{homotopy $K$-theory} $\KH^G(S)$ by the formula
\begin{equation}\label{eqn:HIKT}
	\KH^G(S) = \colim_{n\in\Delta^\op} \mathbb{K}^G(S\times\A^n).
\end{equation}
It was shown in \cite[Corollary 5.6]{KO} that $\mathbb{K}^G$ satisfies Nisnevich excision.
Using the projective bundle formula from \cite{ThomasonK}, one can show that the restriction of $\KH^G$ to $\Sm_S^G$ is represented by an object $\mathrm{KGL}_{S}\in\SH^G(S)$. 
Furthermore, the tameness of $G$ allows us to generalize the results of \cite[\S4.2]{MV} to the equivariant setting. As a consequence,
we deduce that the underlying motivic $G$-space of $\mathrm{KGL}_S$ is equivalent to the group completion of a certain monoid structure on the motivic localization of $\coprod_{n\geq 0}\colim_{\scr E}\mathrm{Gr}_n(\scr E_{S})$, where $\scr E$ ranges over some cofiltered diagram of locally free $G$-modules on $B$, independent of $S$. It follows from this explicit description that $S\mapsto \mathrm{KGL}_S$ is a cocartesian section of $\SH^G(\ph)$. By Corollary~\ref{cor:descent}, we therefore deduce that $\KH^G$ satisfies cdh descent on quasi-projective $G$-schemes.

\subsection{Acknowledgements}

It is my pleasure to thank Joseph Ayoub, David Gepner, Jeremiah Heller, Adeel Khan, Amalendu Krishna, Marc Levine, and Paul Arne Østvær for many helpful and motivating discussions on equivariant motivic homotopy theory and the six operations. I am especially grateful to Elden Elmanto for carefully reading the paper and making numerous comments.
This work was partially supported by the National Science Foundation under grant DMS-1508096.

\subsection{Notation and terminology}
As a matter of terminology, we assume that group schemes are flat and finitely presented, and that locally free sheaves are of finite (but not necessarily constant) rank.
A \emph{vector bundle} over a scheme $X$ is a module scheme over $X$ that is isomorphic to $\V(\scr E)=\Spec(\Sym(\scr E))$ for some locally free sheaf of finite rank $\scr E$ on $X$. Thus, $\V$ is a contravariant equivalence between locally free sheaves and vector bundles. We write $\P(\scr E)=\Proj(\Sym(\scr E))$ for the projective bundle of lines in $\V(\scr E)$.
If $f\colon X\to S$ is a morphism of schemes, we denote by $\Omega_f$, $\Omega_{X/S}$, or even $\Omega_X$ the sheaf of relative differentials of $f$. If $i\colon Z\into X$ is an immersion, we denote by $\scr N_i$ its conormal sheaf.

Starting with the definition of the unstable equivariant motivic homotopy $\infty$-category in \S\ref{sec:unstable}, we will use the language of $\infty$-categories \cite{HTT}. We denote by $\scr S$ the $\infty$-category of $\infty$-groupoids and by $\Cat_\infty$ the $\infty$-category of (possibly large) $\infty$-categories. Categorical terminology must always be understood in the $\infty$-categorical context; for example, a presheaf or sheaf is $\scr S$-valued unless otherwise specified.
We denote by $\PSh(\scr C)$ the $\infty$-category of presheaves on a small $\infty$-category $\scr C$.

\endgroup

\section{Some equivariant geometry}
\label{sec:geometry}

Throughout this section, $B$ is an arbitrary base scheme and $G$ is a flat finitely presented group scheme over $B$.

\subsection{Invariant subschemes}

Let $X$ be a $G$-scheme and let $Y\subset X$ be a subscheme of $X$. We say that $Y$ is \emph{$G$-invariant} if there exists a morphism $G\times Y\to Y$ making the square
\begin{tikzmath}
	\diagram{G\times Y & Y \\ G\times X & X \\};
	\arrows (11-) edge[dashed] (-12) (11) edge[c->] (21) (12) edge[c->] (22) (21-) edge node[above]{$a$} (-22);
\end{tikzmath}
commute, where $a\colon G\times X\to X$ is the action. If such a morphism exists, it is clearly unique and is an action of $G$ on $Y$. 

Note that the action $a\colon G\times X\to X$ is isomorphic to the projection onto the second factor, and hence it is flat and finitely presented. In particular, it is open. Thus, if $U\subset X$ is any open subscheme, then $G\cdot U=a(G\times U)$ is an open subscheme of $X$, and it is clearly the smallest $G$-invariant open subscheme of $X$ containing $U$. Moreover, if $U$ is quasi-compact, then $G\times U$ and hence $G\cdot U$ are also quasi-compact.

\begin{lemma}\label{lem:invariant}
	Let $X$ be a $G$-scheme, $Z\subset X$ a closed subscheme, and $U=X\minus Z$ its open complement. Then $U$ is $G$-invariant if and only if $Z$ admits a $G$-invariant thickening, which may be chosen to be finitely presented if $X$ and $U$ are qcqs.
\end{lemma}

\begin{proof}
	Suppose that $Z$ admits a $G$-invariant thickening $Z'$, so that we have a cartesian square
	\begin{tikzmath}
		\diagram{G\times Z' & Z' \\ G\times X & X\rlap. \\};
		\arrows (11-) edge (-12) (11) edge[c->] (21) (12) edge[c->] (22) (21-) edge node[above]{$a$} (-22);
	\end{tikzmath}
	Then $a$ maps the open complement of $G\times Z'$ to the open complement of $Z'$, \ie, $U$ is $G$-invariant.
	
	For the converse, we consider more generally a groupoid scheme $d_0,d_1\colon X_1\rightrightarrows X_0$ where $d_0$ and $d_1$ are flat and finitely presented. Given a closed subscheme $Z\subset X_0$ with invariant open complement $U$, we will show that there exists an invariant thickening $Z'$ of $Z$ (in this context, a subscheme $Y\subset X_0$ is \emph{invariant} if $d_0^{-1}(Y)=d_1^{-1}(Y)$). Let $Z'$ be the schematic image of $d_0^{-1}(Z)$ by $d_1$.
	Note that there are (isomorphic) cartesian squares of the form:
	\[
	\begin{tikzpicture}
		\diagram{P & X_1 \\ d_0^{-1}(Z) & X_0\rlap, \\};
		\arrows (11-) edge node[above]{$f$} (-12) (11) edge (21) (21-) edge node[above]{$d_1$} (-22) (12) edge node[right]{$d_0$} (22);
	\end{tikzpicture}
	\quad
	\begin{tikzpicture}
		\diagram{P & X_1 \\ d_0^{-1}(Z) & X_0\rlap. \\};
		\arrows (11-) edge node[above]{$f$} (-12) (11) edge (21) (21-) edge node[above]{$d_1$} (-22) (12) edge node[right]{$d_1$} (22);
	\end{tikzpicture}
	\]
	Since flat base change preserves schematic images \cite[\href{http://stacks.math.columbia.edu/tag/081I}{Tag 081I}]{Stacks}, the first cartesian square tells us that the schematic image of $f$ is $d_0^{-1}(Z')$. The second cartesian square then tells us that $d_1^{-1}(Z')=d_0^{-1}(Z')$ as closed subschemes of $X_1$, \ie, that $Z'$ is invariant. It remains to prove that $Z'\cap U=\emptyset$. Since $U\into X$ is flat, $Z'\cap U$ is the schematic image of $d_1\colon d_0^{-1}(Z)\cap d_1^{-1}(U)\to U$, which is empty since $U$ is invariant.
	
	The last statement is proved using noetherian approximation, see \cite[\S8]{EGA4-3} and \cite[Appendix C]{TT}.
	Assume that $X_0$ is qcqs. Then there exists a cartesian morphism of groupoid schemes $X_\bullet\to Y_\bullet$ where $Y_\bullet$ is of finite type over $\Z$ and the maps $Y_1\rightrightarrows Y_0$ are flat and finitely presented \cite[\href{http://stacks.math.columbia.edu/tag/04AI}{Tag 04AI}]{Stacks}. If $U$ is quasi-compact, then $Z$ admits a finitely presented thickening \cite[Lemma 2.6.1 (c)]{TT}, so we can also assume that $Z$ is pulled back from a closed subscheme of $Y_0$ with an invariant open complement. Applying the result of the previous paragraph to $Y_\bullet$ shows that $Z$ admits a finitely presented invariant thickening.
\end{proof}

\begin{lemma}\label{lem:Gimmersion}
	Let $X$ be a $G$-scheme and $V\subset X$ a $G$-invariant subscheme.
	\begin{enumerate}
		\item The largest open subset of $X$ in which $V$ is closed is $G$-invariant.
		If moreover $V$ is quasi-compact, there exists a quasi-compact $G$-invariant open subscheme of $X$ in which $V$ is closed.
		\item 
		If $V\into X$ is quasi-compact, then the schematic closure of $V$ in $X$ is $G$-invariant.
		If moreover $X$ is qcqs and $V\into X$ is finitely presented, there exists a finitely presented $G$-invariant closed subscheme of $X$ in which $V$ is open.
	\end{enumerate}
\end{lemma}

\begin{proof}
	(1) Let $U\subset X$ be an open subset in which $V$ is closed. 
	It will suffice to show that $V$ is closed in $G\cdot U$. Since $V$ is $G$-invariant, $a^{-1}(V)=G\times V$, where $a\colon G\times X\to X$ is the action. Thus, $(G\cdot U)\minus V=a((G\times U)\minus a^{-1}(V))=a(G\times(U\minus V))$, which is open in $X$ since $a$ is open.
	
	(2) As in the proof of Lemma~\ref{lem:invariant}, we consider the more general situation of a groupoid scheme $d_0,d_1\colon X_1\rightrightarrows X_0$, with $d_0$ and $d_1$ flat and finitely presented and $V\subset X_0$ an invariant subscheme. The last part of the statement can then be proved by noetherian approximation.
	Let $\bar V$ be the schematic closure of $V$ in $X_0$. Consider the cartesian squares
	\[
	\begin{tikzpicture}
		\diagram{d_0^{-1}(V) & X_1 \\ V & X_0\rlap, \\};
		\arrows (11-) edge[c->] (-12) (11) edge node[left]{$d_0$} (21) (21-) edge[c->] (-22) (12) edge node[right]{$d_0$} (22);
	\end{tikzpicture}
	\quad
	\begin{tikzpicture}
		\diagram{d_0^{-1}(V) & X_1 \\ V & X\rlap, \\};
		\arrows (11-) edge[c->] (-12) (11) edge node[left]{$d_1$} (21) (21-) edge[c->] (-22) (12) edge node[right]{$d_1$} (22);
	\end{tikzpicture}
	\]
	and recall that flat base change preserves schematic images \cite[\href{http://stacks.math.columbia.edu/tag/081I}{Tag 081I}]{Stacks}. The first cartesian square tells us that the schematic closure of $d_0^{-1}(V)$ in $X_1$ is $d_0^{-1}(\bar V)$. The second cartesian square then tells us that $d_1^{-1}(\bar V)=d_0^{-1}(\bar V)$, \ie, that $\bar V$ is invariant.
\end{proof}

A \emph{$G$-immersion} is a $G$-morphism that is also an immersion. By Lemma~\ref{lem:Gimmersion} (1), a $G$-immersion can be factored canonically as a closed $G$-immersion followed by an open $G$-immersion.

\begin{lemma}\label{lem:locus}
	Let $f\colon X\to S$ be a $G$-morphism and $Y\subset X$ a subset. If $f$ is smooth (\resp{} étale) at the points of $Y$, then $f$ is smooth (\resp{} étale) on a $G$-invariant open neighborhood of $Y$, which may be chosen to be quasi-compact if $Y$ is quasi-compact.
\end{lemma}

\begin{proof}
	Let $U$ be an open subset of $X$ where $f$ is smooth (\resp{} étale). It will suffice to show that $f$ is smooth (\resp{} étale) on $G\cdot U$.
	Let $x\in G\cdot U$ be any point. Then $x=a(y)$ for some $y\in G\times U$,
	where $a\colon G\times X\to X$ is the action.
	Note that the square
	\begin{tikzmath}
		\diagram{G\times X & X \\ G\times S & S \\};
		\arrows (11-) edge node[above]{$a$} (-12) (21-) edge node[above]{$a$} (-22) (11) edge node[left]{$\id\times f$} (21) (12) edgenode[right]{$f$} (22);
	\end{tikzmath}
	is cartesian.
	Since $\id\times f$ is smooth (\resp{} étale) at $y$ and $a$ is flat, it follows from \cite[17.7.1 (ii)]{EGA4-3} that $f$ is smooth (\resp{} étale) at $x$.
\end{proof}

\begin{lemma}\label{lem:cartesian}
	Let
	\begin{tikzmath}
		\diagram{Z & X \\ & Y \\};
		\arrows (11-) edge[c->] node[above]{$s$} (-12) (11) edge[c->] node[below left]{$t$} (22) (12) edge node[right]{$p$} (22);
	\end{tikzmath}
	be a commutative triangle of $G$-schemes where $s$ and $t$ are closed immersions and $p$ is unramified. Then there exists a $G$-invariant open subscheme $U\subset X$ containing $Z$ such that the square
	\begin{tikzmath}
		\diagram{Z & U \\ Z & Y \\};
		\arrows (11-) edge[c->] node[above]{$s$} (-12)
		(11) edge[-,vshift=1pt] (21) edge[-,vshift=-1pt] (21)
		(21-) edge[c->] node[above]{$t$} (22) (12) edge node[right]{$p$} (22);
	\end{tikzmath}
	is cartesian. If $Z$ is quasi-compact, then $U$ may be chosen to be quasi-compact.
\end{lemma}

\begin{proof}
	Since $p$ is unramified, the $G$-immersion $(\id,s)\colon Z\into Z\times_YX$ is open.
	We can then take $U$ to be the complement of $Z\times_YX\minus Z$ in $X$, which is a $G$-invariant open subscheme of $X$ by Lemma~\ref{lem:invariant}.
	If $Z$ is quasi-compact, we can replace $U$ by $G\cdot V$ where $V\subset U$ is a quasi-compact open subset containing $Z$.
\end{proof}

\subsection{Quasi-affine and quasi-projective morphisms}

\begin{definition}\label{def:Gaffine}
	Let $f\colon X\to Y$ be a $G$-morphism.
	\begin{enumerate}
		\item $f$ is \emph{$G$-quasi-affine} (\resp{} \emph{$G$-affine}) if there exists a locally free $G$-module $\scr E$ on $Y$ and a quasi-compact $G$-immersion (\resp{} a closed $G$-immersion) $X\into \bb V(\scr E)$ over $Y$.
		\item $f$ is \emph{$G$-quasi-projective} (\resp{} \emph{$G$-projective}) if there exists a locally free $G$-module $\scr E$ on $Y$ and a quasi-compact $G$-immersion (\resp{} a closed $G$-immersion) $X\into \bb P(\scr E)$ over $Y$.
	\end{enumerate}
\end{definition}

\begin{lemma}\label{lem:Gquasiproj}
	Let $\mathbf C$ be any of the four classes of $G$-morphisms introduced in Definition~\ref{def:Gaffine}.
	\begin{enumerate}
		\item $\mathbf C$ is closed under base change.
		\item If $g$ is separated and $g\circ f$ belongs to $\mathbf C$, so does $f$.
	\end{enumerate}
\end{lemma}

\begin{proof}
	Easy.
\end{proof}

On the other hand, it is not true that the classes of morphisms from Definition~\ref{def:Gaffine} are stable under composition (even if $G$ is trivial!). 

\begin{definition}[{\cite[Definition\ 2.1]{Thomason}}]
	\label{def:resolution}
	Let $X$ be a $G$-scheme. We say that $X$ has the \emph{$G$-resolution property} if, for every finitely generated quasi-coherent $G$-module $\scr M$ on $X$, there exists a locally free $G$-module of finite rank $\scr E$ and an epimorphism $\scr E\onto\scr M$.
\end{definition}

\begin{example}
	\label{exa:resolution}
	Suppose that $B$ 
	is divisorial (\eg, affine, or noetherian, regular, and separated).
	Then $B$ has the $G$-resolution property in the following cases:
	\begin{enumerate}
		\item $G$ is finite locally free;
		\item $G$ is isotrivial; 
		\item $G$ is reductive with isotrivial radical and coradical (\eg, $G$ is semisimple). 
	\end{enumerate}
	This is proved in \cite{Thomason} under the assumption that $B$ is noetherian. An easy noetherian approximation argument implies the general case.
\end{example}

\begin{remark}
	\label{rmk:isotrivial}
	 Any finitely presented group of multiplicative type over $B$ is isotrivial locally in the Nisnevich topology.
	Indeed, if $B$ is a henselian local scheme, any étale cover of $B$ is refined by a finite étale cover, by \cite[Théorème 18.5.11 (c)]{EGA4-4}.
\end{remark}

\begin{lemma}\label{lem:thomason}
	Let $X$ be a qcqs $G$-scheme. Every quasi-coherent $G$-module on $X$ is the colimit of its finitely generated quasi-coherent $G$-submodules.
\end{lemma}

\begin{proof}
	Let $X_\bullet$ be the groupoid scheme defined by the action of $G$ on $X$.
	Since $X$ is qcqs and $G$ is flat and finitely presented, there exists an affine cartesian morphism of groupoid schemes $p\colon X_\bullet\to Y_\bullet$, where $Y_\bullet$ is of finite type over $\Z$ and the maps $Y_1\rightrightarrows Y_0$ are flat and finitely presented.
	Let $\scr M$ be a quasi-coherent $G$-module on $X$. By \cite[\href{http://stacks.math.columbia.edu/tag/09VH}{Tag 09VH}]{Stacks}, $p_*(\scr M)$ has a structure of quasi-coherent module on the groupoid $Y_\bullet$ such that the counit map $p^*p_*(\scr M)\to \scr M$ is a map of $G$-modules.
	Since $p$ is affine, this counit map is an epimorphism. Thus, it suffices to show that $p_*(\scr M)$ is the union of its finitely generated quasi-coherent submodules on the groupoid $Y_\bullet$. This is true by \cite[\href{http://stacks.math.columbia.edu/tag/07TU}{Tag 07TU}]{Stacks}.
\end{proof}

\begin{lemma}\label{lem:QPresolution}
Let $f\colon X\to Y$ be a $G$-quasi-projective $G$-morphism.
If $Y$ is qcqs and has the $G$-resolution property, then $X$ has the $G$-resolution property.
\end{lemma}

\begin{proof}
	One can repeat the proof of \cite[Lemma 2.6]{Thomason}, using Lemma~\ref{lem:thomason}. Details are omitted.
\end{proof}

\begin{lemma}\label{lem:Gquasiaffine}
	Let $f\colon X\to Y$ be a quasi-affine (\resp{} affine) $G$-morphism of finite type. If $Y$ is qcqs and has the $G$-resolution property, then $f$ is $G$-quasi-affine (\resp{} $G$-affine).
\end{lemma}

\begin{proof}
	We repeat the argument from \cite[Theorem 3.6]{Thomason}. Since $f$ is quasi-affine, $f_*(\scr O_X)$ is a quasi-coherent $G$-algebra and we have an open $G$-immersion $X\into \Spec(f_*(\scr O_X))$ over $Y$ (an isomorphism if $f$ is affine). Let $\scr M$ be a finitely generated quasi-coherent $G$-submodule of $f_*(\scr O_X)$ that generates $f_*(\scr O_X)$ as an algebra (such an $\scr M$ exists by Lemma~\ref{lem:thomason}). By the $G$-resolution property, there exists a locally free $G$-module $\scr E$ and an epimorphism $\scr E\onto \scr M$. Then the epimorphism of $G$-algebras $\Sym(\scr E)\onto f_*(\scr O_X)$ induces a closed $G$-immersion $\Spec(f_*(\scr O_X))\into \bb V(\scr E)$ over $Y$. Thus, $f$ is $G$-(quasi)-affine.
\end{proof}

\begin{lemma}\label{lem:composition}
	Let $\mathbf C$ be any of the four classes of $G$-morphisms introduced in Definition~\ref{def:Gaffine} and let $f\colon X\to Y$ and $g\colon Y\to Z$ belong to $\mathbf C$. If $Z$ is qcqs and has the $G$-resolution property, then $g\circ f$ belongs to $\mathbf C$.
\end{lemma}

\begin{proof}
	The case of $G$-(quasi)-affine morphisms follows from Lemma~\ref{lem:Gquasiaffine}. Let us treat the case of $G$-(quasi)-projective morphisms. By definition, there exist locally free $G$-modules $\scr F$ on $Y$ and $\scr G$ on $Z$, and quasi-compact (closed) $G$-immersions $X\into \P(\scr F)$ over $Y$ and $Y\into\P(\scr G)$ over $Z$. Let us abbreviate $\scr F\tens \scr O_Y(1)^{\tens n}$ to $\scr F(n)$.
	By \cite[Proposition 4.6.8]{EGA2} (amended by \cite[1.7.15]{EGA4-1}), there exists $n$ such that the counit $g^*g_*(\scr F(n))\to \scr F(n)$ is an epimorphism. Since $\scr F(n)$ is finitely generated and $g_*(\scr F(n))$ is the union of its finitely generated quasi-coherent $G$-submodules (Lemma~\ref{lem:thomason}), there exists a finitely generated quasi-coherent $G$-submodule $\scr M$ of $g_*(\scr F(n))$ such that $g^*(\scr M)\to \scr F(n)$ is an epimorphism. Then, by the $G$-resolution property, there exists a locally free $G$-module $\scr E$ and an epimorphism $\scr E\onto \scr M$. Thus, we obtain a closed $G$-immersion $\P(\scr F)\into\P(g^*(\scr E))$ over $Y$, whence a quasi-compact (closed) $G$-immersion $X\into \P(\scr E)\times\P(\scr G)$ over $Z$. We conclude using the Segre embedding $\P(\scr E)\times\P(\scr G)\into \P(\scr E\tens\scr G)$, which is a closed $G$-immersion.
\end{proof}

\subsection{Linearly reductive groups}

If $S$ is a $G$-scheme, we denote by $\QCoh^G(S)$ the category of quasi-coherent $G$-modules on $S$.

\begin{definition}[{\cite[Definition 2.2]{AOV}}]
	\label{def:linred}
	A group scheme $G$ over $B$ is called \emph{linearly reductive} if the $G$-fixed-point functor $(\ph)^G\colon\QCoh^G(B)\to\QCoh(B)$ is exact.
\end{definition}

\begin{lemma}[{\cite[Proposition 2.4]{AOV}}]
	\label{lem:localred}
	\leavevmode
	\begin{enumerate}
		\item If $\{f_i\colon U_i\to B\}$ is an fpqc cover and $G_{U_i}$ is linearly reductive for all $i$, then $G$ is linearly reductive.
		\item If $G$ is linearly reductive and $f\colon B'\to B$ is any morphism, then $G_{B'}$ is linearly reductive.
	\end{enumerate}
\end{lemma}

\begin{proof}
Assertion (1) is clear since the family of functors $\{f_i^*\}$ detects exactness and each $f_i^*$ commutes with $(\ph)^G$.
Given (1), it suffices to prove (2) when $f$ is an open immersion and when $f$ is affine. In the latter case the result is clear since $f_*$ is exact and conservative. Suppose that $f$ is an open immersion. Then both $f_*$ and $f^*$ commute with $(\ph)^G$, $f^*$ is exact, and $f_*$ is fully faithful. Let $(\scr M_\alpha)$ be a finite diagram in $\QCoh^G(B')$. Then the computation
\[
	\colim_\alpha \scr M_\alpha^G\simeq\colim_\alpha f^*f_*(\scr M_\alpha^G)\simeq f^*\colim_\alpha f_*(\scr M_\alpha^G)
	\simeq(f^*\colim_\alpha f_*\scr M_\alpha)^G\simeq (\colim_\alpha f^*f_*\scr M_\alpha)^G\simeq (\colim_\alpha\scr M_\alpha)^G
\]
shows that $(\ph)^G\colon\QCoh^G(B')\to\QCoh(B')$ is right exact, \ie, that $G_{B'}$ is linearly reductive.
\end{proof}

\begin{samepage}
\begin{example}\label{exa:linred}
	\leavevmode
	\begin{enumerate}
		\item A finite locally free group of order invertible on $B$ is linearly reductive \cite[Theorem 2.16]{AOV}.
		\item A group of multiplicative type is linearly reductive (combine \cite[Exposé I, Théorème 5.3.3]{SGA3-1} and Lemma~\ref{lem:localred} (1)).
		\item If $B$ has characteristic zero, any reductive group over $B$ is linearly reductive (for fields, see \cite{Nagata}; the general case follows by Lemma~\ref{lem:localred} since any reductive group is étale-locally split). 
	\end{enumerate}
\end{example}
\end{samepage}

Linear reductivity will be used exclusively through the following lemma:

\begin{lemma}\label{lem:linred}
	Suppose that $B$ is affine and that $G$ is linearly reductive, and let $p\colon S\to B$ be an affine $G$-morphism. If $\scr M$ is a locally free $G$-module of finite rank on $S$, then $\scr M$ is projective in $\QCoh^G(S)$.
\end{lemma}

\begin{proof}
	We need to show that $\Hom(\scr M,\ph)\colon \QCoh^G(S)\to\scr A\mathrm b$ is an exact functor. We have
	\[\Hom(\scr M,\ph)\simeq \Gamma_B\circ (\ph)^G\circ p_*\circ \scr H\mathrm{om}(\scr M,\ph),\]
	where:
	\begin{itemize}
		\item $\scr H\mathrm{om}\colon\QCoh^G(S)\times\QCoh^G(S)\to\QCoh^G(S)$ is the internal hom object in $\QCoh^G(S)$;
		\item $p_*\colon \QCoh^G(S)\to\QCoh^G(B)$ is the pushforward along $p$;
		\item $(\ph)^G\colon\QCoh^G(B)\to\QCoh(B)$ is the $G$-fixed point functor;
		\item $\Gamma_B\colon\QCoh(B)\to \scr A\mathrm b$ is the global section functor.
	\end{itemize}
	Now, $\scr H\mathrm{om}(\scr M,\ph)$ is exact because $\scr M$ is locally free of finite rank, $p_*$ is exact because $p$ is affine, $(\ph)^G$ is exact because $G$ is linearly reductive, and $\Gamma_B$ is exact because $B$ is affine. Thus, $\Hom(\scr M,\ph)$ is exact.
\end{proof}

\subsection{Affine resolutions}
\label{sub:jouanolou}

\begin{definition}
	Let $X$ be a $G$-scheme. A \emph{$G$-affine bundle} is a $G$-morphism $Y\to X$ that is a torsor under a $G$-vector bundle $V$ over $X$, such that the action $V\times_XY\to Y$ is $G$-equivariant.
\end{definition}

If $\Aff_r$ is the group scheme of affine automorphisms of $\A^r$ and $\mathbf B\Aff=\coprod_{r\geq 0}\mathbf B\Aff_r$, then $G$-affine bundles over $X$ are classified by morphisms of stacks
\[[X/G]\to \mathbf B\Aff.\]
Using the obvious embeddings $\Aff_r\into \GL_{r+1}$, we see that the category of $G$-affine bundles over a $G$-scheme $X$ is equivalent to the full subcategory of $\QCoh^G(X)_{/\scr O_X}$ spanned by the epimorphisms $\phi\colon \scr E\onto\scr O_X$, where $\scr E$ is locally free of finite rank. The $G$-affine bundle $Y$ corresponding to $\phi$ is the $X$-scheme that classifies splittings of $\phi$. In other words, $Y$ is the preimage of the $1$-section by the epimorphism $\V(\scr E^\vee)\onto \A^1_X$, or equivalently the complement of $\P(\ker(\phi)^\vee)$ in $\P(\scr E^\vee)$, and it is a torsor under $\V(\ker(\phi)^\vee)$.

We now recall Jouanolou's trick, with an equivariant twist. Let $S$ be a $G$-scheme, let $\scr E$ be a locally free $G$-module of rank $\geq 1$ over $S$, and let $p\colon \P(\scr E)\to S$ be the associated projective bundle. We then have a short exact sequence
\begin{equation}\label{eqn:tracezero}
	0\to \scr O_S\to \scr E\tens\scr E^\vee\to \scr Q\to 0
\end{equation}
in $\QCoh^G(S)$, where the first map is the coevaluation. The naturality of the coevaluation implies that the quotient $\scr Q$ is again locally free. We therefore have a $G$-invariant hyperplane
\[\P(\scr Q)\into \P(\scr E\tens\scr E^\vee)\]
whose open complement is an affine $G$-scheme over $S$. Let $\tilde\P(\scr E)$ be the preimage of $\P(\scr E\tens\scr E^\vee)\minus\P(\scr Q)$ by the Segre embedding
\[\segre\colon\P(\scr E)\times_S\P(\scr E^\vee)\into \P(\scr E\tens\scr E^\vee),\]
and let $\pi\colon \tilde\P(\scr E)\to \P(\scr E)$ be the restriction of the projection onto the first factor. By construction, $\tilde \P(\scr E)$ is affine over $S$; it is even $G$-affine over $S$ if $S$ has the $G$-resolution property (Lemma~\ref{lem:Gquasiaffine}) or if the exact sequence~\eqref{eqn:tracezero} splits equivariantly, for example if $B$ and $S$ are affine and $G$ is linearly reductive (Lemma~\ref{lem:linred}).

We claim that $\pi$ is $G$-affine bundle.
Consider the canonical epimorphism
\[\phi\colon p^*(\scr E)(-1)\onto\scr O_{\P(\scr E)}\]
in $\QCoh^G(\P(\scr E))$, whose kernel is the sheaf of differentials $\Omega_{\P(\scr E)/S}$. Over a $T$-point $\alpha\colon\scr E_T\onto\scr L$ of $\P(\scr E)$, $\phi$ is given by the composition
\[\scr E_T\tens \scr L^\vee\xrightarrow{\alpha^\vee} \scr E_T\tens\scr E^\vee_T\to\scr O_T,\]
where the second map is the evaluation.
The $G$-affine bundle corresponding to $\phi$ is then the complement of $\P(\Omega_{\P(\scr E)/S}^\vee)$ in $\P(p^*(\scr E)^\vee)$.
Under the isomorphism $\P(p^*(\scr E)^\vee)\simeq \P(\scr E)\times_S\P(\scr E^\vee)$, it is clear that this complement is exactly $\tilde \P(\scr E)$. This shows that $\pi$ is a $G$-equivariant torsor under the cotangent bundle of $\P(\scr E)$.

\begin{proposition}[Jouanolou's trick]
	\label{prop:jouanolou}
	Let $S$ be a qcqs $G$-scheme with the $G$-resolution property. Then, for every $G$-quasi-projective $S$-scheme $X$, there exists a $G$-affine bundle $\tilde X\to X$ where $\tilde X$ is affine over $S$.
\end{proposition}

\begin{proof}
	If $X=\P(\scr E)$, we have just proved this. It is then clear how to obtain such a bundle if $X$ is $G$-projective.
	In general, it suffices to construct an affine $G$-morphism $X\to P$ where $P$ is $G$-projective: we then let $\tilde X$ be the pullback of $\tilde P$. Choose a locally free $G$-module $\scr F$ and a quasi-compact $G$-immersion $X\into\P(\scr F)$ over $S$. Let $\bar X$ be a $G$-invariant closed subscheme of $\P(\scr F)$ of which $X$ is an open subscheme (Lemma~\ref{lem:Gimmersion} (2)), let $Z=\Spec(\scr O_{\bar X}/\scr I)$ be a $G$-invariant closed subscheme of $\bar X$ with open complement $X$ (Lemma~\ref{lem:invariant}), and let $P$ be the blowup of $\bar X$ at $Z$.
	Since $\bar X$ has the $G$-resolution property, there exists a locally free $G$-module $\scr G$ and an epimorphism $\scr G\onto \scr I$, whence a closed $G$-immersion $P\into \P(\scr G)$ over $\bar X$, so that $P$ is $G$-projective (Lemma~\ref{lem:composition}). Finally, the open $G$-immersion $X\into P$ factors through a closed $G$-immersion $X\into \V(\scr J^\vee)$ over $P$, where $\scr J$ is the ideal of the exceptional divisor in $P$, and in particular it is affine.
\end{proof}

\subsection{Lifting locally free sheaves}

Given a $G$-scheme $X$ and a $G$-invariant closed subscheme $Z\subset X$, we shall say that a $G$-morphism $f\colon X'\to X$ is a \emph{$G$-equivariant étale neighborhood} of $Z$ if it is locally finitely presented, if the induced map $Z\times_X X'\to Z$ is an isomorphism, and if $f$ is étale at all points lying over $Z$.
In that case, by Lemma~\ref{lem:locus}, there exists a $G$-invariant open subscheme $V\subset X'$ such that $f|V$ is étale and is an isomorphism over $Z$. However, it will be important to allow $f$ itself not to be étale, since we often need $X'$ to be affine but we cannot guarantee that $V$ is affine.

The following theorem is an equivariant generalization of a theorem of Arabia \cite[Théorème 1.2.3]{Arabia}.

\begin{theorem}\label{thm:lift-locally-free}
	Suppose that $B$ is affine and that $G$ is linearly reductive. Let $s\colon Z\into X$ be a closed $G$-immersion between affine $G$-schemes and let $\scr N$ be a locally free $G$-module on $Z$. If $X$ has the $G$-resolution property, there exists an affine $G$-equivariant étale neighborhood $X'\to X$ of $Z$ and a locally free $G$-module $\scr M$ on $X'$ lifting $\scr N$. 
\end{theorem}

\begin{proof}
 Since $s_*(\scr N)$ is finitely generated and $X$ has the $G$-resolution property, there exists a locally free $G$-module $\scr E$ on $X$ and an epimorphism $\scr E\onto s_*(\scr N)$, whence an epimorphism $s^*(\scr E)\onto \scr N$. Since $\scr N$ is projective in $\QCoh^G(Z)$ by Lemma~\ref{lem:linred}, there exists a $G$-equivariant idempotent endomorphism $\phi$ of $s^*(\scr E)$ whose image is isomorphic to $\scr N$. 
	Since $\scr E$ is projective in $\QCoh^G(X)$ and the unit map $\scr E\to s_*s^*(\scr E)$ is an epimorphism, there exists a $G$-equivariant endomorphism $\psi$ of $\scr E$ such that $s^*(\psi)=\phi$.
	The idea is now to consider the universal deformation of $\psi$, and the locus where it is idempotent will be the desired étale neighborhood.
	
	Let $\scr I\subset\scr O_X$ be the ideal of $s$.
	By Lemma~\ref{lem:thomason}, $\scr I$ is the union of its finitely generated quasi-coherent $G$-submodules. Since $\phi$ is idempotent, $\psi^2-\psi$ has image in $\scr I\scr E$, hence in $\scr I'\scr E$ for some finitely generated $\scr I'\subset\scr I$. If $Z'\subset X$ is the closed subscheme defined by $\scr I'$, it follows that the idempotent $\phi$ and hence $\scr N$ lift to $Z'$. Replacing $Z$ by $Z'$, we can therefore assume that $\scr I$ is finitely generated. In that case, by the $G$-resolution property, there exists a locally free $G$-module $\scr F$ on $X$ and an epimorphism $\scr F\onto \scr I$.
	Let $\pi$ denote the composition $\scr F\onto \scr I\subset\scr O_X$.
	 By projectivity of $\scr E$, we can lift $\psi^2-\psi\colon \scr E\to\scr I\scr E$ to a morphism $\alpha\colon \scr E\to \scr F\otimes\scr E$, so that $\psi^2-\psi=\pi\alpha$.
	 Let $p\colon V\to X$ be the $G$-vector bundle whose sheaf of sections is $\scr H\mathrm{om}(\scr E,\scr F\otimes\scr E)$, let $\beta\colon p^*(\scr E)\to p^*(\scr F\otimes\scr E)$ be the tautological morphism, and let
	\[
	R= \alpha + (2\psi-\id)\beta + \pi\beta^2 \colon p^*(\scr E) \to p^*(\scr F\otimes\scr E).
	\]
	This should be understood as follows: locally where $\scr F$ is free, $\pi$ is a collection of generators $\pi_1,\dotsc,\pi_n$ of $\scr I$, $\alpha$ is collection of endomorphisms $\alpha_1,\dotsc,\alpha_n$ of $\scr E$ such that $\psi^2-\psi=\sum_i\pi_i\alpha_i$, and $\beta$ is the universal family of $n$ endomorphisms of $\scr E$. 
	Moreover, $\pi R=(\psi+\pi\beta)^2-(\psi+\pi\beta)+[\psi,\pi\beta]$, so the equation $\pi R=0$ expresses the idempotency of $\psi+\pi\beta$ wherever $\psi$ and $\pi\beta$ commute.

	Let $i\colon X'\into V$ be the locus where $R$ and $[\psi,\pi\beta]$ both vanish, so that $\psi+\pi\beta$ is an idempotent endomorphism of $i^*p^*(\scr E)$. Let $\scr M$ be its image. By construction, $\scr M$ and $\scr N$ are isomorphic over $Z\times_XX'$. 
	Since $\phi$ is idempotent, $2\phi-\id$ is an automorphism of $s^*(\scr E)$, and since $R\equiv \alpha+(2\psi-\id)\beta$ modulo $\scr I$, it follows that $p\circ i\colon X'\to X$ induces an isomorphism $Z\times_XX'\simeq Z$.
	It remains to check that $p\circ i$ is étale at all points lying over $Z$. Since this question is local on $X$ and does not involve $G$, we can assume that $\scr E$ and $\scr F$ are free, so that $\scr F$ specifies global generators $\pi_1,\dotsc,\pi_n$ of $\scr I$. By induction on $n$, we can further assume that $n=1$, \ie, that $\scr I$ is a principal ideal generated by $\pi$. 
	Let $X''\subset V$ be the locus where $R$ vanishes, so that $X'\subset X''$ is the locus where $[\psi,\pi\beta]$ vanishes. Note that $X''\to X$ is also an isomorphism over $Z$. In the proof of \cite[Théorème 1.2.3]{Arabia}, Arabia shows that $X''\to X$ is étale over $Z$ and that $[\psi,\pi\beta]$ vanishes in an open neighborhood of $Z$ in $X''$. This implies that $X'\to X$ is also étale over $Z$, as desired.
\end{proof}

\subsection{Lifting smooth quasi-sections}

\begin{theorem}\label{thm:keyLift}
	Suppose that $B$ is affine and that $G$ is linearly reductive.
	Let
	\begin{tikzmath}
		\diagram{X_Z & X \\ Z & S \\};
		\arrows 
		(11-) edge[c->] node[above]{$t$} (-12)
		(11) edge (21) (12) edge node[right]{$p$} (22) (21-) edge[c->] node[above]{$s$} (-22)
		;
	\end{tikzmath}
	be a cartesian square of $G$-schemes where $X$ is affine and $s$ is a closed immersion. Let $V\subset X_Z$ be a finitely presented $G$-invariant closed subscheme of $X_Z$ that is smooth (\resp{} étale) over $Z$.
	Suppose that $p$ is smooth at each point of $V$ and that $X$ has the $G$-resolution property. Then there exists an affine $G$-equivariant étale neighborhood $X'\to X$ of $V$ and a finitely presented $G$-invariant closed subscheme $\hat V\subset X'$ lifting $V$ such that $\hat V\to S$ is smooth (\resp{} étale) at each point of $V$.
\end{theorem}

\begin{proof}
	 Let $i\colon V\into X_Z$ be the inclusion and $\scr I\subset\scr O_{X_Z}$ its ideal.
	 Since $p$ is smooth at $V$, the conormal sheaf $\scr N_i$ is locally free.
	 Replacing $X$ by an affine $G$-equivariant étale neighborhood of $V$, we can assume by Theorem~\ref{thm:lift-locally-free} that there exists a locally free $G$-module $\scr E$ on $X$ and an isomorphism $i^*t^*(\scr E)\simeq\scr N_i$. Let $\phi\colon\scr E\onto t_*i_*(\scr N_i)$ be the adjoint morphism.
	 By Lemma~\ref{lem:linred}, $\scr E$ is projective in $\QCoh^G(X)$. Thus, we can find successive lifts in the diagram
 	\begin{tikzmath}
 		\diagram{\scr E & & \scr O_X \\ t_*i_*(\scr N_i) & t_*(\scr I) & t_*(\scr O_{X_Z})\rlap. \\};
 		\arrows (11-) edge[dashed] node[above]{$\chi$} (-13) (11) edge[->>] node[left]{$\phi$} (21)
		(21-) edge[<<-] (-22) (11) edge[dashed] node[above right=-1pt]{$\psi$} (22)
 		(13) edge[->>] (23) (22-) edge[c->] (-23);
 	\end{tikzmath}
	By Nakayama's lemma, since $\scr I$ is finitely generated, $\psi$ is surjective in a neighborhood of $V$. 
	Let $\scr J\subset\scr O_X$ be the image of $\chi$ and let $\hat V$ be the $G$-invariant closed subscheme of $X$ defined by $\scr J$. By construction, the image of $\scr J$ in $t_*(\scr O_{X_Z})$ equals $t_*(\scr I)$ in a neighborhood of $V$, 
	so that $t^{-1}(\hat V)=V\amalg K$. Using the projectivity of $\scr O_X$ in $\QCoh^G(X)$, we can lift the morphism $(1,0)\colon \scr O_X\to \scr O_V\times\scr O_K$ to an endomorphism of $\scr O_X$, which gives a $G$-invariant function $f$ on $X$ such that the affine open subscheme $X_f$ contains $V$ and is disjoint from $K$. Replacing $X$ by $X_f$, we can therefore assume that $t^{-1}(\hat V)=V$.
	
	It remains to
	show that $\hat V\to S$ is smooth (\resp{} étale) at each point $v\in V$. Let $z\in Z$ be the image of $v$. Since $\hat V_z\simeq V_z$ is smooth (\resp{} étale) over $z$ by assumption, it suffices to show that $\hat V\to S$ is flat at $v$ \cite[17.5.1 (b)]{EGA4-4} (\resp{} \cite[17.6.1 (c$'$)]{EGA4-4}). Let $c$ be the codimension of $V$ in $X_Z$ at $v$, or equivalently the rank of $\scr E$ at $v$. Choose a basis $f_1$, \dots, $f_c$ of the $\kappa(v)$-vector space $\scr I(v)$. For each $1\leq r\leq c$, choose a local section $\omega_r$ of $\scr E$ such that $\psi(\omega_r)(v)=f_r$, and let $g_r=\chi(\omega_r)$.
	Since the dimension of $\scr J(v)$ is at most the rank of $\scr E$, the epimorphism $\scr J(v)\onto\scr I(v)$ is an isomorphism. 
	 Hence, the functions $g_r$ generate $\scr J(v)$ as a $\kappa(v)$-vector space. By Nakayama's lemma, there exists an open neighborhood of $v$ in $X$ where $g_1$, \dots, $g_c$ generate $\scr J$ as an $\scr O_X$-module. On the other hand, since $\scr I$ is a regular ideal \cite[17.12.1]{EGA4-4}, the images of $g_1$, \dots, $g_c$ in the local ring $\scr O_{X_z,v}$ form a regular sequence. It follows from \cite[11.3.8 (c)]{EGA4-3} that $\hat V\to S$ is flat at $v$.
\end{proof}

\begin{corollary}\label{cor:keyLift}
	Suppose that $B$ is affine and that $G$ is linearly reductive.
	Let $S$ be an affine $G$-scheme with the $G$-resolution property, $s\colon Z\into S$ a closed $G$-immersion, and $X$ a smooth (\resp{} étale) affine $G$-scheme over $Z$. Then there exists a finitely presented affine $G$-scheme $\hat X$ over $S$ lifting $X$ such that $\hat X\to S$ is smooth (\resp{} étale) at each point of $X$.
\end{corollary}

\begin{proof}
	By Lemma~\ref{lem:Gquasiaffine}, there exists a locally free $G$-module $\scr E$ over $S$ and a closed $G$-immersion $X\into \bb V(\scr E)$ over $S$.
	In particular, we obtain the following diagram:
	\begin{tikzmath}
		\diagram{X & \bb V(s^*(\scr E)) & \bb V(\scr E) \\ & Z & S\rlap. \\};
		\arrows 
		(11-) edge[c->] (-12) (12-) edge[c->] (-13)
		(12) edge (22) (13) edge (23) (22-) edge[c->] node[above]{$s$} (-23)
		(11) edge (22)
		;
	\end{tikzmath}
	Now apply Theorem~\ref{thm:keyLift}.
\end{proof}

\begin{corollary}\label{cor:retraction}
	Suppose that $B$ is affine and that $G$ is linearly reductive.
	Let
	\begin{tikzmath}
		\diagram{Z & X \\ & S \\};
		\arrows (11-) edge[c->] node[above]{$s$} (-12) (12) edge node[right]{$p$} (22) (11) edge node[below left]{$q$} (22);
	\end{tikzmath}
	be a commutative triangle of $G$-schemes where $s$ is a closed $G$-immersion and $q$ is smooth and separated.
	Suppose that $X\times_SZ$ is affine and has the $G$-resolution property.
	Then there exists an affine $G$-equivariant étale neighborhood $X'\to X$ of $Z$ such that $Z\into X'$ admits a $G$-retraction.
\end{corollary}

\begin{proof}
	Consider the following diagram:
		\begin{tikzmath}
			\diagram{Z & Z\times_S Z & X\times_S Z \\ & Z & X\rlap. \\};
			\arrows (11-) edge[c->] node[above]{$\delta$} (-12) (11) edge[-,vshift=1pt] (22) edge[-,vshift=-1pt] (22)
			(12) edge node[right]{$\pi_1$} (22) (12-) edge[c->] node[above]{$s\times\id$} (-13) (13) edge node[right]{$\pi_1$} (23)
			(22-) edge[c->] node[above]{$s$} (-23)
			;
		\end{tikzmath}
		The vertical projections are smooth since $q$ is smooth, and $\delta$ is a closed immersion since $q$ is separated. Applying Theorem~\ref{thm:keyLift} to this diagram, we obtain a $G$-scheme $X'$ with the desired properties.
\end{proof}

\subsection{Linearizations}

\begin{proposition}\label{prop:linearization}
	Suppose that $B$ is affine and that $G$ is linearly reductive. Let $S$ be an affine $G$-scheme and let $p\colon X\to S$ be a finitely presented affine $G$-morphism with a quasi-regular $G$-section $s\colon S\into X$. Then there exists
	a $G$-morphism $h\colon X\to \bb V(\scr N_s)$ over $S$ such that:
	\begin{enumerate}
		\item the triangle
	\begin{tikzmath}
		\diagram{S & X \\ & \bb V(\scr N_s) \\};
		\arrows (11-) edge[c->] node[above]{$s$} (-12) (11) edge[c->] node[below left]{$z$} (22) (12) edge node[right]{$h$} (22);
	\end{tikzmath}
	commutes, where $z$ is the zero section;
	\item $h$ is étale at each point of $s(S)$.
	\end{enumerate}
\end{proposition}

\begin{proof} 
	Let $\scr I\subset\scr O_X$ be the vanishing ideal of $s$.
	Since $s$ is quasi-regular, $\scr N_s$ is locally free. By Lemma~\ref{lem:linred}, $\scr N_s$ is projective in $\QCoh^G(S)$. Thus, the epimorphism
	\[p_*(\scr I)\onto p_*(\scr I/\scr I^2)\simeq\scr N_s\]
	admits a $G$-equivariant section $\phi\colon \scr N_s\into p_*(\scr I)$.
	Let $\sigma\colon p_*(\scr O_X)\to\scr O_S$ be the morphism of $G$-algebras corresponding to the section $s$. Then $\ker(\sigma)=p_*(\scr I)$, so that the triangle of $G$-modules
	\begin{tikzmath}
		\diagram{\scr O_S & p_*(\scr O_X) \\ & \scr N_s \\};
		\arrows (11-) edge[<-] node[above]{$\sigma $} (-12) (11) edge[<-] node[below left]{$0$} (22) (12) edge[<-] node[right]{$\phi$} (22);
	\end{tikzmath}
	commutes. Let $\psi\colon\Sym(\scr N_s)\to p_*(\scr O_X)$ be the morphism of $G$-algebras induced by $\phi$, and let $h\colon X\to\bb V(\scr N_s)$ be the corresponding $G$-morphism over $S$. It is then clear that the triangle of $G$-algebras
	\begin{tikzmath}
		\diagram{\scr O_S & p_*(\scr O_X) \\ & \Sym(\scr N_s) \\};
		\arrows (11-) edge[<-] node[above]{$\sigma $} (-12) (11) edge[<-] (22) (12) edge[<-] node[right]{$\psi$} (22);
	\end{tikzmath}
	commutes, which proves (1).
	
	Let us check (2). By \cite[17.12.1 (c$'$)]{EGA4-4}, $p$ is smooth at each point of $s(S)$. By \cite[17.11.2 (c$'$)]{EGA4-4}, it therefore suffices to show that the map
	\[s^*(dh)\colon z^*(\Omega_{\V(\scr N_s)/S})\to s^*(\Omega_{X/S})\]
	is an isomorphism. By definition of $h$, this map is the composition
	\[\scr N_s\stackrel\phi\to p_*(\scr I)\to p_*(\scr I/\scr I^2)\simeq s^*(\Omega_{X/S}),\]
	and it is an isomorphism by choice of $\phi$. 
\end{proof}

\subsection{Tame group schemes}

\begin{definition}\label{def:tame}
	A flat finitely presented group scheme $G$ over $B$ is called \emph{tame} if
	the following conditions hold:
	\begin{itemize}
		\item $B$ admits a Nisnevich covering by schemes having the $G$-resolution property (Definition~\ref{def:resolution});
		\item $G$ is linearly reductive (Definition~\ref{def:linred}).
	\end{itemize}
\end{definition}

\begin{example}
	Let $G$ be a group scheme over $B$.
	Then $G$ is tame in the following cases (see Example~\ref{exa:resolution}, Remark~\ref{rmk:isotrivial}, and Example~\ref{exa:linred}):
	\begin{enumerate}
		\item $B$ is arbitrary and $G$ is finite locally free of order invertible on $B$.
		\item $B$ is arbitrary and $G$ is of multiplicative type (\eg, $G=\G_m$ or $G=\mu_n$).
		\item $B$ has characteristic zero and $G$ is reductive (\eg, $G=\GL_n$).
	\end{enumerate}
\end{example}

Note that our definition of tameness does not imply linearity or even affineness (for example, an elliptic curve is tame). However, we will soon restrict our attention to $G$-quasi-projective schemes, and these are only interesting if $G$ can act nontrivially on vector bundles.

\section{Unstable equivariant motivic homotopy theory}
\label{sec:unstable}

As a first step towards the formalism of six operations, we define in this section the unstable equivariant motivic homotopy $\infty$-category $\H^G(S)$ associated with a $G$-scheme $S$. 

\subsection{Preliminaries}
\label{sub:conventions}

For the remainder of this paper, we fix a qcqs base scheme $B$ and a tame group scheme $G$ over $B$. 
From now on we will only consider $G$-schemes $S$ that are finitely presented over $B$ and are Nisnevich-locally $G$-quasi-projective, \ie, for which there exists a Nisnevich cover $\{U_i\to B\}$ such that $S_{U_i}$ is $G$-quasi-projective over $U_i$. We denote by $\Sch_B^G$ the category of such $G$-schemes; this is the category on which the six operations will eventually be defined, although we will make one additional minor simplifying assumption on $G$ at the beginning of \S\ref{sec:SH}.

Note that if $S\in\Sch_B^G$ and $T\to S$ is a finitely presented $G$-quasi-projective $G$-morphism, then $T$ belongs to $\Sch_B^G$, by the definition of tameness, Lemma \ref{lem:Gquasiproj} (1), and Lemma~\ref{lem:composition}.
Moreover, any morphism in $\Sch_B^G$ is $G$-quasi-projective Nisnevich-locally on $B$, by Lemma~\ref{lem:Gquasiproj} (2).

If $S\in\Sch_B^G$, we denote by $\Sch_S^G$ the slice category over $S$, by $\Sm^G_S\subset\Sch^G_S$ the full subcategory spanned by the smooth $S$-schemes, and by $\Et^G_S\subset\Sm^G_S$ the full subcategory spanned by the étale $S$-schemes. It is clear that $\Sch_S^G$ admits finite limits and finite sums, and that they are computed in the usual way. In particular, $\Sm_S^G$ admits finite products and finite sums.

\begin{definition}
	A $G$-scheme $S$ is \emph{small} 
	if there exists a $G$-quasi-projective $G$-morphism $S\to U$ where $U$ is an affine scheme with trivial $G$-action and with the $G$-resolution property. 
\end{definition}

After having defined the $G$-equivariant Nisnevich topology, we will see that any $S\in\Sch^G_B$ admits a Nisnevich covering by small $G$-schemes (Lemma~\ref{lem:small}).
A small $G$-scheme $S$ has several convenient properties:
\begin{itemize}
	\item $S$ is separated and has the $G$-resolution property (Lemma~\ref{lem:QPresolution}).
	\item Any (quasi)-affine $G$-scheme of finite type over $S$ is $G$-(quasi)-affine (Lemma~\ref{lem:Gquasiaffine}).
	\item If $X\to S$ is $G$-quasi-projective, then $X$ is small (Lemma~\ref{lem:composition}).
	\item There exists a $G$-affine bundle $\tilde S\to S$ where $\tilde S$ is affine (Proposition~\ref{prop:jouanolou}).
\end{itemize}

\begin{remark}\label{rmk:Gdiscrete}
	Suppose that $G$ is a finite \emph{discrete} group whose order is invertible on $B$. In that case, many simplifications are possible throughout the paper. In fact, it is possible to remove all quasi-projectivity assumptions and ultimately obtain the formalism of six operations for arbitrary qcqs $B$-schemes with $G$-action, so that we might as well take $B=\Spec\Z[1/\lvert G\rvert]$.
	The main simplification comes from the fact that every qcqs $G$-scheme is locally affine in the $G$-equivariant Nisnevich topology (see Remark~\ref{rmk:nissite}).
	The reader who cares for this generality
	will have no difficulty in adapting the proofs of the main results.
	For further remarks that are relevant to this case, see Remarks~\ref{rmk:nissite}, \ref{rmk:A1invariance}, \ref{rmk:comparison}, \ref{rmk:finitesifted}, \ref{rmk:regular}, \ref{rmk:nagatachow}, and \ref{rmk:generalcdh}.
\end{remark}

\subsection{Homotopy invariance}
\label{sub:homotopy}

In this subsection, we will denote by $\scr C_S$ an arbitrary full subcategory of $\Sch_S^G$ with the property that, if $X\in\scr C_S$ and $Y\to X$ is a $G$-affine bundle, then $Y\in\scr C_S$. For example, $\scr C_S$ can be either $\Sch_S^G$ or $\Sm_S^G$.

\begin{definition}\label{def:htp}
	A presheaf $F$ on $\scr C_S$ is called \emph{homotopy invariant} if every $G$-affine bundle $Y\to X$ in $\scr C_S$ induces an equivalence $F(X)\simeq F(Y)$. We denote by $\PSh_{\htp}(\scr C_S)\subset \PSh(\scr C_S)$ the full subcategory spanned by the homotopy invariant presheaves.
\end{definition}

Since homotopy invariance is defined by a small set of conditions, the inclusion $\PSh_\htp(\scr C_S)\subset \PSh(\scr C_S)$ is an accessible localization. We will denote by $L_\htp$ the corresponding localization endofunctor of $\PSh(\scr C_S)$ whose image is $\PSh_\htp(\scr C_S)$. We say that a morphism $f$ in $\PSh(\scr C_S)$ is a \emph{homotopy equivalence} if $L_\htp(f)$ is an equivalence. Note that a colimit of homotopy invariant presheaves is homotopy invariant. In particular, $L_\htp$ preserves colimits.

A morphism $f\colon X\to Y$ in $\PSh(\scr C_S)$ will be called a \emph{strict $\A^1$-homotopy equivalence} if there exists a morphism $g\colon Y\to X$, a sequence of $\A^1$-homotopies between $g\circ f$ and $\id_X$, and a sequence of $\A^1$-homotopies between $f\circ g$ and $\id_Y$. For example, any $G$-affine bundle possessing a $G$-equivariant section is a strict $\A^1$-homotopy equivalence. It is clear that any strict $\A^1$-homotopy equivalence is a homotopy equivalence.

Our goal for the remainder of this subsection is to obtain an explicit description of the localization functor $L_\htp$. To that end, we first consider a more general situation.
Let $\scr C$ be a small $\infty$-category and let $\scr A$ be a set of morphisms in $\scr C$. A presheaf $F$ on $\scr C$ will be called \emph{$\scr A$-invariant} if it sends morphisms in $\scr A$ to equivalences. We denote by $\PSh_{\scr A}(\scr C)\subset \PSh(\scr C)$ the full subcategory spanned by the $\scr A$-invariant presheaves, and by $L_{\scr A}$ the corresponding localization endofunctor of $\PSh(\scr C)$. A morphism $f$ in $\PSh(\scr C)$ is called an \emph{$\scr A$-equivalence} if $L_{\scr A}(f)$ is an equivalence.

We say that $\scr A$ is \emph{stable under pullbacks} if, for every morphism $X'\to X$ in $\scr C$ and every $f\colon Y\to X$ in $\scr A$, there exists a cartesian square
\begin{tikzmath}
	\diagram{Y' & Y \\ X' & X \\};
	\arrows (11-) edge (-12) (21-) edge (-22) (11) edge node[left]{$f'$} (21) (12) edge node[right]{$f$} (22);
\end{tikzmath}
in $\scr C$ where $f'$ is in $\scr A$. Note that this condition on $\scr A$ does not imply that the class of $\scr A$-equivalences is stable under pullbacks. 

Let $\scr D$ be a presentable $\infty$-category with universal colimits. We recall from \cite[\S1]{GK} that a localization endofunctor $L\colon \scr D\to\scr D$ is called \emph{locally cartesian} if it commutes with local base change, \ie, if the canonical map $L(A\times_BX)\to A\times_BL(X)$ is an equivalence for any span $A\to B\from X$ in $\scr D$ with $A,B\in L(\scr D)$. This is the natural condition that guarantees that the localization $L(\scr D)$ has universal colimits.

\begin{proposition}\label{prop:cartclosed}
	Let $\scr C$ be a small $\infty$-category and $\scr A$ a set of morphisms in $\scr C$ that is stable under pullbacks. Then:
	\begin{enumerate}
		\item The localization functor $L_{\scr A}$ is given by the formula
		\[L_{\scr A}(F)(X)\simeq \colim_{Y\in\scr A_X^\op} F(Y),\]
		where $\scr A_X$ is the full subcategory of $\scr C_{/X}$ spanned by compositions of $\scr A$-morphisms.
		\item $L_{\scr A}$ is a locally cartesian localization functor.
		\item $L_{\scr A}$ preserves finite products.
	\end{enumerate}
\end{proposition}

\begin{proof}
	Let $i_X\colon \scr A_{X}\to \scr C$ be the forgetful functor.
	Given $F\in\PSh(\scr C)$, let $\tilde F$ be the presheaf on $\scr C$ defined by $\tilde F(X)=\colim i_X^*F$.
	To prove (1), we must show that:
	\begin{itemize}
		\item[(i)] $\tilde F$ is $\scr A$-invariant;
		\item[(ii)] the map $F\to \tilde F$ is an $\scr A$-equivalence.
	\end{itemize}
	If $f\colon X'\to X$ is an $\scr A$-morphism,
	 there is an adjunction
	\[f^*: \scr A_{X}^\op\rightleftarrows \scr A_{X'}^\op: f_\sharp,\]
	and $\tilde F(f)\colon \tilde F(X)\to \tilde F(X')$ is the obvious map
	\[
	\colim i_X^*F\to \colim (i_X^*F\circ f_\sharp). 
	\]
	It is an equivalence because precomposition with $f_\sharp$ is left Kan extension along $f^*$ \cite[Lemma 5.2.6.6]{HTT}.
	This proves (i).
	
	We claim that, for every $X\in\scr C$, the restriction functor $i_X^*\colon\PSh(\scr C)\to\PSh(\scr A_{X})$ preserves $\scr A$-invariant presheaves (this is obvious) as well as $\scr A$-equivalences.
	 Since $i_X^*$ preserves colimits, it suffices to show that, for every $\scr A$-morphism $Y'\to Y$,
	$i_X^*(Y')\to i_X^*(Y)$ is an $\scr A$-equivalence in $\PSh(\scr A_{X})$. Since colimits are universal in this $\infty$-category, it suffices to show that, for every $U\in \scr A_{X}$ and every map $U\to Y$ in $\scr C$, the projection $i_X^*(Y')\times_{i_X^*(Y)}U\to U$ is an $\scr A$-equivalence. But this projection can be identified with the morphism $Y'\times_YU\to U$ in $\scr A_{X}$, which is indeed (like any map in $\scr A_{X}$) an $\scr A$-equivalence.
	
	To prove (ii), it therefore suffices to prove that, for every $X\in\scr C$, $i_X^*F\to i_X^*\tilde F$ is an $\scr A$-equivalence. 
	Since $\scr A_X$ has a final object, $\PSh_{\scr A}(\scr A_{X})\subset \PSh(\scr A_{X})$ is the subcategory of constant presheaves, so $L_{\scr A}\colon \PSh(\scr A_{X})\to \PSh(\scr A_{X})$ sends $i_X^*F$ to the constant presheaf with value $\colim_{Y\in\scr A_{X}^\op}F(Y)$. Since $\tilde F$ is $\scr A$-invariant and has the same value on $X$, $i_X^*F\to i_X^*\tilde F$ is an $\scr A$-equivalence, as was to be shown. This completes the proof of (1), which immediately implies (2) by the universality of colimits in $\scr S$.
	Since $\scr A$ is stable under pullbacks, the $\infty$-category $\scr A_{X}$ has finite products and hence is cosifted, so (3) follows from (1).
\end{proof}

\begin{corollary}\label{cor:A1loc}
	The localization functor $L_\htp$ of $\PSh(\scr C_S)$ is locally cartesian and preserves finite products.
\end{corollary}

\begin{proof}
	Since $G$-affine bundles are stable under pullbacks, this is a special case of Proposition~\ref{prop:cartclosed}.
\end{proof}

\subsection{Nisnevich excision}
\label{sub:nisnevich}

Let $X$ be a $G$-scheme. A \emph{Nisnevich square} over $X$ is a cartesian square
\begin{tikzequation}\label{eqn:Nissquare}
	\diagram{W & V \\ U & X \\};
	\arrows (11-) edge[c->] (-12) (11) edge (21)
	(21-) edge[c->] node[above]{$i$} (-22) (12) edge node[right]{$p$} (22);
\end{tikzequation}
of $G$-schemes where $i$ is an open $G$-immersion, $p$ is étale, and $p$ induces an isomorphism $V\times_XZ\simeq Z$, where $Z$ is the reduced closed complement of $U$ in $X$. By Lemma~\ref{lem:invariant}, this implies that there exists a finitely presented $G$-invariant closed subscheme $Z'\subset X$, complementary to $U$, such that $V\times_XZ'\simeq Z'$.

In this subsection, we will denote by $\scr C_S$ an arbitrary full subcategory of $\Sch_S^G$ containing $\emptyset$ and with the property that, if $X\in\scr C_S$ and $Y\to X$ is an étale $G$-morphism, then $Y\in\scr C_S$. For example, $\scr C_S$ can be any of the three categories $\Sch_S^G$, $\Sm_S^G$, and $\Et_S^G$.

\begin{definition}\label{def:Nis}
	Let $S$ be a $G$-scheme. A presheaf $F$ on $\scr C_S$ is called \emph{Nisnevich excisive} if:
	\begin{itemize}
		\item $F(\emptyset)$ is contractible;
		\item for every Nisnevich square $Q$ in $\scr C_S$, $F(Q)$ is cartesian.
	\end{itemize}
	We denote by $\PSh_\Nis(\scr C_S)\subset\PSh(\scr C_S)$ the full subcategory of Nisnevich excisive presheaves.
\end{definition}

Since the property of being Nisnevich excisive is defined by a small set of conditions, $\PSh_\Nis(\scr C_S)$ is an accessible localization of $\PSh(\scr C_S)$. We denote by $L_\Nis$ the corresponding localization endofunctor of $\PSh(\scr C_S)$ whose image is $\PSh_\Nis(\scr C_S)$. We say that a morphism $f$ in $\PSh(\scr C_S)$ is a \emph{Nisnevich equivalence} if $L_\Nis(f)$ is an equivalence. Note that a filtered colimit of Nisnevich excisive presheaves is Nisnevich excisive, so that $L_\Nis$ preserves filtered colimits.

The \emph{Nisnevich topology} on $\scr C_S$ is the coarsest topology for which:
\begin{itemize}
	\item the empty sieve covers $\emptyset$;
	\item for every Nisnevich square~\eqref{eqn:Nissquare} in $\scr C_S$, $\{U\stackrel i\to X, V\stackrel p\to X\}$ generates a covering sieve.
\end{itemize}
A family of $G$-morphisms $\{p_i\colon U_i\to X\}_{i\in I}$  is called a \emph{basic Nisnevich cover} if $I$ is finite, each $p_i$ is étale, and there exists a chain of $G$-invariant finitely presented closed subschemes
\[
\emptyset = Z_0\subset Z_1\subset\dotsb\subset Z_{n-1}\subset Z_n=X
\]
such that, for each $1\leq j\leq n$, the morphism $\coprod_i U_i\to X$ splits $G$-equivariantly over $Z_j\minus Z_{j-1}$. By the proof of \cite[Proposition 2.15]{HeKrOs}, a basic Nisnevich cover is indeed a cover for the Nisnevich topology, and it is then clear that basic Nisnevich covers form a basis for the Nisnevich topology on $\scr C_S$.

If $\mathfrak U=\{U_i\to X\}$ is a family of maps in $\scr C_S$, let $\check C(\mathfrak U)$ denote the Čech nerve of the morphism
\[\coprod_i U_i\to X,\]
where the coproduct is taken in $\PSh(\scr C_S)$. Note that $\colim \check C(\mathfrak U)$ is equivalent to the image of the above morphism, which by definition is the sieve generated by $\mathfrak U$. We say that a presheaf $F$ \emph{satisfies $\mathfrak U$-descent} if $F(X)$ is the limit of the cosimplicial diagram $\Map(\check C(\mathfrak U),F)$, \ie, if $F(X)\simeq \Map(U,F)$, where $U\into X$ is the sieve generated by $\mathfrak U$.

\begin{proposition}\label{prop:topos}
	Let $S$ be a $G$-scheme and let $F$ be a presheaf on $\scr C_S$. The following are equivalent:
	\begin{enumerate}
		\item $F$ is Nisnevich excisive.
		\item $F$ satisfies $\mathfrak U$-descent for every basic Nisnevich cover $\mathfrak U$.
		\item $F$ is a sheaf for the Nisnevich topology.
	\end{enumerate}
\end{proposition}

\begin{proof}
	The equivalence of (1) and (3) is a special case of \cite[Theorem 3.2.5]{AHW}.
	The equivalence of (2) and (3) is a special case of \cite[Corollary C.2]{HoyoisGLV}.
\end{proof}

\begin{corollary}\label{cor:topos}
	The localization functor $L_\Nis$ is left exact and $\PSh_\Nis(\scr C_S)$ is an $\infty$-topos.
\end{corollary}

\begin{remark}\label{rmk:nissite}
	Let $S\in\Sch_B^G$, and let $\scr C$ denote the enlargement of $\Sm_S^G$ consisting of \emph{all} finitely presented smooth $G$-schemes over $S$. Define the Nisnevich topology on $\scr C$ via Nisnevich squares or basic Nisnevich covers, as above.
	Then the inclusion $i\colon \Sm_S^G\into \scr C$ is continuous for the Nisnevich topology (in the strong sense that the restriction functor preserves sheaves of $\infty$-groupoids), but we do not know if it is cocontinuous in general. If $G$ is finite and discrete, however, we claim that $i$ induces an equivalence between the $\infty$-topoi of sheaves.
	 Any $X\in\scr C$ admits a $G$-invariant open cover by $G$-schemes that map to affine open subschemes of $B$, which have the $G$-resolution property. For such a $G$-scheme, any cover by affine schemes has a Čech nerve that belongs to $\Sch_B^G$, by Lemma \ref{lem:Gquasiaffine}.
	 To prove our claim, it therefore suffices to show that any $X\in\scr C$ admits a Nisnevich cover by affine schemes that are smooth over $X$, and we can clearly assume that $X$ is noetherian.
	 For each $x\in X$, let $G_x=\{g\in G\suchthat gx=x\}$ and let $U_x$ be an affine $G_x$-invariant open neighborhood of $x$. By \cite[Corollary 2.19]{HeKrOs}, $\{(G\times U_x)/G_x\to X\}_{x\in X}$ is then a Nisnevich cover of $X$ by affine schemes.
\end{remark}

\begin{lemma}\label{lem:small}
	For every $S\in\Sch^G_B$, there exists a basic Nisnevich cover $\{p_i\colon U_i\to S\}$ where each $U_i$ is small and each $p_i$ is $G$-quasi-affine.
\end{lemma}

\begin{proof}
	Since $G$ is tame, there exists a (nonequivariant) basic Nisnevich cover $\{B_i\to B\}$ where each $B_i$ has the $G$-resolution property. By Lemma~\ref{lem:QPresolution}, we may assume that each $B_i$ is affine, so that each $B_i\to B$ is quasi-affine and hence $G$-quasi-affine. Since $S$ is Nisnevich-locally $G$-quasi-projective, we may further assume that $S\times_BB_i\to B_i$ is $G$-quasi-projective. Then the projections $p_i\colon U_i=S\times_BB_i \to S$ have all the desired properties.
\end{proof}

\subsection{Equivariant motivic spaces}
\label{sub:motivic}

In this subsection, $\scr C_S$ is an arbitrary full subcategory of $\Sch_S^G$ closed under finite coproducts and with the property that, if $X\in\scr C_S$ and $Y\to X$ is either an étale $G$-morphism or a $G$-affine bundle, then $Y\in\scr C_S$. For example, $\scr C_S$ can be either $\Sch^G_S$ or $\Sm^G_S$.
While we are ultimately only interested in the case $\scr C_S=\Sm^G_S$, several proofs in \S\ref{sec:functoriality} require us to consider presheaves on $\Sch^G_S$.

We denote by $\PSh_\mot(\scr C_S)\subset\PSh(\scr C_S)$ the full subcategory of presheaves that are both homotopy invariant and Nisnevich excisive.
It is an accessible localization, and we denote by $L_\mot\colon \PSh(\scr C_S)\to \PSh(\scr C_S)$ the corresponding localization endofunctor, called \emph{motivic localization}.
A morphism $f$ in $\PSh(\scr C_S)$ is called a \emph{motivic equivalence} if $L_\mot(f)$ is an equivalence.
Since the subcategories $\PSh_\htp(\scr C_S)$ and $\PSh_\Nis(\scr C_S)$ of $\PSh(\scr C_S)$ are stable under filtered colimits, we have
\[L_\mot (F) = \colim_{n\to\infty}(L_\htp\circ L_\Nis)^n(F)\]
for every presheaf $F$ on $\scr C_S$. Indeed, it is clear that the right-hand side is both homotopy invariant and Nisnevich excisive.
We will sometimes omit the functor $L_{\mot}$ from the notation, when it is clear from the context that we are working in $\PSh_\mot(\scr C_S)$.

\begin{definition}
	Let $S$ be a $G$-scheme.
	A \emph{motivic $G$-space} over $S$ is a presheaf on $\Sm^G_S$ that is homotopy invariant and Nisnevich excisive. We denote by $\H^G(S)=\PSh_\mot(\Sm^G_S)$ the $\infty$-category of motivic $G$-spaces over $S$.
\end{definition}

\begin{remark}\label{rmk:A1invariance}
	Suppose that every $X\in\scr C_S$ is Nisnevich-locally affine (\eg, $G$ is finite locally free or locally diagonalizable). Then a Nisnevich excisive presheaf on $\scr C_S$ is homotopy invariant if and only if it takes projections $X\times \A^1\to X$ to equivalences. Indeed, because $G$ is linearly reductive, every $G$-affine bundle over a small affine $G$-scheme has a $G$-section, by Lemma~\ref{lem:linred}, and hence is a strict $\A^1$-homotopy equivalence.
\end{remark}

\begin{remark}
	\label{rmk:comparison}
	Combining Remarks~\ref{rmk:nissite} and~\ref{rmk:A1invariance}, we deduce that, if $B$ is noetherian of finite Krull dimension and if $G$ is finite discrete, the $\infty$-category $\H^G(B)$ coincides with that defined in \cite{HeKrOs}. We do not know if this is true for more general tame groups $G$.
\end{remark}

\begin{proposition}\label{prop:localization}
	The localization functor $L_{\mot}$ of $\PSh(\scr C_S)$ is locally cartesian and preserves finite products. In particular, colimits in $\PSh_\mot(\scr C_S)$ are universal.
\end{proposition}

\begin{proof}
	This follows from Corollaries~\ref{cor:A1loc} and~\ref{cor:topos}.
\end{proof}

\begin{proposition}\label{prop:generators}
	Let $S$ be a $G$-scheme. 
	\begin{enumerate}
		\item The $\infty$-category $\PSh_\mot(\scr C_S)$ is generated under sifted colimits by the small affine $G$-schemes in $\scr C_S$.
		\item If a morphism in $\PSh(\scr C_S)$ is an equivalence on small affine $G$-schemes, it is a motivic equivalence.
		\item Every $X\in\scr C_S$ is compact in $\PSh_\mot(\scr C_S)$.
	\end{enumerate}
\end{proposition}

\begin{proof}
	(1) It is clear that the $\infty$-category $\PSh_\mot(\scr C_S)$ is generated under sifted colimits by $\scr C_S$, since Nisnevich excisive presheaves transform finite coproducts in $\scr C_S$ (which exist) into finite products. Let $X\in \scr C_S$, and let $\{U_i\to X\}$ be a $G$-quasi-projective basic Nisnevich cover of $X$ by small $G$-schemes (Lemma~\ref{lem:small}). 
	By Lemma~\ref{lem:composition}, the fiber products $U_{i_0}\times_X\dotsb\times_X U_{i_k}$ are also small and hence are equivalent to small affine schemes by Jouanolou's trick. Hence, the finite sum $\coprod_{i_0,\dotsc,i_k}U_{i_0}\times_X\dotsb\times_X U_{i_k}$ is also equivalent to a small affine scheme. Thus, in $\PSh_\mot(\scr C_S)$, $X$ is the colimit of a simplicial diagram whose terms are small affine schemes.
	
	(2) If $f$ is an equivalence on small affine $G$-schemes, then $L_\htp(f)$ is an equivalence on all small $G$-schemes, by Proposition~\ref{prop:cartclosed} (1), hence $L_\Nis L_\htp(f)$ is an equivalence, by Lemma~\ref{lem:small}.
	
	(3) This follows from the fact that the inclusion $\PSh_\mot(\scr C_S)\subset \PSh(\scr C_S)$ preserves filtered colimits.
\end{proof}

\subsection{Smooth closed pairs}
\label{sub:scp}

Let $S$ be a $G$-scheme. A \emph{smooth closed pair} $(X,Z)$ over $S$ is a smooth $G$-scheme $X\in\Sm^G_S$ together with a $G$-invariant closed subscheme $Z\subset X$ that is also smooth over $S$.
A morphism of smooth closed pairs $f\colon (Y,W)\to (X,Z)$ is a $G$-morphism $f\colon Y\to X$ over $S$ such that $f^{-1}(Z)=W$. We say that $f$ is \emph{Nisnevich} if it is étale and induces an isomorphism $W\simeq Z$.

\begin{definition}
	Let $S$ be a $G$-scheme.
	A morphism of smooth closed pairs $f\colon (X',Z')\to (X,Z)$ over $S$ is \emph{weakly excisive} if the square
	\begin{tikzequation}\label{eqn:wexcisive}
		\diagram{Z' & X'/(X'\minus Z') \\ Z & X/(X\minus Z) \\};
		\arrows (11-) edge[c->] (-12) (11) edge (21) (21-) edge[c->] (-22) (12) edge node[right]{$f$} (22);
	\end{tikzequation}
	in $\PSh(\Sm^G_S)$ is motivically cocartesian, \ie, becomes cocartesian in $\H^G(S)$.
\end{definition}

\begin{remark}
	Let $f\colon (X',Z')\to (X,Z)$ be a morphism of smooth closed pairs. If $f$ induces a motivic equivalence $Z'\to Z$, then $f$ is weakly excisive if and only if $X'/(X'\minus Z')\to X/(X\minus Z)$ is a motivic equivalence. On the other hand, if $f$ is weakly excisive and induces a motivic equivalence $X'\minus Z'\to X\minus Z$, then $Z/Z'\to X/X'$ is a motivic equivalence.
\end{remark}

By definition of the $\infty$-category $\H^G(S)$, it is clear that $G$-affine bundles and Nisnevich morphisms are weakly excisive.
The goal of this subsection is to obtain two other families of weakly excisive morphisms: \emph{blowups} and \emph{deformations to the normal bundle}.

\begin{samepage}
\begin{lemma}\label{lem:wexcisive}
	Let $(X'',Z'')\stackrel g\to (X',Z')\stackrel f\to (X,Z)$ be morphisms of smooth closed pairs over $S$.
	\begin{enumerate}
		\item If $g$ is weakly excisive, then $f$ is weakly excisive iff $f\circ g$ is weakly excisive.
		\item If $f$ and $f\circ g$ are weakly excisive and $f$ induces a motivic equivalence $Z'\to Z$, then $g$ is weakly excisive.
		\item If $\{U_i\to X\}$ is a basic Nisnevich cover of $X$ such that $f_{U_{i_1}\times_X\dotsb\times_XU_{i_k}}$ is weakly excisive for every nonempty family of indices $(i_1,\dotsc,i_k)$, then $f$ is weakly excisive.
	\end{enumerate}
\end{lemma}
\end{samepage}

\begin{proof}
	(1) and (2) are obvious. 
	By Nisnevich descent, the square~\eqref{eqn:wexcisive} associated with $f$ is a colimit in $\H^G(S)$ of squares associated with $f_{U_{i_1}\times_X\dotsb\times_XU_{i_k}}$, whence (3).
\end{proof}

\begin{lemma}\label{lem:triple}
	Let $(X',Z',W')\to (X,Z,W)$ be a morphism of smooth closed triples over $S$. Suppose that $(X',Z')\to (X,Z)$, $(Z',W')\to(Z,W)$ and $(X'\minus W',Z'\minus W')\to(X\minus W,Z\minus W)$ are weakly excisive. Then $(X',W')\to (X,W)$ is weakly excisive.
\end{lemma}

\begin{proof}
	We must show that the boundary of the following diagram is motivically cocartesian:
	\begin{tikzmath}
		\diagram{W' & Z'/(Z'\minus W') & X'/(X'\minus W') \\ W & Z/(Z\minus W) & X/(X\minus W)\rlap. \\};
		\arrows (11-) edge[c->] (-12) (12-) edge[c->] (-13)
		(21-) edge[c->] (-22) (22-) edge[c->] (-23)
		(11) edge (21) (12) edge (22) (13) edge(23);
	\end{tikzmath}
	The first square is motivically cocartesian since $(Z',W')\to (Z,W)$ is weakly excisive. The second square is the cofiber of the obvious morphism from the square associated with $(X'\minus W',Z'\minus W')\to (X\minus W,Z\minus W)$ to the square associated with $(X',Z')\to (X,Z)$, and hence it is also motivically cocartesian.
\end{proof}

Let $(X,Z)$ be a smooth closed pair over $S$. If $p\colon B_ZX\to X$ is the blowup of $X$ at $Z$, then $(B_Z(X), p^{-1}(Z))$ is a smooth closed pair over $S$, depending functorially on $(X,Z)$. Moreover, $p$ is a morphism of smooth closed pair $(B_ZX,p^{-1}(Z))\to (X,Z)$, natural in $(X,Z)$. 
We denote by $N_ZX$ the normal bundle of $Z$ in $X$. The deformation space \cite[Chapter 5]{Fulton} of $(X,Z)$ is a smooth $G$-scheme $D_ZX$ over $S\times\A^1$ whose restrictions to $S\times 1$ and $S\times 0$ are canonically isomorphic to $X$ and $N_ZX$, respectively. Explicitly,
\[D_ZX=B_{Z\times 0}(X\times\A^1)\minus B_{Z\times 0}(X\times 0).\]
The closed $G$-immersion $Z\times \A^1\into X\times\A^1$ lifts uniquely to a closed $G$-immersion $Z\times \A^1\into D_{Z}X$. It is clear that the smooth closed pair $(D_ZX,Z\times\A^1)$ varies functorially with $(X,Z)$, and that the inclusions $ 1\into \A^1\hookleftarrow  0$ induce natural morphisms of pairs
\[i_1:(X,Z)\to (D_ZX,Z\times\A^1)\from (N_ZX,Z):i_0.\]

\begin{lemma}\label{lem:SCP}
	Let $S$ be a $G$-scheme and let $\mathbf P$ be a class of smooth closed pairs over $S$.
	Suppose that the following conditions hold for every smooth closed pair $(X,Z)$:
	\begin{enumerate}
		\item If $\{U_i\to X\}$ is a basic Nisnevich cover of $X$ and $(U_{i_1}\times_X\dotsb\times_X U_{i_k},Z\times_X U_{i_1}\times_X\dotsb\times_X U_{i_k})\in\mathbf P$ for every nonempty family of indices $(i_1,\dotsc,i_k)$, then $(X,Z)\in\mathbf P$.
		\item If $(Y,W)\to (X,Z)$ is a $G$-affine bundle and $(Y,W)\in\mathbf P$, then $(X,Z)\in\mathbf P$.
		\item If $(Y,W)\to (X,Z)$ is a Nisnevich morphism, then $(X,Z)\in\mathbf P$ if and only if $(Y,W)\in\mathbf P$.
		\item If $\scr E$ is a locally free $G$-module on $Z$, then $(\V(\scr E),Z)\in\mathbf P$.
	\end{enumerate}
	Then $\mathbf P$ contains all smooth closed pairs.
\end{lemma}

\begin{proof}
	Let $(X,Z)$ be an arbitrary smooth closed pair over $S$. To show that it belongs to $\mathbf P$, we can assume, by (1) and Lemma~\ref{lem:small}, that $X$ is small. By Jouanolou's trick and (2), we can then also assume that $X$ is affine. By Corollary~\ref{cor:retraction}, we can find a cartesian square
		\begin{tikzmath}
			\diagram{Z & X' \\ Z & X \\};
			\arrows (11-) edge[c->] node[above]{$t$} (-12) (11) edge[-,vshift=1pt] (21) edge[-,vshift=-1pt] (21) (21-) edge[c->] (-22) (12) edge node[right]{$p$} (22);
		\end{tikzmath}
		where $X'$ is affine, $p$ is étale at $t(Z)$, and $t$ admits a $G$-retraction $r\colon X'\to Z$. Since $t$ is a quasi-regular $G$-immersion, we can apply Proposition~\ref{prop:linearization} to $r$ and get a $G$-morphism $h\colon X'\to\V(\scr N_t)$ over $Z$ that is étale at $t(Z)$ and such that $ht=z$, where $z\colon Z\into \V(\scr N_t)$ is the zero section.
		Let $U\subset X'$ be the intersection of the étale locus of $p$ and the étale locus of $h$. By Lemma~\ref{lem:locus}, $U$ is a $G$-invariant open subscheme of $X'$. Shrinking $U$ if necessary, we can assume that it is quasi-compact and that the square
		\begin{tikzmath}
			\diagram{Z & U \\ Z & \bb V(\scr N_t) \\};
			\arrows (11-) edge[c->] node[above]{$t$} (-12) (11) edge[-,vshift=1pt] (21) edge[-,vshift=-1pt] (21) (21-) edge[c->] node[above]{$z$} (-22) (12) edge node[right]{$h$} (22);
		\end{tikzmath}
		is cartesian, by Lemma~\ref{lem:cartesian}.
		We therefore have two étale $G$-morphisms
		\[X\from U\to \V(\scr N_t)\]
		that are isomorphisms over $Z$.
		By (4), the pair $(\bb V(\scr N_t),Z)$ belongs to $\mathbf P$. By two applications of (3), we deduce that $(X,Z)\in\mathbf P$.
\end{proof}

\begin{theorem}
	\label{thm:wexcisive}
	Let $S$ be a $G$-scheme. For every smooth closed pair $(X,Z)$ over $S$, the morphisms
	\[(B_ZX,p^{-1}(Z))\xrightarrow{p} (X,Z)\quad\text{and}\quad (X,Z)\xrightarrow{i_1} (D_ZX,Z\times\A^1)\xleftarrow{i_0} (N_ZX,Z)\]
	are weakly excisive.
\end{theorem}

\begin{proof}
	Let $\mathbf P$ be the class of smooth closed pairs $(X,Z)$ over $S$ for which the conclusion of the theorem holds. We will show that $\mathbf P$ satisfies conditions (1)--(4) of Lemma~\ref{lem:SCP}. Since blowups commute with flat base change, it is clear that the functors $B$, $D$, and $N$ preserve $G$-affine bundles and Nisnevich morphisms. If $f\colon (Y,W)\to(X,Z)$ is such a morphism, it follows from Lemma~\ref{lem:wexcisive} (1,2) that $(X,Z)\in\mathbf P$ iff $(Y,W)\in\mathbf P$, which proves conditions (2) and (3). Condition (1) follows easily from Lemma~\ref{lem:wexcisive} (3).
	
	It remains to show that $\mathbf P$ satisfies condition (4), \ie, that the theorem holds when $X=\V(\scr E)$ for some locally free $G$-module $\scr E$ on $Z$. In that case, $(B_ZX,p^{-1}(Z))$ can be identified with the smooth closed pair $(\V(\scr O_{\P(\scr E)}(1)), \P(\scr E))$. The square~\eqref{eqn:wexcisive} for $p$ is the outer square in the following diagram:
	\begin{tikzmath}
		\def\rowsep{1.5em}
		\def\colsep{1.5em}
		\diagram{p^{-1}(Z) & B_ZX & \frac{B_ZX}{B_ZX\minus p^{-1}(Z)} \\ 
		Z & X & \frac{X}{X\minus Z}\rlap. \\};
		\arrows (11-) edge[c->] (-12) (21-) edge[c->] (-22) (12-) edge (-13) (22-) edge (-23)
		(11) edge (21) (12) edge node[right]{$p$} (22) (13) edge (23);
	\end{tikzmath}
	The first square is motivically cocartesian since the horizontal maps are sections of $G$-vector bundles and hence are homotopy equivalences. The second square is already cocartesian in $\PSh(\Sm^G_S)$, since $B_ZX\minus p^{-1}(Z)\to X\minus Z$ is a monomorphism. This proves that $p$ is weakly excisive.
	
	Similarly, when $X=\V(\scr E)$, $(D_ZX,Z\times\A^1)$ can be identified with the smooth closed pair $(\V(\scr L), \V(\scr L)_Z)$, where $\scr L$ is the restriction to $X$ of $\scr O_{\P(\scr E\oplus\scr O)}(1)$. On the other hand, $(N_ZX,Z)$ is canonically isomorphic to $(X,Z)$. Under these identifications, $i_1$ and $i_0$ are both sections of the $G$-line bundle $(\V(\scr L), \V(\scr L)_Z)\to (X,Z)$ and hence are weakly excisive by Lemma~\ref{lem:wexcisive} (2).
\end{proof}

By the second part of Theorem~\ref{thm:wexcisive}, the morphisms $i_1$ and $i_0$ induce an equivalence
\[\Pi=\Pi_{X,Z}\colon\frac{X}{X\minus Z}\simeq \frac{N_ZX}{N_ZX\minus Z}\]
in $\H^G(S)$, natural in the smooth closed pair $(X,Z)$. It is called the \emph{purity equivalence}.

Let $(X,Z,W)$ be a smooth closed triple over $S$. Then $W\times\A^1$ can be identified with a $G$-invariant closed subscheme of $D_ZX$, and we can form the smooth closed pair $(D_{W\times\A^1}(D_ZX), W\times\A^2)$ over $S\times\A^2$. Pulling it back along the closed immersions
\begin{tikzmath}
	\def\colsep{1em}
	\def\rowsep{2em}
	\diagram{
	(1,1) & 1\times\A^1 & (1,0) \\
	\A^1\times 1 & \A^2 & \A^1\times 0 \\
	(0,1) & 0\times\A^1 & (0,0)\rlap, \\
	};
	\arrows
	(11-) edge (-12) (12-) edge[<-] (-13)
	(21-) edge (-22) (22-) edge[<-] (-23)
	(31-) edge (-32) (32-) edge[<-] (-33)
	(11) edge (21) (21) edge[<-] (31)
	(12) edge (22) (22) edge[<-] (32)
	(13) edge (23) (23) edge[<-] (33);
\end{tikzmath}
we obtain the following commutative diagram of smooth closed pairs over $S$:
\begin{tikzequation}\label{eqn:dd}
	\diagram{
	(X,W) & (D_ZX,W\times\A^1) & (N_ZX,W) \\
	(D_WX,W\times\A^1) & (D_{W\times\A^1}(D_ZX),W\times\A^2) & (D_{W}(N_ZX),W\times\A^1) \\
	(N_WX,W) & (N_{W\times\A^1}(D_ZX),W\times\A^1) & (N_{W}(N_ZX),W)\rlap. \\
	};
	\arrows
	(11-) edge (-12) (12-) edge[<-] (-13)
	(21-) edge (-22) (22-) edge[<-] (-23)
	(31-) edge (-32) (32-) edge[<-] (-33)
	(11) edge (21) (21) edge[<-] (31)
	(12) edge (22) (22) edge[<-] (32)
	(13) edge (23) (23) edge[<-] (33);
\end{tikzequation}
Note that the bottom row of~\eqref{eqn:dd} is canonically isomorphic to
\[(N_WX,W)\xrightarrow{i_1} (D_{N_WZ}(N_WX),W\times\A^1)\xleftarrow{i_0} (N_{N_WZ}(N_WX), W).\]

\begin{corollary}\label{cor:wexcisive}
	For every smooth closed triple $(X,Z,W)$ over $S$, all morphisms in~\eqref{eqn:dd} are weakly excisive.
\end{corollary}

\begin{proof}
	Each column is weakly excisive by Theorem~\ref{thm:wexcisive}.
	The top row is seen to be weakly excisive by applying Lemma~\ref{lem:triple} to the morphisms of triples
	\[(X,Z,W)\to (D_ZX,Z\times\A^1,W\times\A^1)\from (N_ZX,Z,W),\]
	using Theorem~\ref{thm:wexcisive} twice. 
	The remaining rows are automatically weakly excisive by Lemma~\ref{lem:wexcisive} (1,2).
\end{proof}

The top row of~\eqref{eqn:dd} provides an equivalence
\[
\frac{X}{X\minus W} \simeq \frac{N_ZX}{N_ZX\minus W}
\]
in $\H^G(S)$, which we also denote by $\Pi_{X,Z}$. It is more precisely the composition
\[
\frac{X}{X\minus W}=\frac{X/(X\minus Z)}{(X\minus W)/(X\minus Z)}\simeq \frac{N_ZX/(\dots)}{N_{Z\minus W}(X\minus W)/(\dots)}=\frac{N_ZX}{N_ZX\minus W},
\]
where the middle equivalence is $\Pi_{X,Z}$ in the numerator and $\Pi_{X\minus W,Z\minus W}$ in the denominator.

Let $X\in\Sm_S^G$ and let
\[
0\to U\to V\to W\to 0
\]
be a short exact sequence of $G$-vector bundles on $X$.
The top row of~\eqref{eqn:dd} for the smooth closed triple $(V,U,X)$ is a pair of morphisms
\begin{equation}
	\label{eqn:VBdeformation}
(V,X)\to (D_UV,X\times\A^1) \from (U\times_XW,X)
\end{equation}
that are weakly excisive
by Corollary~\ref{cor:wexcisive}.
In particular, we obtain an equivalence
\[
\Psi=\Pi_{V,U}\colon \frac{V}{V\minus X} \simeq \frac{U\times_XW}{(U\times_XW)\minus X}
\]
in $\H^G(S)$. More generally, if $Z\subset X$ is a smooth $G$-invariant closed subscheme, the morphisms~\eqref{eqn:VBdeformation} are also weakly excisive after replacing $X$ by $Z$, so that we have a canonical equivalence
\[
\Psi\colon \frac{V}{V\minus Z} \simeq \frac{U\times_XW}{(U\times_XW)\minus Z}
\]
in $\H^G(S)$, compatible with the previous one.
Note that the diagram~\eqref{eqn:VBdeformation} and hence the equivalence $\Psi$ are natural for universal monomorphisms of short exact sequences of $G$-vector bundles.

\begin{remark}\label{rmk:splitting}
	Suppose given a short exact sequence as above and a splitting $\sigma\colon V\simeq U\times_XW$. Then the equivalence $\Psi$ is induced by $\sigma$. Indeed, identifying $V$ with $U\times_XW$ by means of $\sigma$, the morphisms~\eqref{eqn:VBdeformation} can be identified with sections of the same $G$-line bundle over $V$, so that $\Psi \sigma^{-1}$ is the identity in $\H^G(S)$. 
\end{remark}

With this notation in place, the main point of Corollary~\ref{cor:wexcisive} is that it implies the commutativity of the following square of equivalences in $\H^G(S)$, for every smooth closed triple $(X,Z,W)$:
\begin{tikzmath}
	\diagram{\frac X{X\minus W} & \frac{N_ZX}{N_ZX\minus W} \\ \frac {N_WX}{N_WX\minus X} & \frac{N_WN_ZX}{N_WN_ZX\minus W}\rlap. \\};
	\arrows (11-) edge node[above]{$\Pi_{X,Z}$} (-12) (11) edge node[left]{$\Pi_{X,W}$} (21) (12) edge node[right]{$\Pi_{N_ZX,W}$} (22) (21-) edge node[above]{$\Psi$} (-22);
\end{tikzmath}
Here, $\Psi$ is induced by the canonical short exact sequence
\[
0\to N_WZ\to N_WX \to N_ZX\times_ZW\to 0
\]
of $G$-vector bundles on $W$.

\begin{remark}
	More generally, for a smooth closed tuple $(X_0,\dotsc,X_n)$, one can show using an $n$-fold deformation to the normal bundle that all the equivalences
	\[\frac{X_0}{X_0\minus X_n}\simeq \frac{N_{X_n}N_{X_{n-1}}\dots N_{X_1}X_0}{N_{X_n}N_{X_{n-1}}\dots N_{X_1}X_0\minus X_n}\]
	that one can define using Theorem~\ref{thm:wexcisive} fit in a commutative $n$-cube.
	Note that Corollary~\ref{cor:wexcisive} implies the commutativity of this $n$-cube in the homotopy $1$-category, which is sufficient for many purposes.
\end{remark}

\section{Functoriality}
\label{sec:functoriality}

We discuss the functoriality in $S$ of the $\infty$-category $\H^G(S)$.
We remark that all the results of \S\ref{sub:easy} and~\S\ref{sub:exactness} remain valid if we replace $\Sm^G_S$ by $\Sch^G_S$ (and we can then remove the smoothness assumption in Propositions \ref{prop:smoothBC} and~\ref{prop:projection}). However, the gluing theorem of \S\ref{sub:gluing} uses smoothness in an essential way.

\subsection{Easy results}
\label{sub:easy}

It is clear that the empty scheme, Nisnevich squares, and $G$-affine bundles are preserved by base change. As a result, if $f\colon T\to S$ is a $G$-morphism, then the functor 
\[
f_*\colon\PSh(\Sm^G_T)\to\PSh(\Sm^G_S), \quad f_*(F)(X)=F(X\times_ST),
\]
preserves Nisnevich excisive presheaves, homotopy invariant presheaves, and motivic $G$-spaces. We still denote by $f_*\colon \H^G(T)\to \H^G(S)$ the restriction of $f_*$. Since it preserves limits, it admits a left adjoint.
 We will overload the notation $f^*$ to denote a left adjoint to $f_*$ in any context, and when $f^*$ itself admits a left adjoint we will denote it by $f_\sharp$. Note that $f_*\colon \H^G(T)\to \H^G(S)$ preserves filtered colimits since the inclusion $\H^G(S)\subset \PSh(\Sm_S^G)$ does.
 
 If $f\colon T\to S$ is a smooth $G$-morphism, then the forgetful functor $\Sm^G_T\to \Sm^G_S$ preserves the empty scheme, Nisnevich squares, and $G$-affine bundles. It follows that $f^*\colon \PSh(\Sm^G_S)\to \PSh(\Sm^G_T)$, which in this case is precomposition with the forgetful functor, preserves Nisnevich excisive presheaves, homotopy invariant presheaves, and motivic $G$-spaces. In particular, it restricts to a functor $\H^G(S)\to \H^G(T)$ that is automatically left adjoint to $f_*$. We deduce that, for $f$ smooth, $f^*\colon \H^G(S)\to \H^G(T)$ has a left adjoint $f_\sharp$. 
 
In summary, for a $G$-morphism $f\colon T\to S$, we have commutative squares
\[
\begin{tikzpicture}
	\diagram{\PSh(\Sm_T^G) & \PSh(\Sm_S^G) \\ \H^G(T) & \H^G(S)\rlap, \\};
	\arrows (11-) edge node[above]{$f_*$} (-12) (21-) edge node[below]{$f_*$} (-22)
	(11) edge[<-right hook] (21) (12) edge[<-right hook] (22);
\end{tikzpicture}
\qquad
\begin{tikzpicture}
	\diagram{\PSh(\Sm_T^G) & \PSh(\Sm_S^G) \\ \H^G(T) & \H^G(S)\rlap. \\};
	\arrows (11-) edge[<-] node[above]{$f^*$} (-12) (21-) edge[<-] node[below]{$f^*$} (-22)
	(11) edge (21) (12) edge (22);
\end{tikzpicture}
\]
If $f$ is smooth, we moreover have commutative squares
\[
\begin{tikzpicture}
	\diagram{\PSh(\Sm_S^G) & \PSh(\Sm_T^G) \\ \H^G(S) & \H^G(T)\rlap, \\};
	\arrows (11-) edge node[above]{$f^*$} (-12) (21-) edge node[below]{$f^*$} (-22)
	(11) edge[<-right hook] (21) (12) edge[<-right hook] (22);
\end{tikzpicture}
\qquad
\begin{tikzpicture}
	\diagram{\PSh(\Sm_S^G) & \PSh(\Sm_T^G) \\ \H^G(S) & \H^G(T)\rlap. \\};
	\arrows (11-) edge[<-] node[above]{$f_\sharp$} (-12) (21-) edge[<-] node[below]{$f_\sharp$} (-22)
	(11) edge (21) (12) edge (22);
\end{tikzpicture}
\]

\begin{proposition}[Monoidality]
	Let $f\colon Y\to X$ be a $G$-morphism. Then $f^*\colon\H^G(X)\to\H^G(Y)$ preserves finite products.
\end{proposition}

\begin{proof}
	Since $L_\mot$ preserves finite products (Proposition~\ref{prop:localization}), this follows from the analogous fact at the level of presheaves, which is obvious.
\end{proof}

\begin{proposition}[Smooth base change]
	\label{prop:smoothBC}
	Let
	\begin{tikzmath}
		\diagram{Y' & Y \\ X' & X \\};
		\arrows (11-) edge node[above]{$g$} (-12) (11) edge node[left]{$q$} (21) (21-) edge node[below]{$f$} (-22) (12) edge node[right]{$p$} (22);
	\end{tikzmath}
	be a cartesian square of $G$-schemes where $p$ is smooth. Then the exchange transformations
	\begin{gather*}
		\Ex_\sharp^*\colon q_\sharp g^*\to f^*p_\sharp\colon \H^G(Y)\to\H^G(X'),\\
		\Ex_*^*\colon p^*f_*\to g_*q^*\colon \H^G(X')\to \H^G(Y)
	\end{gather*}
	are equivalences.
\end{proposition}

\begin{proof}
	The second exchange transformation is the mate of the first one, so it suffices to show that $q_\sharp g^*\to f^*p_\sharp$ is an equivalence. This transformation is $L_\mot$ of the analogous exchange transformation at the level of presheaves, which is clearly an equivalence.
\end{proof}

\begin{proposition}[Smooth projection formula]
	\label{prop:projection}
	Let $f\colon Y\to X$ be a smooth $G$-morphism, let $A,C\in\H^G(X)$, and let $B\in\H^G(Y)$. Then the canonical maps
	\begin{gather*}
		f_\sharp(f^*A\times B)\to A\times f_\sharp B, \\
		f^*\Hom(A,C)\to \Hom(f^*A, f^*C)
	\end{gather*}
	are equivalences in $\H^G(X)$ and $\H^G(Y)$, respectively.
\end{proposition}

\begin{proof}
	By adjunction, it suffices to show that the first map is an equivalence. Since $L_\mot$ preserves finite products (Proposition~\ref{prop:localization}), this follows from the projection formula at the level of presheaves. 
\end{proof}

\begin{proposition}[Homotopy invariance]
	\label{prop:hi}
	Let $f\colon T\to S$ be a $G$-affine bundle.
	\begin{enumerate}
		\item The functor $f^*\colon \PSh(\Sm^G_S)\to \PSh(\Sm^G_T)$ detects homotopy and motivic equivalences.
		\item The functor $f^*\colon\H^G(S)\to\H^G(T)$ is fully faithful.
	\end{enumerate}
\end{proposition}

\begin{proof}
	Since $f$ is smooth, $f^*$ has a left adjoint $f_\sharp$ that preserves homotopy and motivic equivalences.
	 If $X\in\Sm^G_S$, the counit $\epsilon_X\colon f_\sharp f^*X\to X$ is the $G$-affine bundle $X\times_ST\to X$. Since $f_\sharp$ and $f^*$ preserve colimits, $\epsilon_X$ is a homotopy equivalence for every $X\in\PSh(\Sm_S^G)$. This easily implies the results.
\end{proof}

\begin{proposition}[Nisnevich separation]
	\label{prop:Nissep}
	Let $\{f_i\colon U_i\to S\}$ be a Nisnevich cover of a $G$-scheme $S$.
	\begin{enumerate}
		\item The family of functors $\{f_i^*\colon\PSh(\Sm_S^G)\to\PSh(\Sm_{U_i}^G)\}$ detects Nisnevich and motivic equivalences.
		\item The family of functors $\{f_i^*\colon\H^G(S)\to\H^G(U_i)\}$ is conservative.
	\end{enumerate}
\end{proposition}

\begin{proof}
	Clearly, (1) implies (2). To prove (1), we may assume that each $f_i$ is smooth.
	Let $h$ be a morphism in $\PSh(\Sm_S^G)$ such that $f_i^*(h)$ is a Nisnevich (\resp{} motivic) equivalence for all $i$.
	Denote by \[f_{i_1\dotso i_n}\colon U_{i_1}\times_S\dotsb\times_SU_{i_n}\to S\] an $n$-fold fiber product of the given covering. Then $f_{i_1\dotso i_n}^*(h)$ is a Nisnevich (\resp{} motivic) equivalence. Consider the augmented simplicial object $C_\bullet \to \id$ whose $n$th term $C_n$ is the endofunctor
	\[
	F\mapsto \coprod_{i_1,\dotsc,i_n}(f_{i_1\dotso i_n})_\sharp f_{i_1\dotso i_n}^*F.
	\]
	For every $X\in\Sm_S^G$, $\colim C_\bullet(X)\to X$ is a Nisnevich covering sieve. Thus, for every $F\in\PSh(\Sm_S^G)$, $\colim C_\bullet(F)\to F$ is a Nisnevich equivalence.
	Since $C_\bullet(h)$ is a degreewise Nisnevich (\resp{} motivic) equivalence, it follows by 2-out-of-3 that $h$ is a Nisnevich (\resp{} motivic) equivalence, as desired.
\end{proof}

\begin{proposition}[Affine resolutions]
	\label{prop:affinecover}
	Let $S$ be a $G$-scheme. Then there exists a family of smooth $G$-quasi-affine morphisms $\{f_i\colon U_i\to S\}$, where each $U_i$ is small and affine, such that:
	\begin{enumerate}
		\item the family of functors $\{f_i^*\colon\PSh(\Sm_S^G)\to\PSh(\Sm_{U_i}^G)\}$ detects motivic equivalences;
		\item the family of functors $\{f_i^*\colon\H^G(S)\to\H^G(U_i)\}$ is conservative.
	\end{enumerate}
\end{proposition}

\begin{proof}
	Let $\{V_i\to S\}$ be a $G$-quasi-affine basic Nisnevich cover of $S$ by small $G$-schemes (Lemma~\ref{lem:small}). By Jouanolou's trick (Proposition~\ref{prop:jouanolou}), there exist $G$-affine bundles $U_i\to V_i$ where $U_i$ is affine. These $U_i$ form the desired family by Propositions~\ref{prop:hi} and \ref{prop:Nissep}.
\end{proof}

\begin{proposition}\label{prop:openimmersion}
	Let $j\colon U\into X$ be an open $G$-immersion. Then the functors
	\begin{gather*}
		j_*\colon\H^G(U)\to\H^G(X), \\
		j_\sharp\colon\H^G(U)\to\H^G(X)
	\end{gather*}
	are fully faithful.
\end{proposition}

\begin{proof}
	It is clear that the unit $\id\to j^*j_\sharp$ is an equivalence.
\end{proof}

To make sense of the next proposition, we remark that one can successively construct functors $(\Sch_B^G)^\op\to \Cat_\infty$ that send a $G$-scheme $S$ to $\Sm_S^G$, $\PSh(\Sm_S^G)$, $\PSh_\Nis(\Sm_S^G)$, and $\H^G(S)$. We refer to \cite[\S9.1]{RobaloThesis} for details on these constructions.

\begin{proposition}[Nisnevich descent]
\label{prop:Nislocal}
\leavevmode
\begin{enumerate}
	\item The assignment $S\mapsto \PSh_\Nis(\Sm^G_S)$, $f\mapsto f^*$, is a Nisnevich sheaf of $\infty$-categories on $\Sch_B^G$.
	\item The assignment $S\mapsto \H^G(S)$, $f\mapsto f^*$, is a Nisnevich sheaf of $\infty$-categories on $\Sch_B^G$.
\end{enumerate}
\end{proposition}

\begin{proof}
	(1) It suffices to show that $V\mapsto \PSh_\Nis(\Sm^G_V)$ is a Nisnevich sheaf on $\Et^G_S$ for all $S\in\Sch_B^G$.
	Note that if $V\in\Et^G_S$, then $\Sm_V^G\simeq (\Sm_S^G)_{/V}$, since if $gf$ is smooth and $g$ is étale then $f$ is smooth.
	Thus,
	\[
	\PSh_\Nis(\Sm_V^G)\simeq \PSh_\Nis(\Sm_S^G)_{/V}.
	\]
	The fact that $V\mapsto \PSh_\Nis(\Sm_S^G)_{/V}$ is a Nisnevich sheaf on $\Et^G_S$ (and even on $\Sm^G_S$) follows from Corollary~\ref{cor:topos} and general descent theory for $\infty$-topoi \cite[Theorem 6.1.3.9 (3)]{HTT}.
	
	(2) Again, it suffices to show that $V\mapsto \H^G(V)$ is a Nisnevich sheaf on $\Et^G_S$ for all $S\in\Sch_B^G$. The inclusions $\H^G(V)\subset \PSh_\Nis(\Sm^G_V)$ are the components of a natural transformation on $\Et^G_S$, since all morphisms in $\Et^G_S$ are smooth. Given (1) and the fact that limits of fully faithful functors are fully faithful, the sheaf condition is reduced to the following statement: for every basic Nisnevich cover $\{p_i\colon U_i\to V\}$ in $\Et_S^G$ and every $F\in \PSh_\Nis(\Sm^G_V)$, if $p_i^*(F)$ is homotopy invariant for all $i$, then $F$ is homotopy invariant. This is clear, since $G$-affine bundles are stable under base change.
\end{proof}

\subsection{Exactness properties of pushforwards}
\label{sub:exactness}

\begin{lemma}\label{lem:bundlelift}
	Let $f\colon T\to S$ be an affine $G$-morphism where $S$ has the $G$-resolution property, and let $Y\to T$ be a $G$-affine bundle. Then there exists a $G$-affine bundle $V\to S$ and a $G$-morphism $V_T\to Y$ over $T$.
\end{lemma}

\begin{proof}
	Let $\scr E\onto \scr O_T$ be the epimorphism of locally free $G$-modules corresponding to $Y$ (see \S\ref{sub:jouanolou}), and let $\scr M$ be defined by the cartesian square
	\begin{tikzmath}
		\diagram{\scr M & \scr O_S \\ f_*(\scr E) & f_*(\scr O_T) \\};
		\arrows (11-) edge (-12) (11) edge (21) (21-) edge (-22) (12) edge node[right]{$\eta$} (22);
	\end{tikzmath}
	in $\QCoh^G(S)$. Since $f$ is affine, the horizontal arrows are epimorphisms. Since $\scr M$ is the colimit of its finitely generated quasi-coherent $G$-submodules (Lemma~\ref{lem:thomason}), there exists a finitely generated quasi-coherent $G$-submodule $\scr N\subset\scr M$ such that $\scr N\to\scr M\to\scr O_S$ is surjective. By the $G$-resolution property, there exists a locally free $G$-module $\scr F$ and an epimorphism $\scr F\onto\scr N$. Then $\scr F\onto \scr O_S$ defines a $G$-affine bundle $V$ over $S$ and the morphism $f^*(\scr F)\to\scr E$ over $\scr O_T$ defines a $G$-morphism $V_T\to Y$ over $T$, as desired.
\end{proof}

\begin{lemma}\label{lem:affinePF}
	Let $f\colon T\to S$ be an affine $G$-morphism where $S$ has the $G$-resolution property. Then the functor $f_*\colon \PSh(\Sm^G_T)\to\PSh(\Sm^G_S)$ preserves homotopy equivalences.
\end{lemma}

\begin{proof}
	Since $f_*$ preserves colimits, it suffices to show that it sends $G$-affine bundles to homotopy equivalences. We prove more generally that, if $Y\to X$ is a $G$-affine bundle in $\Sch^G_T$, then $f_*Y\to f_*X$ is a homotopy equivalence in $\PSh(\Sch^G_S)$. By universality of colimits in $\PSh(\Sch^G_S)$, it suffices to show that, for every $G$-morphism $p\colon U\to S$ and every $U_T\to X$, the projection $f_*Y\times_{f_*X}U\to U$ is a homotopy equivalence. Consider the cartesian square
	\begin{tikzmath}
		\diagram{U_T & U \\ T & S\rlap. \\};
		\arrows (11-) edge node[above]{$g$} (-12) (11) edge node[left]{$q$} (21) (12) edge node[right]{$p$} (22) (21-) edge node[below]{$f$} (-22);
	\end{tikzmath}
	We then have $p^* f_*\simeq g_* q^*$ and, since $\Sch^G_U$ has pullbacks, the projection $f_*Y\times_{f_*X}U\to U$ can be identified with $p_\sharp$ of the projection $g_* Y_U\times_{g_* X_U}U\to U$ in $\PSh(\Sch^G_U)$. Since $p_\sharp$ preserves homotopy equivalences, we may as well assume that $U=S$, so that $U_T=T$. Since $f_*$ preserves limits, we then have
	\[f_*Y\times_{f_*X}S\simeq f_* (Y\times_XT).\]
	We are thus reduced to proving the following statement: if $Y\to T$ is a $G$-affine bundle, then $f_*Y$ is homotopically contractible in $\PSh(\Sch^G_S)$.
	 By Lemma~\ref{lem:bundlelift}, there exists a $G$-affine bundle $r\colon V\to S$ and a $G$-morphism $V_T\to Y$ over $T$. By Proposition~\ref{prop:hi} (1), it remains to prove that $r^*f_*Y$ is homotopically contractible in $\PSh(\Sch^G_V)$. This presheaf can be identified with the pushforward along $V_T\to V$ of the $G$-affine bundle $V_T\times_TY\to V_T$. But by choice of $V$, this $G$-affine bundle has a $G$-section and hence is a strict $\A^1$-homotopy equivalence, and it is clear that pushforwards preserve strict $\A^1$-homotopy equivalences.
\end{proof}

Let $\scr C$ be a small $\infty$-category with an initial object $\emptyset$ (\resp{} with finite coproducts).
We denote by $\PSh_\emptyset(\scr C)$ (\resp{} $\PSh_\Sigma(\scr C)$) the full subcategory of $\PSh(\scr C)$ consisting of those presheaves $F$ such that $F(\emptyset)\simeq *$ (\resp{} that transform finite coproducts into finite products).
Recall that the $\infty$-category $\PSh_\Sigma(\scr C)$ is freely generated by $\scr C$ under sifted colimits \cite[Proposition 5.5.8.15]{HTT}. 
We say that an $\infty$-category $\scr A$ is \emph{weakly contractible} if the map $\scr A\to *$ is cofinal, or, equivalently, if the $\infty$-groupoid completion of $\scr A$ is contractible; colimits indexed by such $\infty$-categories will be called \emph{weakly contractible colimits}.
Clearly, a presheaf $F\in\PSh(\scr C)$ belongs to $\PSh_\emptyset(\scr C)$ if and only if its $\infty$-category of elements is weakly contractible. It follows that $\PSh_\emptyset(\scr C)$ is freely generated by $\scr C$ under weakly contractible colimits.
We denote by $a_\emptyset\colon\PSh(\scr C)\to \PSh_\emptyset(\scr C)$ the left adjoint to the inclusion.

\begin{lemma}\label{lem:closedexact}
	Let $i\colon Z\into S$ be a closed $G$-immersion between $G$-schemes that are affine over $B$. Then $i_*\colon \PSh_\emptyset(\Sm^G_Z)\to\PSh_\emptyset(\Sm^G_S)$ preserves Nisnevich equivalences. 
\end{lemma}

\begin{proof}
	We start with a preliminary observation. Let $\scr C$ be an $\infty$-category with an initial object. If $K$ is any simplicial set and $p\colon K\to\scr C$ is a diagram, then a colimit of $p$ is the same thing as a colimit of an extended diagram $K^{\triangleleft}\to\scr C$ that takes the initial vertex of $K^{\triangleleft}$ to an initial object of $\scr C$, and the simplicial set $K^{\triangleleft}$ is of course weakly contractible.
	
	By Nisnevich separation and smooth base change, we may assume that $B$ is affine and has the $G$-resolution property.
	Let $C\subset\PSh(\Sm^G_Z)^{\Delta^1}$ be the full subcategory consisting of:
	\begin{itemize}
		\item equivalences;
		\item the map $0\to\emptyset$, where $0$ is the empty presheaf and $\emptyset$ is the empty scheme;
		\item the map $K_Q\to X$ for every Nisnevich square $Q$ over $X$ in $\Sm^G_Z$, where $K_Q$ denotes the pushout of $Q$ in $\PSh(\Sm_Z^G)$.
	\end{itemize}
	Note that $C$ contains the initial object of $\PSh(\Sm^G_Z)^{\Delta^1}$, since it is an equivalence. By the preceding observation, the class of Nisnevich equivalences in $\PSh_\emptyset(\Sm^G_Z)$ is the closure of $a_\emptyset(C)$ under 2-out-of-3, pushouts, and weakly contractible colimits. Since $i_*\colon \PSh_\emptyset(\Sm^G_Z)\to\PSh_\emptyset(\Sm^G_S)$ preserves weakly contractible colimits, it suffices to show that $i_*a_\emptyset$ sends elements of $C$ to Nisnevich equivalences. 
	This is obvious for equivalences and for $0\to\emptyset$, since $a_\emptyset(0)=\emptyset$. 
	To conclude, we prove more generally that for every Nisnevich square $Q$ over $X$ in $\Sch^G_Z$, the functor $i_*\colon \PSh(\Sch^G_Z)\to\PSh(\Sch^G_S)$ sends $K_Q\to X$ to a Nisnevich equivalence.
	By universality of colimits in $\PSh(\Sch^G_S)$, it suffices to show that, for every $G$-morphism $p\colon U\to S$ and every map $U_Z\to X$, the projection
	\begin{equation*}\label{eqn:98we}
		i_*K_Q\times_{i_*X}U\to U
	\end{equation*}
	is a Nisnevich equivalence in $\PSh(\Sch^G_S)$. Consider the cartesian square
	\begin{tikzmath}
		\diagram{U_Z & U \\ Z & S\rlap. \\};
		\arrows (11-) edge[c->] node[above]{$i'$} (-12) (11) edge node[left]{$p'$} (21) (21-) edge[c->] node[above]{$i$} (-22) (12) edge node[right]{$p$} (22);
	\end{tikzmath}
	We then have $i'_*p^{\prime *}\simeq p^* i_*$ and, since $\Sch^G_Z$ has pullbacks, the above projection can be identified with $p_\sharp$ of the projection
	\[i'_*(K_{p^{\prime *}Q})\times_{i'_*p^{\prime *}X} U\to U.\]
	Since $p_\sharp$ preserves Nisnevich equivalences, we may as well assume that $U=S$, so that $U_Z=Z$. It is then clear that
	\[i_*K_Q\times_{i_*X}S\simeq i_*K_{Q\times_XZ}.\]
	Thus, we are reduced to proving the following statement:
	\begin{itemize}
		\item[(\textasteriskcentered)] For every Nisnevich square $Q$ over $Z$, $i_*K_Q$ is Nisnevich contractible in $\PSh(\Sch^G_S)$.
	\end{itemize}
	
	Let $Q$ be the square
	\begin{tikzmath}
		\diagram{U\times_ZV & V \\ U & Z\rlap. \\};
		\arrows (11-) edge[c->] (-12) (21-) edge[c->] (-22) (11) edge (21) (12) edge (22);
	\end{tikzmath}
	Let us first prove (\textasteriskcentered) with ``Nisnevich'' replaced by ``Zariski'', so that $V\to Z$ is an open $G$-immersion. Then $S\minus(Z\minus U)$ and $S\minus (Z\minus V)$ form a Zariski cover of $S$. By the Zariski version of Proposition~\ref{prop:Nissep}, it suffices to prove that the restrictions of $i_*K_Q$ to these two open subsets are Zariski contractible. Those restrictions are pushforwards of the restrictions of $K_Q$ to $U$ and $V$, which are clearly final objects in $\PSh(\Sch^G_U)$ and $\PSh(\Sch^G_V)$, respectively. Since $i_*$ preserves final objects, this shows that $i_*K_Q$ is Zariski contractible. Altogether, we have proved that $i_*\colon\PSh_\emptyset(\Sch^G_Z)\to \PSh_\emptyset(\Sch^G_S)$ preserves Zariski equivalences.
	
	We now prove (\textasteriskcentered) in general. By Nisnevich separation, we may assume that $V$ is $G$-quasi-projective, since this is true Nisnevich-locally on $B$. Let $W\into Z$ be a $G$-invariant closed complement of $U$ in $Z$ such that $V\times_ZW\simeq W$. Since $B$ is affine and has the $G$-resolution property and since $V$ is $G$-quasi-projective, we can use Jouanolou's trick (Proposition~\ref{prop:jouanolou}) to obtain a $G$-affine bundle $\tilde V\to V$ where $\tilde V$ is affine. Let $\tilde W=\tilde V\times_VW$. We now consider the following diagram, the dashed parts of which will be described below:
	\begin{tikzmath}
		\def\colsep{4em}
		\diagram{
		W & & \hat W & \\
		& \tilde W & \tilde V & \hat V \\
		& W & Z & S\rlap. \\
		};
		\arrows
		(11-) edge[c->,dashed] (-13)
		(22-) edge[c->] (-23) (23-) edge[c->,dashed] (-24)
		(32-) edge[c->] (-33) (33-) edge[c->] (-34)
		(11) edge[c->,dashed] (22) (13) edge[c->,dashed] (24)
		(11) edge[-,vshift=1pt] (32) edge[-,vshift=-1pt] (32) (13) edge[dashed] node[below left, near end]{$q$} (34)
		(22) edge (32) (23) edge (33) (24) edge[dashed] node[right]{$p$} (34)
		;
	\end{tikzmath}
	First of all, since $\tilde V\to Z$ is smooth, Corollary~\ref{cor:keyLift} implies that there exists an affine $G$-scheme $\hat V$ over $S$ fitting into a cartesian square as above and such that $p$ is smooth at $\tilde V$. Since $B$ and $W$ are affine and $G$ is linearly reductive, the $G$-affine bundle $\tilde W\to W$ admits a $G$-section (Lemma~\ref{lem:linred}), which is a closed $G$-immersion $W\into \tilde W$. We now apply Theorem~\ref{thm:keyLift}: replacing if necessary $\hat V$, $\tilde V$, and $\tilde W$ by affine $G$-equivariant étale neighborhoods of $W$, we obtain an affine $G$-scheme $\hat W$ such that all squares in the above diagram are cartesian, and such that $q\colon\hat W\to S$ is étale at $W$.
	By Lemma~\ref{lem:locus}, there exists a quasi-compact $G$-invariant open neighborhood $\hat W^\circ$ of $W$ on which $q$ is étale.
	 The open $G$-immersion $j\colon S\minus W\into S$ and the étale $G$-morphism $\hat W^\circ\to S$ then form a Nisnevich cover of $S$. By Proposition~\ref{prop:Nissep}, it will therefore suffice to prove that $j^*i_*K_Q$ and $p^*i_*K_Q$ are Nisnevich contractible. We have $j^*i_*K_Q\simeq i'_*K_{Q_U}$ and $p^*i_*K_Q\simeq i^{\prime\prime}_*K_{Q_{\tilde V}}$, where $i'\colon U\into S\minus W$ and $i''\colon\tilde V\into \hat V$. Now, $K_{Q_U}$ is clearly a final object in $\PSh(\Sch^G_U)$. It remains to prove that $i^{\prime\prime}_*K_{Q_{\tilde V}}\in\PSh(\Sch^G_{\hat V})$ is Nisnevich contractible. The Nisnevich square $Q_{\tilde V}$ is isomorphic to the square
	\begin{tikzmath}
		\diagram{\tilde V_U\amalg \tilde Y & \tilde V\amalg \tilde Y \\ \tilde V_U & \tilde V\rlap, \\};
		\arrows (11-) edge[c->] (-12) (21-) edge[c->] (-22) (11) edge (21) (12) edge (22);
	\end{tikzmath}
	where $Y$ is the complement of the diagonal in $V\times_ZV$ and $\tilde Y=Y\times_V\tilde V$. It is then clear that $K_{Q_{\tilde V}}\to \tilde V$ is a Zariski equivalence. As we have already proved that $i^{\prime\prime}_*$ preserves Zariski equivalences, we are done.
\end{proof}

\begin{theorem}\label{thm:exactness}
	Let $i\colon Z\into S$ be a closed $G$-immersion. Then $i_*\colon \PSh_\emptyset(\Sm^G_Z)\to\PSh_\emptyset(\Sm^G_S)$ preserves motivic equivalences.
\end{theorem}

\begin{proof}
	By Proposition~\ref{prop:affinecover} and smooth base change, we may assume that $B$ and $S$ are affine and that $S$ has the $G$-resolution property.
	The observation at the beginning of the proof of Lemma~\ref{lem:closedexact} shows that the class of motivic equivalences in $\PSh_\emptyset(\Sm^G_Z)$ is the closure under 2-out-of-3, pushouts, and weakly contractible colimits of the union of the classes of homotopy equivalences and of Nisnevich equivalences.
	The result now follows from Lemmas \ref{lem:affinePF} and~\ref{lem:closedexact}.
\end{proof}

\begin{corollary}\label{cor:wccolim}
	Let $i\colon Z\into S$ be a closed $G$-immersion. Then $i_*\colon\H^G(Z)\to\H^G(S)$ preserves weakly contractible colimits. 
\end{corollary}

\begin{remark}
	Unless $S\minus Z$ is empty, $i_*\colon \H^G(Z)\to\H^G(S)$ does not preserve the initial object. This inconvenience disappears when we pass to the $\infty$-category of \emph{pointed} motivic $G$-spaces. We will see some remarkable consequences of this fact in \S\ref{sub:pointed}.
\end{remark}

\begin{remark}\label{rmk:finitesifted}
	Suppose that $G$ is finite locally free. Then, if $f\colon Y\to X$ is a finite $G$-morphism, the functor $f_*\colon \PSh_\Sigma(\Sm_Y^G)\to \PSh_\Sigma(\Sm_X^G)$ preserves Nisnevich equivalences.
	 Modulo a noetherian approximation argument, this can be proved exactly as in \cite[\S3, Proposition 1.27]{MV}, using the fact that, when $B$ is noetherian of finite Krull dimension, the family of functors 
	\[\PSh(\Sm_S^G)\to\scr S,\quad F\mapsto F(X\times_{X/G}(X/G)^h_x)\quad (X\in\Sm^G_S\text{ and }x\in X/G),\]
	which is well-defined by Lemma~\ref{lem:Gquasiaffine}, detects Nisnevich equivalences \cite[Proposition 13]{DeligneNote}. It follows that $f_*\colon \H^G(Y)\to\H^G(X)$ preserves sifted colimits. We do not know if these facts hold for more general $G$.
\end{remark}

\subsection{Gluing}
\label{sub:gluing}

Let $Z\into S$ be a closed $G$-immersion with open complement $U\subset S$. If $X\in\Sch^G_S$ and $t\colon Z\into X_Z$ is a $G$-equivariant section of the projection $X_Z\to Z$, we define a presheaf of sets $\Phi_S(X,t)$ on $\Sch^G_S$ by:
\[\Phi_S(X,t)(Y)=\begin{cases} \Hom_S(Y,X)\times_{\Hom_Z(Y_Z,X_Z)}* & \text{if $Y_Z\neq\emptyset$,} \\ * & \text{if $Y_Z=\emptyset$,}\end{cases}\]
where the map $*\to \Hom_Z(Y_Z,X_Z)$ hits the composition $Y_Z\to Z\stackrel t\into X_Z$. More succinctly,
\begin{equation*}\label{eqn:Phi}
	\Phi_S(X,t)=(X\coprod_{X_U}U)\times_{i_*X_Z}S.
\end{equation*}
Note that $\Phi_S(X,t)$ is functorial in the pair $(X,t)$. Moreover, if $f\colon T\to S$ is a $G$-morphism, there is a natural isomorphism
\begin{equation*}\label{eqn:PhiPB}
f^*\Phi_S(X,t)\simeq \Phi_T(X_T,t_T)
\end{equation*}
in $\PSh(\Sch^G_T)$.

\begin{lemma}\label{lem:gluing1}
	Let $p\colon X'\to X$ be a $G$-morphism in $\Sch^G_S$ and let $t\colon Z\into X_Z$ and $t'\colon Z\into X'_Z$ be closed $G$-sections such that $pt'=t$.
Suppose that $p$ is étale at each point of $t'(Z)$.
Then $\Phi_S(p)\colon\Phi_S(X',t')\to\Phi_S(X,t)$ is a Nisnevich equivalence in $\PSh(\Sch^G_S)$.
\end{lemma}

\begin{proof}
	By Lemmas~\ref{lem:locus} and \ref{lem:cartesian} and 2-out-of-3, we can assume that $p$ is étale and that the square
	\begin{tikzmath}
		\diagram{Z & X_Z' \\ Z & X_Z \\};
		\arrows (11-) edge[c->] node[above]{$t'$} (-12) (11) edge[-,vshift=1pt] (21) edge[-,vshift=-1pt] (21) (21-) edge[c->] node[above]{$t$} (-22) (12) edge node[right]{$p$} (22);
	\end{tikzmath}
	is cartesian. We must show that $L_\Nis\Phi_S(p)$ is an equivalence.
	Since it is a map between $0$-truncated objects in an $\infty$-topos, it suffices to show that it is $1$-connective, \ie, that both it and its diagonal are effective epimorphisms.

	Let $f\colon Y\to X$ be an element of $\Phi_S(X,t)(Y)$ and let $Y'=Y\times_X{X'}$. The pair of $G$-morphisms $\pi_1\colon Y'\to Y$ and $Y\times_X(X\minus t(Z))\into Y$ form a Nisnevich cover of $Y$, and we claim that $f$ lifts to $\Phi_S(X',t')$ on that cover. The $G$-morphism $\pi_2\colon Y'\to X'$ defines an element of $\Phi_S(X',t')(Y')$ lifting $f\pi_1$. By definition of $\Phi_S(X,t)$, the open $G$-immersion $Y_U\into Y\times_X(X\minus t(Z))$ is an isomorphism. Thus, $f$ lifts trivially on $Y\times_X(X\minus t(Z))$. This shows that $L_\Nis\Phi_S(p)$ is an effective epimorphism.

	Let $f,g\colon Y\to X'$ be elements in $\Phi_S(X',t')(Y)$ that become equal in $\Phi_S(X,t)(Y)$. Then $f$ and $g$ induce two $G$-sections of the étale $G$-morphism $X'\times_XY\to Y$. Their equalizer $V$ is a quasi-compact $G$-invariant open subscheme of $Y$. By definition of $\Phi_S(X',t')$, $V$ contains $Y_Z$. Thus, $V$ and $Y_U$ form a Zariski cover of $Y$ on which $f$ and $g$ agree. This shows that the diagonal of $L_\Nis\Phi_S(p)$ is an effective epimorphism.
\end{proof}

\begin{lemma}\label{lem:gluing2}
	Let $\scr E$ be a locally free $G$-module on $S$ and let $t\colon S\to \V(\scr E)$ be the zero section. Then $\Phi_S(\V(\scr E),t_Z)\to S$ is a homotopy equivalence in $\PSh(\Sch^G_S)$.
\end{lemma}

\begin{proof}
	The map
	\[\A^1\times\Phi_S(\V(\scr E),t_Z)\to \Phi_S(\V(\scr E),t_Z),\quad (a,f)\mapsto af,\]
	is an $\A^1$-homotopy between the identity and the composition
	\[\Phi_S(\V(\scr E),t_Z)\to S\stackrel t\to \Phi_S(\V(\scr E),t_Z).\qedhere\]
\end{proof}

\begin{theorem}[Gluing]
	\label{thm:gluing}
	Let $i\colon Z\into S$ be a closed $G$-immersion with open complement $j\colon U\into S$. Then, for every $F\in\H^G(S)$, the square
	\begin{tikzmath}
		\def\colsep{2em}
		\diagram{j_\sharp j^* F & F \\ U & i_*i^* F \\};
		\arrows (11-) edge node[above]{$\epsilon$} (-12) (11) edge node[left]{$!$} (21) (21-) edge node[below]{$!$} (-22) (12) edge node[right]{$\eta$} (22);
	\end{tikzmath}
	is cocartesian, where $!$ denotes a unique map.
\end{theorem}

\begin{proof}
	By Nisnevich separation and smooth base change, we can assume that $S$ is separated. All functors of $F$ in this square preserve weakly contractible colimits: this is obvious for the top left, top right, and bottom left corner, and it follows from Corollary~\ref{cor:wccolim} for the bottom right corner. By Proposition~\ref{prop:generators} (1) and our assumption that $S$ is separated, it suffices to prove the theorem when $F$ is the motivic localization of a scheme $X\in\Sm^G_S$ that is affine over $S$. By Theorem~\ref{thm:exactness}, $i_*i^*F$ is then the motivic localization of $i_*X_Z\in\PSh(\Sm^G_S)$. It therefore suffices to show that the canonical map
	\[X\coprod_{X_U} U\to i_*X_Z\]
	in $\PSh(\Sm^G_S)$ is a motivic equivalence. 
	We prove more generally that it is a motivic equivalence in $\PSh(\Sch^G_S)$.
	By universality of colimits, it suffices to show that for every $Y\in\Sch^G_S$ and every map $t\colon Y_Z\to X_Z$ in $\Sch^G_Z$, the projection
	\[(X\coprod_{X_U} U)\times_{i_*X_Z}Y\to Y\]
	is a motivic equivalence. 
	Using that $\Sch^G_Y$ has pullbacks as in the proofs of Lemmas \ref{lem:affinePF} and~\ref{lem:closedexact}, we may assume without loss of generality that $Y=S$, so that $Y_Z=Z$. We are thus reduced to proving the following statement:
	\begin{itemize}
		\item[(\textasteriskcentered)] For every $X\in\Sm^G_S$ that is affine over $S$ and every $G$-section $t\colon Z\into X_Z$, $\Phi_S(X,t)$ is motivically contractible in $\PSh(\Sch^G_S)$.
	\end{itemize}
	
	By Proposition~\ref{prop:affinecover}, it suffices to prove (\textasteriskcentered) when $S$ is small and affine. Moreover, by Lemma~\ref{lem:gluing1}, we are free to replace $X$ by any affine $G$-equivariant étale neighborhood of $t(Z)$.
	By Theorem~\ref{thm:keyLift}, we can thus find a commutative diagram of $G$-schemes
	\begin{tikzmath}
		\def\colsep{4em}
		\diagram{Z & V & \\ & X_Z & X \\ & Z & S\rlap, \\};
		\arrows (11-) edge[c->,dashed] (-12) (12) edge[c->,dashed] (23)
		(11) edge[c->] node[above right]{$t$} (22) (22-) edge[c->] (-23)
		(22) edge (32) (23) edge (33) (32-) edge[c->] node[below]{$i$} (-33)
		(11) edge[-,vshift=1pt] (32) edge[-,vshift=-1pt] (32) (12) edge[dashed] node[below left, near end]{$p$} (33);
	\end{tikzmath}
	where $V\into X$ is a closed $G$-immersion, $p\colon V\to S$ is étale at $Z$, and all squares are cartesian. 
	By Lemma~\ref{lem:locus}, there exists a quasi-compact $G$-invariant open neighborhood $V^\circ$ of $Z$ on which $p$ is étale.
	Denote by $v\colon V^\circ\into V$ the inclusion.
	The open $G$-immersion $j\colon U\into S$ and the étale $G$-morphism $pv\colon V^\circ\to S$ then form a Nisnevich cover of $S$. By Nisnevich separation, it will suffice to show that $j^*\Phi_S(X,t)$ and $v^*p^*\Phi_S(X,t)$ are motivically contractible. Since $Z_U=\emptyset$, $j^*\Phi_S(X,t)\simeq \Phi_{U}(X_U,t_U)$ is a final object in $\PSh(\Sch^G_U)$. 
	Since $v$ and $pv$ are smooth, we have $v^*p^*\Phi_S(X,t)\simeq v^*\Phi_V(X_V,t_V)$, so it suffices to show that $\Phi_V(X_V,t_V)$ is motivically contractible in $\PSh(\Sch^G_V)$.
	But by construction of $V$, there exists a $G$-section $V\into X_V$ extending $t_V\colon Z\to X_V$. Thus, it remains to prove (\textasteriskcentered) when $S$ is small and affine and $X\to S$ has a $G$-section $s$ extending $t$. By Proposition~\ref{prop:linearization}, there exists a $G$-morphism $h\colon X\to \V(\scr N_s)$ that is étale at $s(S)$ and such that $hs=z$, where $z$ is the zero section.
 By Lemma~\ref{lem:gluing1}, $h$ induces a motivic equivalence $\Phi_S(X,t)\to \Phi_S(\bb V(\scr N_s),z_Z)$. Finally, $\Phi_S(\V(\scr N_s),z_Z)$ is motivically contractible in $\PSh(\Sch^G_S)$ by Lemma~\ref{lem:gluing2}.
\end{proof}

\begin{corollary}\label{cor:fullyfaithful}
	Let $i\colon Z\into S$ be a closed $G$-immersion. Then the functor $i_*\colon\H^G(Z)\to \H^G(S)$ is fully faithful.
\end{corollary}

\begin{proof}
	Let $j\colon U\into S$ be the complementary open $G$-immersion. It is clear that $i^*j_\sharp\colon \H^G(U)\to \H^G(Z)$ is the constant functor with value the initial object. Its right adjoint $j^*i_*$ is thus the constant functor with value the terminal object. Let $F\in\H^G(S)$. Applying Theorem~\ref{thm:gluing} to $i_*F$, we obtain a cocartesian square
	\begin{tikzmath}
		\def\colsep{2em}
		\diagram{U & i_*F \\ U & i_*i^*i_* F\rlap, \\};
		\arrows (11-) edge (-12) (11) edge[-,vshift=1pt] (21) edge[-,vshift=-1pt] (21) (21-) edge (-22) (12) edge node[right]{$\eta$} (22);
	\end{tikzmath}
	in $\H^G(S)$, showing that $\eta i_*$ is an equivalence. It follows from a triangle identity that $i_*\epsilon$ is an equivalence. It remains to show that $i_*$ is conservative. By Proposition~\ref{prop:affinecover} and smooth base change, we can assume that $S$ is small and affine. Let $f$ be a morphism in $\H^G(Z)$ such that $i_*(f)$ is an equivalence, and let us show that $f$ is an equivalence. By Proposition~\ref{prop:generators} (1), it suffices to show that $f$ is an equivalence on every affine scheme $X\in\Sm^G_Z$. This follows at once from Corollary~\ref{cor:keyLift} and Lemma~\ref{lem:locus}.
\end{proof}

\begin{corollary}[Closed base change]
	\label{cor:closedBC}
	Let
	\begin{tikzmath}
		\diagram{W & Y \\ Z & X \\};
		\arrows (11-) edge[c->] node[above]{$k$} (-12) (11) edge node[left]{$g$} (21) (21-) edge[c->] node[below]{$i$} (-22) (12) edge node[right]{$f$} (22);
	\end{tikzmath}
	be a cartesian square of $G$-schemes where $i$ is a closed $G$-immersion. Then the exchange transformation
	\begin{equation*}
		\Ex^*_*\colon f^*i_*\to k_*g^*\colon \H^G(Z)\to\H^G(Y)
	\end{equation*}
	is an equivalence.
\end{corollary}

\begin{proof}
	By Corollary~\ref{cor:fullyfaithful}, it suffices to show that $\Ex_*^*i^*$ is an equivalence.
	This follows easily from Theorem~\ref{thm:gluing} and smooth base change.
\end{proof}

\section{Ambidexterity for smooth projective morphisms}
\label{sec:duality}

\subsection{Pointed equivariant motivic spaces}
\label{sub:pointed}

Let $S$ be a $G$-scheme. A \emph{pointed motivic $G$-space} over $S$ is a motivic $G$-space $X$ over $S$ equipped with a global section $S\to X$. We denote by $\H^G_\pt(S)$ the $\infty$-category of pointed motivic $G$-spaces, \ie, the undercategory $\H^G(S)_{S/}$. We denote by $(\ph)_+\colon \H^G(S)\to\H^G_\pt(S)$ the left adjoint to the forgetful functor.

By \cite[Proposition 4.8.2.11]{HA}, the cartesian symmetric monoidal structure on $\H^G(S)$ extends uniquely to a symmetric monoidal structure on $\H^G_\pt(S)$ that is compatible with colimits. Its tensor product will be denoted by $\tens$, and its unit by $\1_S$. 

If $f\colon T\to S$ is a $G$-morphism, both $f_*\colon \H^G(T)\to \H^G(S)$ and its left adjoint $f^*$ preserve the final object, and hence they lift to an adjunction
\[
f^* : \H_\pt^G(S) \rightleftarrows \H_\pt^G(T): f_*.
\]
If moreover $f$ is smooth, then $f^*$ preserves limits and we therefore have an adjunction
\[
f_\sharp : \H_\pt^G(T) \rightleftarrows \H_\pt^G(S) : f^*.
\]
The left adjoint $f_\sharp$ is characterized by $f_\sharp (X_+) \simeq (f_\sharp X)_+$. The following pointed version of the smooth projection formula follows immediately from Proposition~\ref{prop:projection}: for $B\in\H_\pt^G(S)$ and $C\in\H_\pt^G(T)$, the canonical map
\[
f_\sharp(f^*B \tens C) \to B\tens f_\sharp C
\]
is an equivalence. All the other results of \S\ref{sub:easy} have obvious pointed analogs.

\begin{proposition}\label{obs:support}
	Let $i\colon Z\into S$ be a closed $G$-immersion. Then the functor $i_*\colon \H^G_\pt(Z)\to\H^G_\pt(S)$ preserves colimits.
\end{proposition}

\begin{proof}
	One deduces immediately from Corollary~\ref{cor:wccolim} that $i_*$ preserves weakly contractible colimits. It also trivially preserves the initial object, and hence it preserves all colimits.
\end{proof}

It follows from Proposition~\ref{obs:support} that, for $i\colon Z\into S$ a closed $G$-immersion, $i_*$ has a right adjoint
\[i^!\colon \H^G_\pt(S)\to \H^G_\pt(Z).\]

\begin{proposition}[Pointed gluing]
	\label{prop:pointedgluing}
	Let $i\colon Z\into S$ be a closed $G$-immersion with open complement $j\colon U\into S$. For every $X\in\H^G_\pt(S)$,
	\[
	j_\sharp j^* X \to X \to i_*i^* X
	\]
	is a cofiber sequence, and
	\[
	i_* i^! X \to X \to j_*j^* X
	\]
	is a fiber sequence.
\end{proposition}

\begin{proof}
	The second statement follows from the first one by adjunction.
	Denote by $u\colon \H_\pt^G(S)\to\H^G(S)$ the forgetful functor.
	We consider the following diagram in $\H^G(S)$:
	\begin{tikzmath}
		\def\colsep{1.8em}
		\def\rowsep{2em}
		\diagram{j_\sharp j^* u(X) &j_\sharp j^* u(X)\coprod_US & u(X) \\
		U & S & i_* i^*u(X)\rlap. \\};
		\arrows (11-) edge (-12) (12-) edge (-13)
		(21-) edge node[above]{$j$} (-22) (22-) edge (-23)
		(11) edge (21) ([yshift=1em] 12.south) edge (22) (13) edge (23);
	\end{tikzmath}
	The outside square is cocartesian by Theorem~\ref{thm:gluing}, and the first square is formally cocartesian. Thus, the second square is also cocartesian. But the second square is the image by $u$ of the given sequence. It remains to observe that $u$ reflects weakly contractible colimits.
\end{proof}

\begin{corollary}[Smooth–closed base change]
	\label{cor:closedsupp}
	Let
	\begin{tikzmath}
		\diagram{W & T \\ Z & S \\};
		\arrows (11-) edge[c->] node[above]{$t$} (-12) (11) edge node[left]{$q$} (21) (21-) edge[c->] node[below]{$s$} (-22) (12) edge node[right]{$p$} (22);
	\end{tikzmath}
	be a cartesian square of $G$-schemes, where $p$ is smooth and $s$ is a closed immersion. Then the exchange transformations \begin{gather*}
	\Ex_{\sharp *}\colon p_\sharp t_*\to s_* q_\sharp\colon\H^G_\pt(W)\to \H^G_\pt(S),\\
	\Ex^{*!}\colon q^*s^!\to t^!p^*\colon \H^G_\pt(S)\to\H^G_\pt(W)
\end{gather*}
 are equivalences.
\end{corollary}

\begin{proof}
	The second transformation is the mate of the first one, so it suffices to show that the first transformation is an equivalence. 
	Since $t_*$ is fully faithful, it suffices to show that $\Ex_{\sharp *}t^*$ is an equivalence.
	This follows easily from Proposition~\ref{prop:pointedgluing} and smooth base change.
\end{proof}

\begin{corollary}[Closed projection formula]
	\label{cor:closedproj}
	Let $i\colon Z\into S$ be a closed $G$-immersion and let $A\in\H^G_\pt(S)$. For every $B\in\H^G_\pt(Z)$ and every $C\in \H^G_\pt(S)$, the canonical maps
	\begin{gather*}
		A\tens i_*B\to i_*(i^*A\tens B),\\
		\Hom(i^*A, i^!C)\to i^!\Hom(A,C)
	\end{gather*}
	are equivalences.
\end{corollary}

\begin{proof}
	The first equivalence follows at once from Proposition~\ref{prop:pointedgluing} and the smooth projection formula. The second equivalence follows by adjunction.
\end{proof}

\subsection{Spheres, Thom spaces, and purity}

Let $S$ be a $G$-scheme and let $\scr M$ be a locally free $G$-module on $S$. Let $p\colon\V(\scr M)\to S$ be the associated vector bundle with zero section $s$. We denote by
\[\Sigma^{\scr M}: \H^G_\pt(S)\rightleftarrows\H^G_\pt(S):\Omega^{\scr M}\]
the adjunction $p_\sharp s_*\dashv s^!p^*$.
We call $\Sigma^{\scr M}X$ the \emph{$\scr M$-suspension} of $X$ and $\Omega^{\scr M}X$ the \emph{$\scr M$-loop space} of $X$.
The pointed $G$-space $\Sigma^{\scr M}\1_S\in \H_\pt^G(S)$
is called the \emph{$\scr M$-sphere} and is also denoted by $\s^{\scr M}$.

Let $f\colon T\to S$ be a $G$-morphism and let $\scr M$ be a locally free $G$-module on $S$. By smooth and closed base change, we obtain equivalences
\[f^*\Sigma^{\scr M} \simeq \Sigma^{f^*(\scr M)}f^*\quad\text{and}\quad \Omega^{\scr M}f_*\simeq f_*\Omega^{f^*(\scr M)}.\]
If moreover $f$ is smooth (\resp{} a closed immersion), we also have an equivalence
\[
\Sigma^{\scr M}f_\sharp \simeq f_\sharp \Sigma^{f^*(\scr M)} \quad(\text{\resp{} } f^!\Omega^{\scr M} \simeq \Omega^{f^*(\scr M)}f^!)
\]
by Corollary~\ref{cor:closedsupp}.

By the smooth and closed projection formulas, we have canonical equivalences
\[\Sigma^\scr M\simeq \s^{\scr M}\tens (\ph)\quad\text{and}\quad \Omega^{\scr M}\simeq \Hom(\s^{\scr M},\ph).\]
By Proposition~\ref{prop:pointedgluing}, we have
\[
\s^{\scr M}\simeq \frac{\V(\scr M)}{\V(\scr M)\minus S},
\]
\ie, there is a cofiber sequence
\[(\V(\scr M)\minus S)_+\to \V(\scr M)_+\to \s^{\scr M}\]
in $\H_\pt^G(S)$.
It follows from this description that the assignments $\scr M\mapsto \Sigma^{\scr M}$ and $\scr M\mapsto \Omega^{\scr M}$ are functors on the category of locally free $G$-modules and epimorphisms (this functoriality can also be described more directly using exchange transformations).

Let 
\begin{equation*}
	\label{eqn:SESmod}
	0\to \scr N\to\scr M\to\scr P\to 0
\end{equation*}
be a short exact sequence of locally free $G$-modules on a $G$-scheme $S$. In \S\ref{sub:scp}, we defined a canonical equivalence
\[
\Psi\colon \s^{\scr M}\simeq \s^{\scr N}\tens \s^{\scr P}. 
\]
We will also denote by $\Psi$ the induced equivalence of functors $\Sigma^{\scr M}\simeq \Sigma^{\scr N}\Sigma^{\scr P}$.

\begin{remark}
	Given a short exact sequence as above, we can form the diagram
\begin{tikzmath}
	\diagram{
	S & & \\
	\V(\scr P) & \V(\scr M) & \\
	S & \V(\scr N) & S\rlap, \\
	};
	\arrows (11) edge node[left]{$t$} (21) edge node[above right]{$s$} (22)
	(22) edge (32) edge node[above right]{$p$} (33)
	(21-) edge node[above]{$a$} (-22) (21) edge node[left]{$q$} (31)
	(31-) edge node[below]{$u$} (-32) (22) edge node[right]{$b$} (32)
	(32-) edge node[below]{$r$} (-33);
\end{tikzmath}
where $p,q,r$ are the structure maps of the associated vector bundles and $s,t,u$ are their respective zero sections.
By Corollary~\ref{cor:closedsupp}, the transformation $\Ex_{\sharp *}\colon b_\sharp a_*\to u_*q_\sharp$ is an equivalence.
Using Remark~\ref{rmk:splitting} and the fact that any short exact sequence of locally free $G$-modules splits when pulled back to an appropriate $G$-affine bundle, it is easy to show that the equivalence $\Psi$ agrees with the transformation
\[\Sigma^{\scr M}=p_\sharp s_*\simeq r_\sharp b_\sharp a_* t_*\xrightarrow{\Ex_{\sharp *}} r_\sharp u_*q_\sharp t_*=\Sigma^{\scr N}\Sigma^{\scr P}.\]
We will not need this alternative description.
\end{remark}

If $\scr M$ is a locally free $G$-module on $X\in\Sm^G_S$, we let
\[\Th_{X}(\scr M)=p_\sharp \s^{\scr M}\in\H^G_\pt(S),\]
where $p\colon X\to S$ is the structure map. The pointed $G$-space $\Th_{X}(\scr M)$ is called the \emph{Thom space} of $\scr M$.

We will now recast the purity equivalence $\Pi$ of \S\ref{sub:scp} in a functorial setting.
Let
\begin{tikzequation}\label{eqn:scp}
	\diagram{ Z & X \\ & S \\};
	\arrows (11-) edge[c->] node[above]{$s$} (-12) (11) edge node[below left]{$q$} (22) (12) edge node[right]{$p$} (22); 
\end{tikzequation}
be a commutative triangle of $G$-schemes where $p$ and $q$ are smooth and $s$ is a closed immersion. To such a smooth closed pair $(X,Z)$ over $S$ we can associate the functor $p_\sharp s_*$, and to a morphism of smooth closed pairs $f\colon (X',Z')\to (X,Z)$ we can associate the natural transformation
\[
\psi(f)\colon p'_\sharp s'_* f^* \to p_\sharp s_*
\]
(with obvious notation),
which is adjoint to the composition
\[
s'_*f^*\xleftarrow{\Ex^*_*} f^*s_* \xrightarrow\eta f^*p^{*}p_\sharp s_* \simeq p^{\prime *}p_\sharp s_*.
\]
By Proposition~\ref{prop:pointedgluing}, $p_\sharp s_*\1_Z$ is the pointed space $X/(X\minus Z)$, and $\psi(f)_{\1_Z}$ is simply the map ${X'}/({X'}\minus {Z'})\to X/(X\minus Z)$ induced by $f$.

Consider the deformation to the normal bundle:
	\[
	\begin{tikzpicture}
		\diagram{Z & X \\ & S \\};
		\arrows (11-) edge[c->] node[above]{$s$} (-12) (11) edge (22) (12) edge node[right]{$p$} (22);
	\end{tikzpicture}
	\quad
	\begin{tikzpicture}
		\draw[c->,font=\scriptsize] (0,0) -- node[above]{$i_1$} (1,0);
	\end{tikzpicture}
	\quad
	\begin{tikzpicture}
		\diagram{Z\times\A^1 & D_ZX \\ & S\times\A^1 \\};
		\arrows (11-) edge[c->] node[above]{$\hat s$} (-12) (11) edge (22) (12) edge node[right]{$\hat p$} (22);
	\end{tikzpicture}
	\quad
	\begin{tikzpicture}
		\draw[left hook->,font=\scriptsize] (1,0) -- node[above]{$i_0$} (0,0);
	\end{tikzpicture}
	\quad
	\begin{tikzpicture}
		\diagram{Z & \bb V(\scr N_s) \\ & S\rlap, \\};
		\arrows (11-) edge[c->] node[above]{$s_0$} (-12) (11) edge (22) (12) edge node[right]{$p_0$} (22);
	\end{tikzpicture}
	\]
	and denote by $r\colon S\times\A^1\to S$ and $r\colon Z\times\A^1\to Z$ the projections.
We then obtain a zig-zag of natural transformations
\[
p_\sharp s_*i_1^* \xrightarrow{\psi(i_1)} r_\sharp\hat p_\sharp \hat s_* \xleftarrow{\psi(i_0)} q_\sharp\Sigma^{\scr N_s} i_0^*.
\]
Precomposing with $r^*$, we obtain the zig-zag
\[
p_\sharp s_* \to r_\sharp\hat p_\sharp \hat s_* r^* \from q_\sharp\Sigma^{\scr N_s}.
\]

\begin{proposition}[Purity]
	\label{prop:purity}
	With the above notation, the four natural transformations
	\begin{gather*}
		p_\sharp s_* \to r_\sharp \hat p_\sharp \hat s_* r^* \from q_\sharp\Sigma^{\scr N_s},\\
		s^!p^* \from r_*\hat s^!\hat p^* r^* \to \Omega^{\scr N_s}q^*
	\end{gather*}
	are equivalences.
\end{proposition}
	
\begin{proof}
	The second zig-zag is the mate of the first one, so it suffices to show that the first zig-zag is an equivalence.
		Since $s_*$ is fully faithful and all functors involved preserve colimits, it suffices to show that it is an equivalence on $s^*(Y_+)$ for $Y\in\Sm_X^G$. But in that case, by Proposition~\ref{prop:pointedgluing}, it can be identified with the zig-zag
	\[
	\frac{X\times_XY}{(X\minus Z)\times_XY} \to \frac{D_ZX\times_XY}{(D_ZX\minus (Z\times\A^1))\times_XY}\from \frac{N_ZX\times_XY}{(N_ZX\minus Z)\times_XY},
	\]
	which is an equivalence in $\H_\pt^G(S)$ by Theorem~\ref{thm:wexcisive}.
\end{proof}
	
With every triangle~\eqref{eqn:scp} are therefore associated canonical equivalences
\[
	\Pi_s\colon p_\sharp s_*\simeq q_\sharp \Sigma^{\scr N_s}\quad\text{and}\quad\Pi_s\colon s^!p^*\simeq \Omega^{\scr N_s}q^*,
\]
which will be called the \emph{purity equivalences}, generalizing their namesake from \S\ref{sub:scp}.
The two basic properties of the latter (naturality in the smooth closed pair and compatibility with smooth closed triples) can be promoted to properties of these natural transformations. 
We do not make them explicit since in the sequel we will only use the simpler form of these results discussed in \S\ref{sub:scp}.

\subsection{The Pontryagin–Thom collapse}
\label{sub:pt}

Let $S$ be a $G$-scheme and let $f\colon X\to S$ be a smooth $G$-projective morphism. Fix a locally free $G$-module $\scr E$ on $S$, everywhere of rank $\geq 1$, and a closed $G$-immersion $i\colon X\into \P(\scr E)$ over $S$.
The goal of this subsection is to define the \emph{Pontryagin–Thom collapse map}
\[\eta=\eta_{(X,i)}\colon \s^{\scr M}\to \Th_X(\scr N)\]
in $\H^G_\pt(S)$. Here, $\scr M=\scr M_{(X,i)}$ and $\scr N=\scr N_{(X,i)}$ are locally free $G$-modules on $S$ and $X$, related by an explicit equivalence
\[\alpha=\alpha_{(X,i)}\colon \s^{\Omega_f}\tens \s^{\scr N}\simeq \s^{f^*(\scr M)}\]
in $\H_\pt^G(X)$.
When $G$ is trivial and $S$ is the spectrum of a field, $\eta_{(X,i)}$ coincides with the map defined by Voevodsky in \cite[Theorem 2.11]{VV}, up to some $\scr O_S^n$-suspension.

We warn the reader that the construction of $\eta$ is rather involved, so we start with an informal discussion. 
If $X$ is instead a smooth compact manifold over $S=*$, the classical Pontryagin–Thom collapse map is defined as follows. We first choose an embedding of $X$ into a ball $B$ inside some large Euclidean space $E$. The collapse map $\eta\colon E/(E\minus B)\to E/(E\minus X)$ is then a map from a sphere $\s^E$ to the Thom space of the normal bundle of $X$ in $E$ (by the tubular neighborhood theorem). To replicate this construction for $X$ a smooth $G$-projective $S$-scheme, one would need an embedding of $X$ into a projective bundle that misses a hyperplane, which is of course impossible (unless $X$ is finite). The key idea is that we can nevertheless find such an embedding ``up to homotopy'': there exists a zig-zag
\[X\stackrel p\longleftarrow \tilde X\stackrel s\longinto \tilde\P(\scr F)\stackrel q\longrightarrow \P(\scr F)\]
in $\Sm_S^G$, where $\scr F$ is a locally free $G$-module on $S$, $s$ is a closed $G$-immersion, and $p$ and $q$ are compositions of $G$-affine bundles, such that the following two properties hold. 
The first is that the image of $\tilde X$ in $\P(\scr F)$ misses a $G$-invariant hyperplane $\P(\scr G)$; this allows us to define a collapse map 
\begin{equation}
	\label{eqn:PT1}
	\s^{\scr G}\simeq \frac{\P(\scr F)}{\P(\scr G)}\to\Th_{\tilde X}(\scr N_s)
\end{equation}
 in $\H^G_\pt(S)$ (by the purity theorem). The second is that
there exists a locally free $G$-module $\scr N$ on $X$ such that $\s^{p^*(\scr N)}\simeq \s^{\scr N_s}$,
 at least after tensoring with a sphere defined over $S$. This gives us a ``stable'' equivalence
\begin{equation}
	\label{eqn:PT2}
\Th_{\tilde X}(\scr N_s) \simeq \Th_X(\scr N).
\end{equation}
The Pontryagin–Thom collapse map $\eta$ is then the composition of \eqref{eqn:PT1} and \eqref{eqn:PT2}.
Some further care must be taken in the actual construction of $\eta$ below,
and although it does fit this basic sketch, this will not be made explicit.

We first give the construction of $\eta$ and $\alpha$ in the special case where $i\colon X\into\P(\scr E)$ is the identity. 
Recall from \S\ref{sub:jouanolou} the short exact sequence
\begin{equation}\label{eqn:Q}
	0\to\scr O_S\to\scr E\tens\scr E^\vee\to \scr Q\to 0,
\end{equation}
where the first map is the coevaluation and $\scr Q$ is locally free.
We let
\[ P=\P(\scr E\otimes\scr E^\vee),\quad H=\P(\scr Q),\quad V = P\minus H.\]
The map $v\colon V\to S$ is the $G$-affine bundle that universally splits the short exact sequence~\eqref{eqn:Q}. By Proposition~\ref{prop:hi}, the functor
\[
v^*\colon \H_\pt^G(S)\to \H_\pt^G(V)
\]
is fully faithful.
Thus, implicitly replacing $S$ by $V$, we may assume given a splitting
\[
\scr E\otimes\scr E^\vee\simeq \scr O_S\oplus \scr Q
\]
of~\eqref{eqn:Q},
so that $V$ is the $G$-vector bundle $\V(\scr Q)$ and $P$ is its projective completion. In particular, we have the usual zig-zag
\[
\s^{\scr Q}=\frac{V}{V\minus 0}\to \frac{P}{P\minus 0}\from \frac{P}{H},
\]
where the first map is a Zariski equivalence and the second is a homotopy equivalence ($H\into P\minus 0$ being the zero section of the tautological line bundle on $\P(\scr Q)$).

Let $X^\vee=\P(\scr E^\vee)$ and let $\segre\colon X\times X^\vee\into P$ be the Segre embedding (here and in what follows, products are taken over $S$).
Consider the cartesian squares
\begin{tikzmath}
	\diagram{\tilde X & X\times X^\vee & K \\
	V & P & H\rlap. \\};
	\arrows (11-) edge[c->] node[above]{$j$} (-12) (-13) edge[left hook->] (12-)
	(21-) edge[c->] (-22) (-23) edge[left hook->] (22-)
	(11) edge[c->] node[right]{$\segre$} (21) (12) edge[c->] node[right]{$\segre$} (22) (13) edge[c->] node[right]{$\segre$} (23);
\end{tikzmath}
We denote by $\pi\colon \tilde X\to X$ and $\rho\colon \tilde X\to X^\vee$ the canonical projections; recall from \S\ref{sub:jouanolou} that $\pi$ is a $G$-affine bundle.
We also consider the cartesian squares
\begin{tikzmath}
	\diagram{\tilde K & \tilde X\times X^\vee & \tilde X \\
	K & X\times X^\vee & X\rlap, \\};
	\arrows (11-) edge[c->] (-12) (12-) edge (-13)
	(21-) edge[c->] (-22) (22-) edge (-23)
	(11) edge node[right]{$\pi$} (21) (12) edge node[right]{$\pi$} (22) (13) edge node[right]{$\pi$} (23);
\end{tikzmath}
and we let 
\[
s\colon \tilde X\into \tilde X\times X^\vee,\quad x\mapsto (x,\rho x),
\]
be the graph of $\rho$. 
We then have the following zig-zag
\begin{tikzmath}
	\diagram{X & & \\
	\tilde X & \tilde X\times X^\vee & \\
	& X\times X^\vee & P\rlap, \\};
	\arrows
	(21) edge node[left]{$\pi$} (11)
	(21-) edge[c->] node[above]{$s$} (-22) (22) edge node[right]{$\pi$} (32)
	(32-) edge[c->] node[above]{$\segre$} (-33)
	(21) edge[c->] node[below left]{$j$} (32);
\end{tikzmath}
where the vertical maps are $G$-affine bundles and the horizontal maps are closed $G$-immersions.
Note that $s(\tilde X)\cap \tilde K=\emptyset$ since $j(\tilde X)\cap K=\emptyset$.
We therefore obtain the following collapse map in $\H_\pt^G(S)$:
\begin{multline}
	\label{eqn:thom1}
	\s^{\scr Q}\simeq \frac PH \to \frac P{H\cup (P\minus\segre(X\times X^\vee))}\stackrel\Pi\simeq \frac{\Th_{X\times X^\vee}(\scr N_\segre)}{\Th_{K}(\scr N_\segre)}\\
	\stackrel\pi \from \frac{\Th_{\tilde X\times X^\vee}(\pi^*(\scr N_\segre))}{\Th_{\tilde K}(\pi^*(\scr N_\segre))}
	\to \frac{\V(\pi^*(\scr N_\segre))}{\V(\pi^*(\scr N_\segre))\minus s(\tilde X)}\stackrel\Pi\simeq \Th_{\tilde X}(\scr N_s\oplus j^*(\scr N_\segre)).
\end{multline}

Consider the canonical short exact sequence
\begin{equation*}
	\label{eqn:canE}
	0\to \Omega_{X} \to f^*(\scr E)(-1) \to \scr O_{X} \to 0
\end{equation*}
in $\QCoh^G(X)$.
Tensoring it with its dual, we obtain the following diagram of short exact sequences:
\begin{tikzequation}
	\label{eqn:grid}
	\def\colsep{2em}
	\diagram{
	 & 0 & 0 & 0 & \\
	 0 & \Omega_X & f^*(\scr E)(-1) & \scr O_X & 0 \\
	 0 & \Omega_X \tens f^*(\scr E^\vee)(1) & f^*(\scr E\tens \scr E^\vee) & f^*(\scr E^\vee)(1) & 0 \\
	 0 & \Omega_X \tens\Omega_X^\vee & f^*(\scr E)(-1)\tens\Omega_X^\vee & \Omega_X^\vee & 0\rlap. \\
	  & 0 & 0 & 0 & \\
	};
	\arrows
	(21-) edge (-22) (22-) edge (-23) (23-) edge (-24) (24-) edge (-25)
	(31-) edge (-32) (32-) edge (-33) (33-) edge (-34) (34-) edge (-35)
	(41-) edge (-42) (42-) edge (-43) (43-) edge (-44) (44-) edge (-45)
	(12) edge (22) (22) edge (32) (32) edge (42) (42) edge (52)
	(13) edge (23) (23) edge (33) (33) edge (43) (43) edge (53)
	(14) edge (24) (24) edge (34) (34) edge (44) (44) edge (54)
	;
\end{tikzequation}
We define
\[
\scr R=f^*(\scr E)(-1)\tens\Omega_X^\vee. 
\]
The short exact sequences \eqref{eqn:grid} define an equivalence
\begin{equation}\label{eqn:Kpath1}
	\Sigma^{\scr O}\Sigma^{\Omega_X}\Sigma^{\scr R} \simeq \Sigma^{f^*(\scr E\tens\scr E^\vee)}
\end{equation}
of endofunctors of $\H_\pt^G(X)$.

\begin{remark}
	In the sequel, we will never use the \emph{definitions} of $\scr R$ and of the equivalence~\eqref{eqn:Kpath1}. Thus, any choice of a locally free $G$-module $\scr R$ on $X$ together with such an equivalence can be used to define a Pontryagin–Thom collapse map. Different choices will be shown to yield ``stably'' equivalent maps. 
\end{remark}

The Segre embedding $\segre$ induces the short exact sequence
\begin{equation}
	\label{eqn:canE3}
	0\to \scr N_\segre \to \segre^*(\Omega_{P}) \to \Omega_{X\times X^\vee} \to 0.
\end{equation}
With the obvious isomorphisms
\[
	\Omega_V\simeq v^*(\scr Q),\quad \Omega_{\tilde X}\simeq \pi^*(\Omega_X)\oplus\rho^*(\Omega_{X^\vee}),\quad \scr N_s\simeq \rho^*(\Omega_{X^\vee}),
\]
the short exact sequence $j^*\eqref{eqn:canE3}$ becomes
\begin{equation*}
	\label{eqn:qwerty}
	0\to j^*(\scr N_\segre) \to \tilde f^*(\scr Q)\to \pi^*(\Omega_X)\oplus \scr N_s\to 0,
\end{equation*}
where $\tilde f=f\pi\colon\tilde X\to S$ is the structure map.
This induces an equivalence
\begin{equation}\label{eqn:Kpath2}
	\Sigma^{\pi^*(\Omega_X)}\Sigma^{j^*(\scr N_\segre)\oplus \scr N_s}\simeq \Sigma^{\tilde f^*(\scr Q)}
\end{equation}
of endofunctors of $\H_\pt^G(\tilde X)$. Combining \eqref{eqn:Kpath1} and~\eqref{eqn:Kpath2}, we obtain an equivalence
\begin{equation}
	\label{eqn:Kpath3}
	\Sigma^{\scr O}\Sigma^{\pi^*(\Omega_X)}\Sigma^{\pi^*(\scr R)}\simeq \Sigma^{\scr O}\Sigma^{\pi^*(\Omega_X)}\Sigma^{j^*(\scr N_\segre)\oplus \scr N_s}.
\end{equation}
Combining \eqref{eqn:Kpath3} and~\eqref{eqn:Kpath1}, we further obtain an equivalence
\begin{equation}\label{eqn:Kpath4}
	\Sigma^{\tilde f^*(\scr E\tens\scr E^\vee)}\Sigma^{\pi^*(\scr R)}\simeq \Sigma^{\tilde f^*(\scr E\tens\scr E^\vee)}\Sigma^{j^*(\scr N_\segre)\oplus \scr N_s}.
\end{equation}
We let
\[
\scr M=(\scr E\tens\scr E^\vee)\oplus \scr Q,\quad \scr N=f^*(\scr E\tens\scr E^\vee)\oplus\scr R.
\]
From \eqref{eqn:Kpath1} we obtain the equivalence $\alpha=\alpha_{(\P(\scr E),\id)}$ as follows:
\begin{equation}\label{eqn:Kpath5}
\Sigma^{\scr O}\Sigma^{\Omega_X}\Sigma^{\scr R}\simeq \Sigma^{\scr O}\Sigma^{f^*(\scr Q)}
\;\Longrightarrow\;
\underbrace{\Sigma^{\Omega_X}\Sigma^{\scr R}\Sigma^{\scr O}}_{\textstyle\Sigma^{f^*(\scr E\tens\scr E^\vee)}}\Sigma^{\Omega_X}\Sigma^{\scr R}\simeq \underbrace{\Sigma^{\Omega_X}\Sigma^{\scr R}\Sigma^{\scr O}}_{\textstyle\Sigma^{f^*(\scr E\tens\scr E^\vee)}}\Sigma^{f^*(\scr Q)}
\;\Longrightarrow\;
\Sigma^{\Omega_X}\Sigma^{\scr N}\stackrel\alpha\simeq\Sigma^{f^*(\scr M)}.
\end{equation}
Finally, the Pontryagin–Thom collapse map $\eta=\eta_{(\P(\scr E),\id)}$ is defined as the composition
\[
\s^{\scr M}\simeq \Sigma^{\scr E\tens\scr E^\vee}\s^{\scr Q}
\xrightarrow{\eqref{eqn:thom1}}
\Sigma^{\scr E\tens\scr E^\vee}\Th_{\tilde X}(\scr N_s\oplus j^*(\scr N_\segre))
\stackrel{\eqref{eqn:Kpath4}}\simeq
\Th_{\tilde X}(\pi^*(\scr N))\simeq\Th_X(\scr N).
\]

This concludes the construction of $\alpha_{(X,i)}$ and $\eta_{(X,i)}$ when $i$ is the identity.
In general,
let $\scr M_{(X,i)}$ be $\scr M_{(\P(\scr E),\id)}$, and let $\scr N_{(X,i)}$ be the conormal sheaf of the immersion $X\into \V(\scr N_{(\P(\scr E),\id)})$, which is canonically identified with $\scr N_i\oplus i^*(\scr N_{(\P(\scr E),\id)})$. Then the equivalence $\alpha_{(X,i)}$ is given by
\[
\Sigma^{\Omega_X}\Sigma^{\scr N_{(X,i)}}=\Sigma^{\Omega_X}\Sigma^{\scr N_i}\Sigma^{i^*(\scr N_{(\P(\scr E),\id)})} \stackrel\Psi\simeq \Sigma^{i^*(\Omega_{\P(\scr E)})}\Sigma^{i^*(\scr N_{(\P(\scr E),\id)})} \stackrel{\alpha_{(\P(\scr E),\id)}}\simeq \Sigma^{f^*(\scr M)},
\]
and the Pontryagin–Thom collapse map $\eta_{(X,i)}$ is the composition
\[
\s^{\scr M_{(X,i)}}\xrightarrow{\eta_{(\P(\scr E),\id)}} \Th_{\P(\scr E)}(\scr N_{(\P(\scr E),\id)})\to \frac{\V(\scr N_{(\P(\scr E),\id)})}{\V(\scr N_{(\P(\scr E),\id)})\minus i(X)}\stackrel\Pi\simeq \Th_X(\scr N_{(X,i)}).
\]

Now that $\eta=\eta_{(X,i)}$ has been defined, we can upgrade it to a natural transformation
\[
\eta\colon\Sigma^{\scr M}\to f_\sharp\Sigma^{\scr N}f^*
\]
using the projection formulas
\[
\Sigma^{\scr M}\simeq \s^{\scr M}\tens(\ph)\quad\text{and}\quad f_\sharp\Sigma^{\scr N}f^*\simeq \Th_X(\scr N)\tens(\ph).
\]

We denote by
\[
\epsilon\colon f^*f_\sharp \to \Sigma^{\Omega_f}
\]
the composition
\begin{equation}\label{eqn:pure}
f^*f_\sharp\simeq \pi_{2\sharp}\pi_1^*\to \pi_{2\sharp}\delta_*\delta^*\pi_1^*\simeq \pi_{2\sharp}\delta_*\stackrel\Pi\simeq \Sigma^{\scr N_\delta}\stackrel\nu\simeq\Sigma^{\Omega_f},
\end{equation}
where $\pi_{1,2}\colon X\times X\rightrightarrows X$ are the projections, $\delta\colon X\into X\times X$ is the diagonal, and $\nu\colon\scr N_\delta\to\Omega_f$ is the isomorphism sending the class of $x\tens 1-1\tens x$ to $dx$. Note that, unlike $\eta$, $\epsilon$ does not depend on a choice of embedding $i\colon X\into \P(\scr E)$.
For $Y\in\Sm_X^G$, the component of $\epsilon$ at $Y_+$ is given more explicitly by the map
\[
(Y\times X)_+ \to \frac{Y\times X}{(Y\times X)\minus Y}\stackrel\Pi\simeq \Sigma^{\scr N_\delta}Y_+\stackrel\nu\simeq \Sigma^{\Omega_f}Y_+
\]
collapsing the complement of the graph of $Y\to X$, where $Y\times X$ belongs to $\Sm_X^G$ via the second projection. 

\begin{lemma}\label{lem:symmetry}
	Let $S$ be a $G$-scheme and $\scr E$ a locally free $G$-module on $S$. Then the transposition on $\s^{\scr E}\tens \s^{\scr E}$ in $\H_\pt^G(S)$ is homotopic to $\langle -1\rangle\tens\id$, where $\langle -1\rangle\colon \s^{\scr E}\to \s^{\scr E}$ is induced by the linear automorphism $\scr E\to\scr E$, $x\mapsto -x$.
\end{lemma}

\begin{proof}
	The matrices
	\[
	\begin{pmatrix} 0 & 1 \\ 1 & 0 \end{pmatrix}\quad\text{and}\quad\begin{pmatrix} -1 & 0 \\ 0 & 1 \end{pmatrix}
	\]
	are related by elementary transformations and hence are $\A^1$-homotopic in $\mathrm{SL}_2(\mathbb Z)$.
	We conclude using the action of $\mathrm{SL}_2(\mathbb Z)$ on $\V(\scr E\oplus \scr E)/(\V(\scr E\oplus\scr E)\minus 0)$.
\end{proof}

The following theorem is the last nontrivial result in our approach to the formalism of six operations:

\begin{theorem}[Unstable ambidexterity]
	\label{thm:duality}
	Let $\scr E$ be a locally free $G$-module of rank $\geq 1$ on $S$ and let
	\begin{tikzmath}
		\diagram{X & \P(\scr E) \\ & S \\};
		\arrows (11-) edge[c->] node[above]{$i$} (-12) (11) edge node[below left]{$f$} (22) (12) edge (22);
	\end{tikzmath}
	be a commuting triangle of $G$-schemes where $i$ is a closed immersion and $f$ is smooth.
	Then the compositions
	\begin{gather}
		\label{eqn:qwerty1}
	f^*\Sigma^{\scr M}\xrightarrow{\eta} f^*f_\sharp\Sigma^{\scr N}f^* \xrightarrow{\epsilon} \Sigma^{\Omega_f} \Sigma^{\scr N} f^*\stackrel{\alpha}\simeq \Sigma^{f^*(\scr M)}f^*\simeq f^*\Sigma^{\scr M}, \\
	\label{eqn:qwerty2}
	\Sigma^{\scr M}f_\sharp \xrightarrow\eta f_\sharp\Sigma^{\scr N}f^* f_\sharp \xrightarrow\epsilon f_\sharp\Sigma^{\scr N}\Sigma^{\Omega_f}\stackrel{\alpha}\simeq f_\sharp \Sigma^{f^*(\scr M)} \simeq \Sigma^{\scr M}f_\sharp
	\end{gather}
	are the identity.
\end{theorem}

\begin{proof}
	We keep using the notation introduced above in the construction of $\eta$. In addition, if $e\colon Z\into Y$ is a closed $G$-immersion, we will abbreviate $Y\minus e(Z)$ to $e^c$ or $Z^c$. Let $\iota$ be the transformation~\eqref{eqn:qwerty1} evaluated on $\1_S$.
It is not difficult to check that \eqref{eqn:qwerty1} can be identified with $\iota\tens f^*(\ph)$ and \eqref{eqn:qwerty2} with $f_\sharp(\iota\tens(\ph))$.
Thus, it will suffice to show that $\iota$ is the identity in $\H_\pt^G(X)$.
	
	Let us first reduce to the case where $i\colon X\into \P(\scr E)$ is the identity.
	In fact, we claim that $\iota_X=i^*(\iota_{\P(\scr E)})$ in $\H_\pt^G(X)$.
	Recall that $\scr N_{(X,i)}=\scr N_i\oplus i^*(\scr N_{(\P(\scr E),\id)})$. We consider the following diagram in $\H_\pt^G(X)$:
	\begin{tikzmath}
		\def\colsep{4em}
		\diagram{
		\Th_X(f^*(\scr M)) & \frac{\V(\scr N_{\P(\scr E)})\times X}{(\P(\scr E)\times X)^c} &[-3em] \frac{\V(\scr N_{\P(\scr E)})\times X}{(X\times X)^c} & \frac{\V(\scr N_{X})\times X}{(X\times X)^c} \\
		& & \frac{\V(\scr N_{\P(\scr E)})\times X}{\Delta_X^c} & \frac{\V(\scr N_{X})\times X}{\Delta_X^c} \\
		& & \Th_X(\scr N_{(i\times\id)\delta_X} \oplus i^*(\scr N_{\P(\scr E)})) & \Th_X(\scr N_{\delta_X}\oplus\scr N_{X}) \\
		& & \Th_X(i^*(\scr N_{\delta_{\P(\scr E)}})\oplus i^*(\scr N_{\P(\scr E)})) & \Th_X(\Omega_X\oplus\scr N_{X}) \\
		& & \Th_X(i^*(\Omega_{\P(\scr E)})\oplus i^*(\scr N_{\P(\scr E)})) & \Th_X(f^*(\scr M))\rlap. \\
		};
		\arrows
		(11-) edge node[above]{$f^*(\eta_{\P(\scr E)})$} (-12) (12-) edge (-13) (13-) edge node[above]{$\Pi_i$} (-14)
		(12) edge (23) (13) edge (23) (14) edge (24) (23-) edge node[above]{$\Pi_i$} (-24)
		(23) edge node[left]{$\Pi_{(i\times\id)\delta}$} (33) (24) edge node[right]{$\Pi_{\delta}$} (34)
		(33-) edge node[above]{$\Psi$} (-34) (33) edge node[left]{$\simeq$} (43)
		(34) edge node[right]{$\nu$} (44) (43) edge node[left]{$\nu$} (53)
		(44) edge node[right]{$\alpha_X$} (54)
		(53-) edge node[below]{$i^*(\alpha_{\P(\scr E)})$} (-54)
		(53) edge node[above]{$\Psi$} (44)
		;
	\end{tikzmath}
	The lower composition is $i^*(\iota_{\P(\scr E)})$ and the upper composition is $\iota_X$. The middle rectangle commutes by Corollary~\ref{cor:wexcisive} applied to the smooth closed triple $(\V(\scr N_{\P(\scr E)})\times X,X\times X,\Delta_X)$. The commutativity of the trapezoid follows from the commutative diagram of canonical short exact sequences
	\begin{tikzmath}
		\def\rowsep{1em}
		\diagram{
		\delta_X^*(\scr N_{i\times\id}) & \scr N_{(i\times\id)\delta_X} & \scr N_{\delta_X}  \\
		& i^*(\scr N_{\delta_{\P(\scr E)}}) &\\
		\scr N_i & i^*(\Omega_{\P(\scr E)}) & \Omega_X\rlap. \\
		};
		\arrows (11-) edge[c->] (-12) (12-) edge[->>] (-13)
		(31-) edge[c->] (-32) (32-) edge[->>] (-33)
		(11) edge node[left]{$\simeq$} (31) (12) edge node[left]{$\simeq$} (22) (22) edge node[left]{$\nu$} (32) (13) edge node[right]{$\nu$} (33)
		;
	\end{tikzmath}
	Finally, the commutativity of the lower triangle is the definition of $\alpha_{(X,i)}$ in terms of $\alpha_{(\P(\scr E),\id)}$.
	
	From now on, we therefore assume that $X=\P(\scr E)$. As in the definition of $\eta$, we implicitly pull back everything along the $G$-affine bundle $v\colon V\to S$, so as to have a canonical isomorphism $\scr E\tens\scr E^\vee\simeq \scr O_S\oplus\scr Q$.
	Let $\gamma\colon \tilde X\into P\times \tilde X$ be the graph of the $G$-immersion $\sigma j\colon\tilde X\into P$.
	From the short exact sequence
	\[
	0\to \scr N_\gamma\to \sigma^*(\Omega_V)\oplus \Omega_{\tilde X}\xrightarrow{d\sigma+\id}\Omega_{\tilde X}\to 0,
	\]
	we obtain an isomorphism 
	\[\mu\colon \tilde f^*(\scr Q)=\sigma^*(\Omega_V)\simeq\scr N_\gamma\]
	sending a section $x$ of $\sigma^*(\Omega_V)$ to $(x,-(d\sigma)(x))$.
	We claim that the following rectangle commutes in $\H_\pt^G(X)$, where the first row is $\iota$:
	\begin{tikzequation}\label{eqn:keyduality}
		\diagram{\s^{f^*(\scr M)} & f^*(\Th_X(\scr N)) & \Sigma^{\Omega_f}\s^{\scr N} & \s^{f^*(\scr M)} \\
		\Sigma^{f^*(\scr M)}\tilde X_+ & \Sigma^{f^*(\scr E\tens\scr E^\vee)}\frac{P\times \tilde X}{H\times \tilde X} & \Sigma^{f^*(\scr E\tens\scr E^\vee)}\frac{P\times \tilde X}{\gamma^c} & \Sigma^{f^*(\scr M)}\tilde X_+\rlap. \\};
		\arrows (11-) edge node[above]{$f^*(\eta)$} (-12) (12-) edge node[above]{$\epsilon$} (-13) (13-) edge node[above]{$\alpha$} (-14)
		(21-) edge node[above]{$\simeq$} (-22) (22-) edge (-23) (23-) edge node[above]{$\mu\Pi$} (-24)
		(11) edge[<-] node[left]{$\pi$} node[right]{$\simeq$} (21) (14) edge[<-] node[right]{$\pi$} node[left]{$\simeq$} (24);
	\end{tikzequation}
Assuming this for the moment, let us conclude the proof of the theorem by showing that the bottom row of~\eqref{eqn:keyduality} is the identity.
	 Let $\zeta\colon \tilde X\into V\times\tilde X\subset P\times \tilde X$ be the zero section.
Consider the map
	 \[
	 \phi\colon \A^1\times V\times\tilde X\to V\times\tilde X,\quad (t,v,x)\mapsto (v-t\sigma(x),x).
	 \]
	 This is an $\A^1$-family of linear automorphisms of the $G$-vector bundle $V\times \tilde X$ over $\tilde X$.
	 Note that $\phi_0$ is the identity and $\phi_1\gamma=\zeta$.
	 Its projective completion is an $\A^1$-family
	 \[
	 \hat\phi\colon \A^1\times P\times\tilde X\to P\times\tilde X.
	 \]
	 We now consider the following diagram:
	 \begin{tikzmath}
	 	\diagram{\Sigma^{f^*(\scr Q)}\tilde X_+ & \frac{P\times\tilde X}{H\times\tilde X} & \frac{P\times \tilde X}{\gamma^c} & \Sigma^{f^*(\scr Q)}\tilde X_+ \\
		\Sigma^{f^*(\scr Q)}\tilde X_+ & \frac{P\times\tilde X}{H\times\tilde X} & \frac{P\times \tilde X}{\zeta^c} & \Sigma^{f^*(\scr Q)}\tilde X_+\rlap, \\};
		\arrows (11-) edge node[above]{$\simeq$} (-12) (12-) edge (-13) (13-) edge node[above]{$\Pi$} (-14)
		(21-) edge node[above]{$\simeq$} (-22) (22-) edge (-23) (23-) edge node[above]{$\Pi$} (-24)
		(11) edge node[left]{$\phi_1$} (21) (12) edge node[left]{$\hat\phi_1$} (22) (13) edge node[left]{$\hat\phi_1$} (23) (14) edge[-,vshift=1pt] (24) edge[-,vshift=-1pt] (24);
	 \end{tikzmath}
	 where for the bottom purity equivalence we use the obvious isomorphism $\scr N_\zeta\simeq\tilde f^*(\scr Q)$. The commutativity of the first two squares is clear.
	 The last square commutes because the isomorphism $\scr N_\gamma\simeq \scr N_\zeta$ restriction of $\gamma^*(d\phi_1)$ is exactly the composite $\scr N_\gamma\simeq \tilde f^*(\scr Q)\simeq \scr N_\zeta$ of the given isomorphisms.
	 The left vertical arrow is $\A^1$-homotopic to the identity via $\phi$, and it is clear that the lower row is the identity, by definition of the equivalence $\s^{\scr Q}\simeq P/H$.
	 
	 It remains to prove the commutativity of the rectangle~\eqref{eqn:keyduality}.
	 The proof is mostly straightforward but there are a few subtle points. 
	  We first note that there is a commutative square
	  \begin{tikzmath}
		\diagram{\s^{f^*(\scr M)} & \s^{f^*(\scr M)} \\ \Sigma^{f^*(\scr M)}\tilde X_+ & \Sigma^{f^*(\scr M)}\tilde X_+ \\};
		\arrows (11-) edge node[above]{$\iota$} (-12)
		(21) edge node[left]{$\pi$} (11)
		(22) edge node[right]{$\pi$} (12)
		(21-) edge node[above]{$\pi_\sharp\pi^*(\iota)$} (-22);
	  \end{tikzmath}
	  in $\H_\pt^G(X)$,
	  and that $\pi_\sharp\pi^*(\iota)$ is the composition of the following four maps:
	  \begin{enumerate}
	  	\item the first part of $\pi_\sharp\tilde f^*(\eta)$: $\pi_\sharp\tilde f^*$ of the $(\scr E\tens\scr E^\vee)$-suspension of~\eqref{eqn:thom1};
	  \item the second part of $\pi_\sharp\tilde f^*(\eta)$: $\pi_\sharp\tilde f^*\tilde f_\sharp$ of $\eqref{eqn:Kpath4}\1_{\tilde X}$, followed by the projection
	  \[\pi\times\id \colon \Th_{\tilde X\times\tilde X}(\pi_1^*\pi^*(\scr N))\to \Th_{X\times\tilde X}(\pi_1^*(\scr N));\]
	\item $\pi_\sharp\pi^*(\epsilon)$, which collapses the complement of $\tilde\delta\colon \tilde X\into X\times\tilde X$: 
	\[
	\Th_{X\times\tilde X}(\pi_1^*(\scr N))\to\frac{\V(\pi_1^*(\scr N))}{\tilde\delta^c}\stackrel\Pi\simeq \Th_{\tilde X}(\scr N_{\tilde\delta}\oplus\pi^*(\scr N))\stackrel\nu\simeq \Sigma^{\Omega_f\oplus\scr N}\tilde X_+;
	\]
	\item finally, the equivalence $\alpha$:
	\[
	\pi_\sharp\pi^*\eqref{eqn:Kpath5}\1_X: \Sigma^{\Omega_f\oplus\scr N}\tilde X_+\simeq \Sigma^{f^*(\scr M)}\tilde X_+.
	\]
	  \end{enumerate}
	  We must prove that $\pi_\sharp\pi^*(\iota)$ coincides with the lower row of~\eqref{eqn:keyduality}.
	  Let us first simplify steps (2)--(4); this is where the precise definitions of \eqref{eqn:Kpath4} and \eqref{eqn:Kpath5} come into play. Let ${\bar\delta}\colon \tilde X\times_X\tilde X\into \tilde X\times\tilde X$ be the obvious closed immersion, which is the pullback of $\tilde\delta$ along $\pi\times\id$.
	  We contemplate the following diagram in $\H_\pt^G( X)$, where $\Sigma'$ stands for $\Sigma^{f^*(\scr E\tens\scr E^\vee)}$:
	  \begin{tikzmath}
		  \def\colsep{1.5em}
	  	\diagram{\Sigma'\Th_{\tilde X\times\tilde X}(\pi_1^*(\scr N_s\oplus j^*(\scr N_\segre)))
		& \Sigma' \Th_{\tilde X\times\tilde X}(\pi_1^*\pi^*(\scr R)) & \Sigma'\Th_{X\times\tilde X}(\pi_1^*(\scr R)) \\
		\Sigma'\frac{\V(\pi_1^*(\scr N_s\oplus j^*(\scr N_\segre)))}{\bar\delta^c} & \Sigma'\frac{\V(\pi_1^*\pi^*(\scr R))}{\bar\delta^c} & \Sigma'\frac{\V(\pi_1^*(\scr R))}{\tilde\delta^c} \\
		\Sigma'\Th_{\tilde X\times_X\tilde X}(\scr N_{\bar\delta}\oplus \pi_1^*(\scr N_s\oplus j^*(\scr N_\segre))) & \Sigma'\Th_{\tilde X\times_X\tilde X}(\scr N_{\bar\delta}\oplus \pi_1^*\pi^*(\scr R)) & \Sigma'\Th_{\tilde X}(\pi^*(\scr N_{\delta})\oplus\pi^*(\scr R)) \\
				\Sigma'\Th_{\tilde X}(\pi^*(\Omega_X)\oplus \scr N_s\oplus j^*(\scr N_\segre)) & &
				\Sigma'\Th_{\tilde X}(\pi^*(\Omega_X)\oplus\pi^*(\scr R))
				\\
				\Sigma'\Th_{\tilde X}(\pi^*(\Omega_X)\oplus \scr N_s\oplus j^*(\scr N_\segre)) & & \Sigma'\Th_{\tilde X}(\tilde f^*(\scr Q))\rlap. \\};
		\arrows
		(11-) edge node[above]{\eqref{eqn:Kpath4}} (-12)
		(21-) edge node[above]{\eqref{eqn:Kpath4}} (-22)
		(31-) edge node[above]{\eqref{eqn:Kpath4}} (-32)
		(12-) edge node[above]{$\pi\times\id$} (-13)
		(22-) edge node[above]{$\pi\times\id$} (-23)
		(11) edge (21) (12) edge (22) (13) edge (23)
		(21) edge node[left]{$\Pi_{\bar\delta}$} (31) (22) edge node[left]{$\Pi_{\bar\delta}$} (32) (23) edge node[right]{$\Pi_{\tilde\delta}$} (33)
		
		(31) edge node[left]{$\nu\pi_1$} (41)
		(41-) edge node[above]{\eqref{eqn:Kpath4}} (-43)
		(32-) edge[-top,vshift=\dbl] node[above=\dbl]{$\pi_2$} (-33) edge[-bot,vshift=-\dbl] node[below=\dbl]{$\pi_1$} (-33)
		(51-) edge node[above]{\eqref{eqn:Kpath2}} (-53)
		(43) edge node[right]{\eqref{eqn:Kpath5}} (53)
		(33) edge node[right]{$\nu$} (43)
		(41) edge[-,vshift=1pt] (51) edge[-,vshift=-1pt] (51)
		;
	  \end{tikzmath}
	  The composition of the top row and the right column is the composition of steps (2)–(4).
	  Note that the parallel equivalences $\pi_1$ and $\pi_2$ are retractions of the same map and hence are homotopic in $\H_\pt^G(X)$. The commutativity of each square is clear except the last one.
	  Unfolding the definitions of the three equivalences involved, we see that its commutativity is equivalent to that of the following rectangle, where $\tau$ exchanges the two occurrences of $\Sigma^{\pi^*(\Omega)}$:
	  \begin{tikzmath}
		  \def\colsep{2em}
	  	\diagram{\Sigma^{\pi^*(\Omega)}\Sigma^{\pi^*(\scr R)}\Sigma^{\scr O}\Sigma^{\pi^*(\Omega)}\Sigma^{\scr N_s\oplus j^*(\scr N_\segre)} & \Sigma^{\pi^*(\Omega)}\Sigma^{\pi^*(\scr R)}\Sigma^{\scr O}\Sigma^{\tilde f^*(\scr Q)} & \Sigma^{\pi^*(\Omega)}\Sigma^{\pi^*(\scr R)}\Sigma^{\scr O}\Sigma^{\pi^*(\Omega)}\Sigma^{\pi^*(\scr R)} \\
		\Sigma^{\pi^*(\Omega)}\Sigma^{\pi^*(\scr R)}\Sigma^{\scr O}\Sigma^{\pi^*(\Omega)}\Sigma^{\scr N_s\oplus j^*(\scr N_\segre)} & \Sigma^{\pi^*(\Omega)}\Sigma^{\pi^*(\scr R)}\Sigma^{\scr O}\Sigma^{\tilde f^*(\scr Q)} & \Sigma^{\pi^*(\Omega)}\Sigma^{\pi^*(\scr R)}\Sigma^{\scr O}\Sigma^{\pi^*(\Omega)}\Sigma^{\pi^*(\scr R)}\rlap. \\};
		\arrows (11) edge node[left]{$\tau$} (21) (13) edge node[right]{$\tau$} (23)
		(11-) edge node[above]{\eqref{eqn:Kpath2}} (-12) (12-) edge node[above]{\eqref{eqn:Kpath1}} (-13)
		(21-) edge node[above]{\eqref{eqn:Kpath2}} (-22) (22-) edge node[above]{\eqref{eqn:Kpath1}} (-23);
	  \end{tikzmath}
	  This follows immediately from Lemma~\ref{lem:symmetry}.
	  
	  Thus, the composition of (2)–(4) coincides with the $f^*(\scr E\tens\scr E^\vee)$-suspension of the following composition:
	  \begin{multline}
		\label{eqn:asdf2}
	  \Th_{\tilde X\times\tilde X}(\pi_1^*(\scr N_s\oplus j^*(\scr N_\segre)))\to\frac{\V(\pi_1^*(\scr N_s\oplus j^*(\scr N_\segre)))}{\bar\delta^c}\stackrel\Pi\simeq \Th_{\tilde X\times_X\tilde X}(\scr N_{\bar\delta}\oplus \pi_1^*(\scr N_s\oplus j^*(\scr N_\segre)))
	  \\
	  \stackrel{\pi_1}\to \Th_{\tilde X}(\pi^*(\scr N_\delta)\oplus \scr N_s\oplus j^*(\scr N_\segre))\stackrel{\eqref{eqn:Kpath2}\nu}\simeq \Sigma^{f^*(\scr Q)} \tilde X_+.
	  \end{multline}
	  On the other hand, (1) is the $f^*(\scr E\tens\scr E^\vee)$-suspension of the following composition:
	  \begin{multline}
		  \label{eqn:asdf1}
	  	\Sigma^{f^*(\scr Q)}\tilde X_+\simeq \frac {P\times\tilde X}{H\times \tilde X} \to \frac {P\times\tilde X}{(H\cup \segre^c)\times\tilde X}\stackrel\Pi\simeq \frac{\Th_{X\times X^\vee\times\tilde X}(\pi_1^*(\scr N_{\segre}))}{\Th_{K\times\tilde X}(\pi_1^*(\scr N_{\segre}))}\\
	  	\stackrel\pi \from \frac{\Th_{\tilde X\times X^\vee\times\tilde X}(\pi^*\pi_1^*(\scr N_{\segre}))}{\Th_{\tilde K\times\tilde X}(\pi^*\pi_1^*(\scr N_{\segre}))}
	  	\to \frac{\V(\pi^*\pi_1^*(\scr N_{\segre}))\times\tilde X}{s^c\times\tilde X}\stackrel\Pi\simeq \Th_{\tilde X\times\tilde X}(\pi_1^*(\scr N_{s}\oplus j^*(\scr N_{\segre}))).
	  \end{multline}
	  We must therefore show that the ${f^*(\scr E\otimes\scr E^\vee)}$-suspension of $\eqref{eqn:asdf2}\circ\eqref{eqn:asdf1}$ coincides with the bottom row of~\eqref{eqn:keyduality}. Consider the closed immersion
	  \[
	  \beta=(j,\id)\colon \tilde X\into X\times X^\vee\times\tilde X,\quad x\mapsto (\pi x,\rho x, x),
	  \]
	  and the two closed immersions
	  \begin{gather*}
	  	\beta_1\colon \tilde X\times_X\tilde X\into \tilde X\times X^\vee\times \tilde X,\quad (x,y)\mapsto (x,\rho x, y),\\
	  	\beta_2\colon \tilde X\times_X\tilde X\into \tilde X\times X^\vee\times \tilde X,\quad (x,y)\mapsto (x,\rho y, y).
	  \end{gather*}
	   Note that $\beta_1$ and $\beta_2$ do not define the same closed subscheme of $\tilde X\times X^\vee\times\tilde X$. Instead, there is a commutative square
	  \begin{tikzmath}
	  	\diagram{\tilde X\times X^\vee\times\tilde X & \tilde X\times X^\vee\times\tilde X\\
		\tilde X\times_X\tilde X & \tilde X\times_X\tilde X\rlap, \\};
		\arrows (11-) edge node[above]{$\tau$} (-12) (21-) edge node[above]{$\tau$} (-22)
		(21) edge[c->] node[left]{$\beta_1$} (11) (22) edge[c->] node[right]{$\beta_2$} (12);
	  \end{tikzmath}
	  where $\tau$ exchanges the two copies of $\tilde X$, inducing an isomorphism $\tau^*(\scr N_{\beta_1})\simeq \scr N_{\beta_2}$. In particular, we have canonical isomorphisms
	  \[
	  \pi_1^*(\scr N_\beta)\simeq\scr N_{\beta_1}\quad\text{and}\quad \pi_2^*(\scr N_\beta)\simeq\scr N_{\beta_2}.
	  \]
	  This is used to define the arrows labeled $\pi_2$ and $\tau$ in the following diagram:
	  \begin{tikzmath}[ampersand replacement=\&]
		  \def\colsep{1.2em}
		  \def\rowsep{2em}
	  	\diagram{
		\frac{P\times\tilde X}{(H\cup\segre^c)\times \tilde X}
		\& \frac{P\times \tilde X}{\gamma^c}
		\& \Th_{\tilde X}(\scr N_\gamma)
		\\
		\frac{\Th_{X\times X^\vee\times\tilde X}(\pi_1^*(\scr N_\segre))}{\Th_{K\times X}(\pi_1^*(\scr N_\segre))}
		\& \frac{\V(\pi_1^*(\scr N_\segre))}{\beta^c}
		\& \Th_{\tilde X}(\scr N_\beta\oplus j^*(\scr N_\segre))
		\\
		\frac{\Th_{\tilde X\times X^\vee\times\tilde X}(\pi_1^*\pi^*(\scr N_\segre))}{\Th_{\tilde K\times X}(\pi_1^*\pi^*(\scr N_\segre))}
		\& \frac{\V(\pi_1^*\pi^*(\scr N_\segre))}{\beta_2^c}
		\& \Th_{\tilde X\times_X\tilde X}(\scr N_{\beta_2}\oplus \pi_2^*j^*(\scr N_\segre))
		\\
		\frac{\V(\pi_1^*\pi^*(\scr N_\segre))}{s^c\times\tilde X}
		\& \frac{\V(\pi_1^*\pi^*(\scr N_\segre))}{\beta_1^c}
		\& \Th_{\tilde X\times_X\tilde X}(\scr N_{\beta_1}\oplus \pi_1^*j^*(\scr N_\segre))
		\\
		\Th_{\tilde X\times\tilde X}(\pi_1^*(\scr N_s\oplus j^*(\scr N_\segre)))
		\& \frac{\V(\pi_1^*(\scr N_s\oplus j^*(\scr N_\segre)))}{\bar\delta^c}
		\& \Th_{\tilde X\times_X\tilde X}(\scr N_{\bar\delta}\oplus\pi_1^*(\scr N_s\oplus j^*(\scr N_\segre)))
		\\
		\& \& \Th_{\tilde X}(\pi^*(\scr N_\delta)\oplus \scr N_s\oplus j^*(\scr N_\sigma)) \\
		\& \& \Th_{\tilde X}(\tilde f^*(\scr Q))\rlap. \\
		};
		\arrows
		(11) edge node[left]{$\Pi_\segre$} (21) (21) edge[<-] node[left]{$\pi$} (31)
		(12) edge node[left]{$\Pi_\segre$} (22) (22) edge[<-] node[left]{$\pi$} (32)
		(13) edge node[right]{$\Psi$} (23) (23) edge[<-] node[right]{$\pi_2$} (33)
		(31) edge (41)
		(32) edge[draw=none] node[font=\normalsize]{(\textasteriskcentered)} (42)
		(33) edge node[right]{$\tau$} (43)
		(41) edge node[left]{$\Pi_{s\times\id}$} (51)
		(42) edge node[left]{$\Pi_{s\times\id}$} (52)
		(43) edge node[right]{$\Psi$} (53)
		(11-) edge (-12) (21-) edge (-22) (31-) edge (-32) (41-) edge (-42) (51-) edge (-52)
		(12-) edge node[above]{$\Pi_\gamma$} (-13)
		(22-) edge node[above]{$\Pi_\beta$} (-23)
		(32-) edge node[above]{$\Pi_{\beta_2}$} (-33)
		(42-) edge node[above]{$\Pi_{\beta_1}$} (-43)
		(52-) edge node[above]{$\Pi_{\bar\delta}$} (-53)
		(53) edge node[right]{$\pi_1$} (63) (63) edge node[right]{\eqref{eqn:Kpath2}$\nu$} (73)
		;
	  \end{tikzmath}
	  All the unlabeled arrows are quotient maps, and all the labeled ones are equivalences. The commutativity of the two top squares and of the two bottom squares follows from Corollary~\ref{cor:wexcisive}. The commutativity of the other two small squares is obvious. To complete the proof, we will show that:
	  \begin{itemize}
	  	\item[(a)] the composition of the right column is induced by the isomorphism $\mu\colon\scr N_\gamma\simeq\tilde f^*(\scr Q)$;
		\item[(b)] despite appearances, the rectangle (\textasteriskcentered) commutes.
	  \end{itemize}
	  
	  Let us prove (a). The first equivalence labeled $\Psi$ in the above diagram is induced by the short exact sequence of conormal sheaves
	  \begin{equation}
		  \label{eqn:SESgamma}
	  0\to j^*(\scr N_\segre) \to \scr N_\gamma \to \scr N_\beta \to 0
	  \end{equation}
	  associated with the triangle $\gamma=(\segre\times\id)\circ\beta$.
	  The second one is similarly induced by the short exact sequence associated with the triangle $\beta_1=(s\times\id)\circ\bar\delta$; it is the pullback by $\pi_1$ of the short exact sequence
	  \begin{equation}
		  \label{eqn:SESbeta}
	  0\to \scr N_s \to \scr N_\beta\to \pi^*(\scr N_\delta)\to 0
	  \end{equation}
	  associated with the triangle $\beta = (\id\times s)\circ \tilde\delta$. Thus, we have a commuting square
	  \begin{tikzmath}
	  	\diagram{
		\Th_{\tilde X\times_X\tilde X}(\scr N_{\beta_1}\oplus \pi_1^*j^*(\scr N_\segre)) & \Th_{\tilde X\times_X\tilde X}(\scr N_{\bar\delta}\oplus\pi_1^*(\scr N_s\oplus j^*(\scr N_\segre))) \\
		\Th_{\tilde X}(\scr N_\beta\oplus j^*(\scr N_\segre)) & \Th_{\tilde X}(\pi^*(\scr N_\delta)\oplus \scr N_s\oplus j^*(\scr N_\sigma))\rlap, \\
		};
		\arrows (11-) edge node[above]{$\Psi$} (-12) (21-) edge node[below]{$\Psi$} (-22)
		(11) edge node[left]{$\pi_1$} (21) (12) edge node[right]{$\pi_1$} (22);
	  \end{tikzmath}
	  allowing us to commute $\Psi$ and $\pi_1$. The resulting composition $\pi_1\tau\pi_2^{-1}$ is clearly the identity. To prove (a), it remains to show that the composition
	  \[
	  \Sigma^{\scr N_\gamma}\stackrel{\eqref{eqn:SESgamma}} \simeq \Sigma^{\scr N_\beta\oplus j^*(\scr N_\segre)}\stackrel{\eqref{eqn:SESbeta}}\simeq \Sigma^{\pi^*(\scr N_\delta)\oplus \scr N_s\oplus j^*(\scr N_\segre)}\stackrel\nu\simeq \Sigma^{\pi^*(\Omega_X)\oplus \scr N_s\oplus j^*(\scr N_\segre)}\stackrel{\eqref{eqn:Kpath2}}\simeq \Sigma^{\tilde f^*(\scr Q)}
	  \]
	  is induced by $\mu$. Recall that~\eqref{eqn:Kpath2} is induced by the short exact sequence
	  \[
	  0\to j^*(\scr N_\segre) \to \segre^*(\Omega_V) \xrightarrow{d\segre} \Omega_{\tilde X} \to 0
	  \]
	  and the obvious isomorphisms $\Omega_{\tilde X}\simeq \pi^*(\Omega_X)\oplus \rho^*(\Omega_{X^\vee})$ and $\rho^*(\Omega_{X^\vee})\simeq \scr N_s$.
	  Note that $\scr N_\beta$ is the conormal sheaf of the diagonal of $\tilde X$; let
	  \[
	  \nu'\colon \scr N_\beta\to \Omega_{\tilde X}
	  \]
	  be the isomorphism sending the class of $x\tens 1-1\tens x$ to $dx$. The claim now follows from the following isomorphisms of short exact sequences:
	  \[
	  \begin{tikzpicture}
		  \diagram{
		  j^*(\scr N_\segre) & \scr N_\gamma & \scr N_\beta \\
		  j^*(\scr N_\segre) & \segre^*(\Omega_V) & \Omega_{\tilde X}\rlap, \\
		  };
		  \arrows (11-) edge[c->] (-12) (12-) edge[->>] (-13)
		  (21-) edge[c->] (-22) (22-) edge[->>] (-23)
		  (11) edge[-,vshift=1pt] (21) edge[-,vshift=-1pt] (21) (12) edge node[right]{$\mu$} (22) (13) edge node[right]{$\nu'$} (23)
		  ;
		\end{tikzpicture}
	\qquad
  \begin{tikzpicture}
	  \diagram{
	  \scr N_s & \scr N_\beta & \pi^*(\scr N_\delta) \\
	  \rho^*(\Omega_{X^\vee}) & \Omega_{\tilde X} & \pi^*(\Omega_X)\rlap. \\
	  };
	  \arrows (11-) edge[c->] (-12) (12-) edge[->>] (-13)
		(21-) edge[c->] (-22) (22-) edge[->>] (-23)
		(11) edge node[left]{$\simeq$} (21) (12) edge node[right]{$\nu'$} (22) (13) edge node[right]{$\nu$} (23)
		;
	\end{tikzpicture}
	  \]
	  
	  To prove (b), note that both $\beta_1$ and $\beta_2$ factor through the closed $G$-immersion
	  \[
	  \lambda\colon \tilde X\times_X\tilde X\times X^\vee \into \tilde X\times X^\vee\times \tilde X.
	  \]
	  Let us write $\beta_{1,2} = \lambda\circ\zeta_{1,2}$, where 
	  \[\zeta_{1,2}\colon \tilde X\times_X\tilde X\into \tilde X\times_X\tilde X\times X^\vee,\]
	  and let $\kappa\colon \tilde X\times_X\tilde X\times X^\vee\to X\times X^\vee$ be the projection.
	  We can then break up the rectangle (\textasteriskcentered) as follows:
	  \begin{tikzmath}
		  \def\colsep{-2em}
	  	\diagram{
		\frac{\V(\pi_1^*\pi^*(\scr N_\segre))}{\beta_2^c} & \frac{\V(\pi_1^*\pi^*(\scr N_\segre))}{\lambda^c} & \frac{\V(\pi_1^*\pi^*(\scr N_\segre))}{\beta_1^c} \\
		\frac{\V(\scr N_\lambda\oplus \kappa^*(\scr N_\segre))}{\zeta_2^c} & \Th_{\tilde X\times_X\tilde X\times X^\vee}(\scr N_\lambda\oplus \kappa^*(\scr N_\segre)) & \frac{\V(\scr N_\lambda\oplus \kappa^*(\scr N_\segre))}{\zeta_1^c} \\
		\Th_{\tilde X\times_X\tilde X}(\scr N_{\zeta_2}\oplus\zeta_2^*(\scr N_\lambda)\oplus \pi_2^*j^*(\scr N_\segre)) & & \Th_{\tilde X\times_X\tilde X}(\scr N_{\zeta_1}\oplus\zeta_1^*(\scr N_\lambda)\oplus \pi_1^*j^*(\scr N_\segre)) \\
		\Th_{\tilde X\times_X\tilde X}(\scr N_{\beta_2}\oplus \pi_2^*j^*(\scr N_\segre)) & & \Th_{\tilde X\times_X\tilde X}(\scr N_{\beta_1}\oplus \pi_1^*j^*(\scr N_\segre))\rlap. \\
		};
		\arrows (11-) edge[<-] (-12) (12-) edge (-13) (21-) edge[<-] (-22) (22-) edge (-23)
		(31-) edge node[above]{$\tau$} (-33) (41-) edge node[above]{$\tau$} (-43)
		(11) edge node[left]{$\Pi_\lambda$} (21) (12) edge node[left]{$\Pi_\lambda$} (22) (13) edge node[right]{$\Pi_\lambda$} (23)
		(21) edge node[left]{$\Pi_{\zeta_2}$} (31) (23) edge node[right]{$\Pi_{\zeta_1}$} (33)
		(31) edge node[left]{$\Psi$} (41) (33) edge node[right]{$\Psi$} (43)
		;
	  \end{tikzmath}
	  Only the commutativity of the middle rectangle is not clear. However, it is clear that the middle rectangle commutes if we add in the automorphism $\tau$ of $\Th_{\tilde X\times_X\tilde X\times X^\vee}(\scr N_\lambda\oplus \kappa^*(\scr N_\segre))$ that permutes the two factors of $\tilde X$ (defined on the Thom space since $\lambda\tau=\tau\lambda$ and $\tau\kappa=\kappa$). The claim follows from the observation that this automorphism is homotopic to the identity, since the diagonal $\tilde X\into \tilde X\times_X\tilde X$ is a section of a $G$-vector bundle.
\end{proof}

\section{Stable equivariant motivic homotopy theory}
\label{sec:SH}

In this section, we construct the stable equivariant motivic $\infty$-category $\SH^G(S)$ for a $G$-scheme $S$.
Recall that we have fixed a qcqs base scheme $B$ and a tame group scheme $G$ over $B$.
Throughout this section, we moreover assume that one of the following conditions holds:
\begin{itemize}
	\item $G$ is finite locally free; or
	\item $B$ has the $G$-resolution property.
\end{itemize}
This restriction is not essential, but the definition of $\SH^G(S)$ in general is more complicated; it is however determined by the requirement that $\SH^G(\ph)$ be a Nisnevich sheaf. 

Let $\Pr^{\mathrm L}$ denote the $\infty$-category whose objects are presentable $\infty$-categories and whose morphisms are colimit-preserving functors.
Recall that $\Pr^{\mathrm L}$ admits limits and colimits, and that the former are computed in $\Cat_\infty$ \cite[\S5.5.3]{HTT}.
It also admits a symmetric monoidal structure $\Pr^{\mathrm L,\tens}$ \cite[Proposition 4.8.1.14]{HA}, and we call a commutative algebra in $\Pr^{\mathrm L,\otimes}$ a \emph{presentably symmetric monoidal} $\infty$-category. For $\scr C\in\CAlg(\Pr^{\mathrm L,\tens})$, a \emph{$\scr C$-module} will always mean a $\scr C$-module in the symmetric monoidal $\infty$-category $\Pr^{\mathrm L,\tens}$; we denote by $\Mod_\scr C$ the $\infty$-category of $\scr C$-modules.

\subsection{Equivariant motivic spectra}
\label{sub:SH}

To transform the ``suspended adjunction'' of Theorem~\ref{thm:duality} into a genuine adjunction, we need to make the spheres $\s^{\scr M}$ invertible for the tensor product of pointed motivic $G$-spaces. 
To that end, we use the formalism developed by Robalo in \cite[\S2.1]{Robalo}.
Given a presentably symmetric monoidal $\infty$-category $\scr C$
and a set of objects $X$ in $\scr C$,
 there exists a functor
\[
\scr C\to \scr C[X^{-1}]
\]
with the following universal property in the $\infty$-category $\CAlg(\Pr^{\mathrm L,\otimes})$: any functor $f\colon \scr C\to \scr D$ such that $f(x)$ is invertible for all $x\in X$ factors uniquely through $\scr C[X^{-1}]$.
In particular, the $\infty$-category of $\scr C[X^{-1}]$-modules (in $\Pr^{\mathrm L,\otimes}$) is the full subcategory of $\scr C$-modules on which all the objects of $X$ act by equivalences.
In \emph{loc.\ cit.}, this construction is considered only when $X$ has a single element, but it is clear that the filtered colimit
\[
\scr C[X^{-1}]=\colim_{\substack{F\subset X \\ F\text{ finite}}}\scr C[(\bigotimes F)^{-1}]
\]
satisfies the required universal property.

If $\scr M$ is a $\scr C$-module and $x\in\scr C$, let $\Stab_x(\scr M)$ denote the colimit of the sequence
\[
\scr M\xrightarrow{\ph\tens x} \scr M\xrightarrow{\ph\tens x}\scr M\xrightarrow{\ph\tens x}\dotsb
\]
in $\Mod_{\scr C}$. Note that the underlying $\infty$-category of $\Stab_x(\scr M)$ is
the limit of the tower
\[
\dotsb\xrightarrow{\Hom(x,\ph)}\scr M\xrightarrow{\Hom(x,\ph)}\scr M\xrightarrow{\Hom(x,\ph)}\scr M.
\]
More generally, we define the $\scr C$-module $\Stab_X(\scr M)$ as follows. Let $L$ denote the $1$-skeleton of the nerve of the poset $\N$, and let $L(X)$ be the simplicial set of almost zero maps $X\to L$. Note that $L(X)$ is a filtered simplicial set and is the union of the simplicial subsets $L(F)$ for finite subsets $F\subset X$.
Using the symmetric monoidal structure on $\scr C$, we can construct a diagram $L(X)\to\Mod_{\scr C}$ sending each vertex to $\scr M$ and each edge in the $x$-direction to the functor $(\ph)\tens x$.\footnote{More precisely, this amounts to defining a functor $P_\mathrm{fin}(X)\to\Mod_{\scr C}$, where $P_\mathrm{fin}(X)$ is the poset of finite subsets of $X$. Considering the universal case, we have to define a functor $P_\mathrm{fin}(X)\to B(\coprod_{n\geq 0}X^n/\Sigma_n)$, where the target is a $(2,1)$-category. We can construct the desired functor as a simplicial map from the nerve of $P_\mathrm{fin}(X)$ to the simplicial bar construction on the monoidal groupoid $\coprod_{n\geq 0}X^n/\Sigma_n$, sending a $k$-simplex $Y_0\subset\dotsb\subset Y_k$ to $(Y_1\minus Y_0,\dotsc,Y_k\minus Y_{k-1})$.}
 We let $\Stab_X(\scr M)$ be the colimit of this diagram.
 Cofinality considerations show that
\[
\Stab_X(\scr M)=\colim_{\substack{F\subset X \\ F\text{ finite}}}\Stab_{\bigotimes F}(\scr M).
\]
Informally speaking, an object in $\Stab_X(\scr M)$ is an object of $\scr M$ equipped with compatible $w$-deloopings for $w$ any finite tensor product of elements of $X$.

Recall from \cite[Remark 2.20]{Robalo} that an object $x\in\scr C$ is \emph{$n$-symmetric} if the cyclic permutation of $x^{\tens n}$ is homotopic to the identity. If $x\in \scr C$ is $n$-symmetric for some $n\geq 2$ and $\scr M$ is a $\scr C$-module, then $x$ acts on $\Stab_x(\scr M)$ by an equivalence. If this holds for every $x\in X$, there results a canonical map of $\scr C[X^{-1}]$-modules
\[
\scr M\tens_{\scr C}\scr C[X^{-1}] \to \Stab_X(\scr M),
\]
which is an equivalence \cite[Corollary 2.22]{Robalo}. In particular, we obtain the following explicit description of $\scr C[X^{-1}]$ as a $\scr C$-module:
\[
\scr C[X^{-1}]\simeq \Stab_X(\scr C).
\]

If $\scr C\to\scr D$ is a morphism of presentably symmetric monoidal $\infty$-categories, there is an induced base change functor
\[
(\ph)\tens_{\scr C}\scr D \colon \Mod_{\scr C}\to\Mod_{\scr D}
\]
between the $\infty$-categories of modules. Below we will need to know that an adjoint pair
$
f:\scr M \rightleftarrows \scr N:g
$
in $\Mod_{\scr C}$ (\ie, both $f$ and $g$ are morphisms of $\scr C$-modules in $\Pr^{\mathrm L,\tens}$) gives rise to an adjoint pair
$
\scr M\tens_{\scr C}\scr D \rightleftarrows \scr N\tens_{\scr C}\scr D
$
in $\Mod_{\scr D}$. This follows from the fact that the above base change functor is in fact an $(\infty,2)$-functor and thus preserves adjunctions. Unfortunately, we do not know a reference for this fact, so we provide an alternative argument in the special case that is relevant for our purposes. Suppose that $\scr D=\scr C[X^{-1}]$ and that each $x\in X$ is $n$-symmetric for some $n\geq 2$, so that $\scr M\tens_{\scr C}\scr D\simeq\Stab_{X}(\scr M)$ as $\scr C$-modules. The claim is then that the functors $\Stab_{X}(f)$ and $\Stab_X(g)$ form an adjoint pair. If we write $\Stab_X(\scr M)$ and $\Stab_X(\scr N)$ as cofiltered limits of $\infty$-categories indexed by $L(X)^\op$, the right adjoint to $\Stab_X(f)$ is the functor induced in the limit by the right adjoint to $f$, and similarly for $\Stab_X(g)$.
 The claim thus follows from the fact that limits of $\infty$-categories preserve adjunctions, since they preserve unit transformations \cite[Definition 5.2.2.7]{HTT}. In particular, if $f$ or $g$ is fully faithful, so is $\Stab_X(f)$ or $\Stab_X(g)$.

We are now ready to define the stable equivariant motivic homotopy $\infty$-category.
If $S$ is a $G$-scheme, we denote by $\Sph_S$ the collection of all spheres $\s^{\scr E}$ in $\H_\pt^G(S)$, where $\scr E$ is a locally free $G$-module of finite rank on $S$. 

\begin{definition}
	Let $S$ be a $G$-scheme with structure map $p\colon S\to B$.
	The symmetric monoidal $\infty$-category of \emph{motivic $G$-spectra} over $S$ is defined by
	\[
	\SH^G(S)=\H_\pt^G(S)[p^*(\Sph_B)^{-1}].
	\]
	We denote by
	\[\Sigma^\infty: \H^G_\pt(S)\rightleftarrows\SH^G(S) : \Omega^\infty\]
	the canonical adjunction, where $\Sigma^\infty$ is symmetric monoidal.
\end{definition}

 We will see in Corollary~\ref{cor:thomequiv} below that $\SH^G(S)=\H_\pt^G(S)[\Sph_S^{-1}]$, which is the intended definition, but the above definition makes it easier to extend the functorial properties of $\H_\pt^G(\ph)$ to $\SH^G(\ph)$. If $\scr E$ is a locally free $G$-module on $S$,
we will denote by $\s^{-\scr E}$ the $\tens$-inverse of $\s^{\scr E}=\Sigma^\infty\s^{\scr E}$ in $\SH^G(S)$. Note that $\SH^G(S)$ is stable since $\s^{\scr O}\simeq\Sigma(\A^1\minus 0,1)$.

For every $G$-morphism $f\colon T\to S$, the universal property of $\SH^G(S)$ yields an adjunction
\[
f^* : \SH^G(S)\rightleftarrows \SH^G(T): f_*
\]
where the left adjoint is symmetric monoidal, such that $f^*\Sigma^\infty=\Sigma^\infty f^*$ and $\Omega^\infty f_*=f_*\Omega^\infty$.

\begin{lemma}\label{lem:SHBC}
	For every $G$-morphism $f\colon T\to S$, the functor
	\[
	\H_\pt^G(T)\tens_{\H_\pt^G(S)} \SH^G(S)\to \SH^G(T)
	\]
	induced by $f^*\colon \SH^G(S)\to \SH^G(T)$
	is an equivalence of symmetric monoidal $\infty$-categories.
\end{lemma}

\begin{proof}
	Compare universal properties.
\end{proof}

\begin{lemma}\label{lem:cyclic}
	Let $S$ be a $G$-scheme. Every sphere $\s^{\scr E}\in \H_\pt^G(S)$ is $3$-symmetric.
\end{lemma}

\begin{proof}
	This follows at once from Lemma~\ref{lem:symmetry}.
\end{proof}

\begin{proposition}
\label{prop:SHgeneration}
Let $S$ be a $G$-scheme with structure map $p\colon S\to B$.
\begin{enumerate}
	\item There is a canonical equivalence of $\H_\pt^G(S)$-modules
	\[
	\SH^G(S) \simeq \Stab_{p^*(\Sph_B)}\H_\pt^G(S).
	\]
	\item The $\infty$-category $\SH^G(S)$ is generated under sifted colimits by $E^{-1}\tens \Sigma^\infty X_+$, where $E\in p^*(\Sph_B)$ and $X\in\Sm_S^G$ is small and affine.
	\item For every $X\in\Sm_S^G$, $\Sigma^\infty X_+$ is compact in $\SH^G(S)$.
	\item For every $G$-morphism $f\colon T\to S$, the functor $f_*\colon \SH^G(T)\to\SH^G(S)$ preserves colimits.
\end{enumerate}
\end{proposition}

\begin{proof}
	(1) This follows from Lemma~\ref{lem:cyclic} and the above discussion.
	
	(2) By (1) and \cite[Lemma 6.3.3.6]{HTT}, every $E\in\SH^G(S)$ can be written as a filtered colimit of objects of the form $E^{-1}\tens \Sigma^\infty X$ with $X\in \H^G_\pt(S)$. Since the adjunction $\H^G(S)\rightleftarrows \H^G_\pt(S)$ is monadic, every such $X$ is a simplicial colimit of objects in the image of the left adjoint, and we conclude with Proposition~\ref{prop:generators} (1).
	
	(3) By Proposition~\ref{prop:generators} (3), the objects $X_+$ for $X\in\Sm_S^G$ are compact generators of $\H^G_\pt(S)$. In particular, every Thom space $\Th_X(\scr E)$ is compact in $\H_\pt^G(S)$, being a pushout of compact objects, and it follows that $\Omega^{\scr E}=\Hom(\s^{\scr E},\ph)$ preserves filtered colimits. By (1), this implies that $\Omega^\infty$ preserves filtered colimits and hence that $\Sigma^\infty$ preserves compact objects.
	
	(4) This follows from (2) and (3).
\end{proof}

It follows from
Lemma~\ref{lem:SHBC}
that, if $f\colon X\to S$ is a smooth $G$-morphism, the $\H_\pt^G(S)$-module adjunction
\[
f_\sharp : \H_\pt^G(X) \rightleftarrows \H_\pt^G(S) : f^*
\]
induces by base change along $\H_\pt^G(S)\to\SH^G(S)$ an $\SH^G(S)$-module adjunction
\[
f_\sharp : \SH^G(X) \rightleftarrows \SH^G(S): f^*.
\]
Similarly, if $i\colon Z\into S$ is a closed $G$-immersion, the $\H_\pt^G(S)$-module adjunction
\[
i^* : \H_\pt^G(S) \rightleftarrows \H_\pt^G(Z) : i_*
\]
induces by base change along $\H_\pt^G(S)\to\SH^G(S)$ an $\SH^G(S)$-module adjunction
\[
i^* : \SH^G(S) \rightleftarrows \SH^G(Z) : i_*.
\]
In particular, $i_*$ preserves colimits and we also have an adjunction
\[
i_* : \SH^G(Z) \rightleftarrows \SH^G(S) : i^!.
\]

In summary, the functors $f_\sharp$ and $i_*$, for $f$ smooth and $i$ a closed immersion, extend to $\SH^G(\ph)$ together with the smooth and closed projection formulas. 
Moreover, the following properties are immediately deduced from their unstable (pointed) analogs using the functoriality of base change along $\Sigma^\infty$: smooth base change, closed base change,
smooth–closed base change,
 and the full faithfulness of $f^*$, $j_*$, $j_\sharp$, and $i_*$ for $f$ a $G$-affine bundle, $j$ an open $G$-immersion, and $i$ a closed $G$-immersion.
The gluing and purity theorems follow easily from their unstable versions and Proposition~\ref{prop:SHgeneration} (2), and the Nisnevich separation property from the fact that there is, for every $S$, a conservative family of functors $\SH^G(S)\to \H^G_\pt(S)$ that commute with smooth base change.
Thus, all the functoriality discussed so far for $\H_\pt^G(\ph)$ extends to $\SH^G(\ph)$. What is perhaps less obvious is that arbitrary spheres are invertible in $\SH^G(S)$:

\begin{proposition}
	\label{prop:thomequiv}
	Let $p\colon X\to S$ be a smooth $G$-morphism with a closed $G$-section $s$. Then the adjunction 
	\[p_\sharp s_*: \SH^G(S)\rightleftarrows \SH^G(S): s^!p^*\]
	is an equivalence of $\infty$-categories. 
\end{proposition}

\begin{proof}
	By purity and the projection formulas, the left adjoint can be identified with $\s^{\scr E}\tens (\ph)$ where $\scr E$ is the conormal sheaf of $s$. Thus, the proposition is equivalent to the statement that $\s^{\scr E}$ is invertible in $\SH^G(S)$.
	By Proposition~\ref{prop:affinecover} and the fact that smooth base change is closed symmetric monoidal, we may assume that $S$ is small and affine.
	Let $r\colon S\to B$ be the structure map.
	We claim that there exists a locally free $G$-module $\scr F$ on $B$ and an epimorphism $r^*(\scr F)\onto\scr E$.
	Assuming this claim, we have $\s^{r^*(\scr F)}\simeq \s^{\scr E}\tens \s^{\scr G}$ in $\H_\pt^G(S)$, where $\scr G$ is the kernel of $r^*(\scr F)\onto\scr E$, and since $\s^{r^*(\scr F)}$ is invertible in $\SH^G(S)$, so is $\s^{\scr E}$.
	
	It therefore suffices to establish the claim.
	If $G$ is finite, we can find an epimorphism $r^*(\scr O_G^n)\onto \scr E$, where $\scr O_G$ is the regular representation of $G$ over $B$.
	If $G$ is not finite, then by assumption $B$ has the $G$-resolution property.
	As $r\colon S\to B$ is quasi-affine, $r^*r_*(\scr E)\to\scr E$ is an epimorphism.
	By Lemma~\ref{lem:thomason}, $r_*(\scr E)$ is the union of its finitely generated quasi-coherent $G$-submodules. Hence, there exists a finitely generated quasi-coherent $G$-module $\scr M$ on $B$ and an epimorphism $r^*(\scr M)\onto \scr E$. By the $G$-resolution property, there exists a locally free $G$-module $\scr F$ on $B$ and an epimorphism $\scr F\onto \scr M$. We therefore obtain an epimorphism $r^*(\scr F)\onto \scr E$, as desired.
\end{proof}

\begin{remark}
	The use of purity in the proof of Proposition~\ref{prop:thomequiv} is not essential: after reducing to the case where $S$ is small and affine and $p$ is $G$-quasi-projective, one can obtain an equivalence $X/(X\minus S)\simeq \s^{\scr E}$ directly from the results of \S\ref{sec:geometry} (purity implies that this equivalence is independent of choices).
\end{remark}

\begin{corollary}
	\label{cor:thomequiv}
	Let $p\colon S\to B$ be a $G$-morphism. Then the canonical functor
	\[
	\SH^G(S)=\H_\pt^G(S)[p^*(\Sph_B)^{-1}]\to \H_\pt^G(S)[\Sph_S^{-1}]
	\]
	is an equivalence of symmetric monoidal $\infty$-categories.
\end{corollary}

\begin{remark}\label{rmk:regular}
	If $G$ is finite locally free, we also have
	\[
	\SH^G(S)\simeq \H_\pt^G(S)[(\s^{p^*(\scr O_G)})^{-1}],
	\]
	where $\scr O_G$ is the regular representation of $G$.
	Indeed, the proof of Proposition~\ref{prop:thomequiv} shows that any sphere $\s^{\scr E}$ becomes invertible in the right-hand side. 
\end{remark}

Let $f\colon X\to S$ be a smooth separated $G$-morphism. As in \S\ref{sub:pt}, we can define a natural transformation
\[\epsilon\colon f^*f_\sharp\Sigma^{-\Omega_f} \to \id\colon \SH^G(X) \to \SH^G(X)\]
by the composition~\eqref{eqn:pure}. By adjunction, we obtain a natural transformation $f_\sharp\Sigma^{-\Omega_f} \to f_*$.

\begin{theorem}[Ambidexterity]
	\label{thm:stableduality}
	Let $f\colon X\to S$ be a smooth proper $G$-morphism. Then the transformation
	\[f_\sharp\Sigma^{-\Omega_f} \to f_*\colon \SH^G(X)\to \SH^G(S)\]
	is an equivalence.
\end{theorem}

\begin{proof}
	Any proper morphism in $\Sch_B^G$ is $G$-projective Nisnevich-locally on $B$. By Nisnevich separation and smooth base change, we can therefore assume that $f$ is smooth and $G$-projective.
	In that case, Theorem~\ref{thm:duality} shows that the transformation $\epsilon\colon f^*f_\sharp\Sigma^{-\Omega_f} \to \id$ is the counit of an adjunction $f^*\dashv f_\sharp\Sigma^{-\Omega_f}$ (between the homotopy $1$-categories), whence the result.
\end{proof}

\begin{corollary}[Proper base change]
	\label{cor:properBC}
	Let
	\begin{tikzmath}
		\diagram{Y' & Y \\ X' & X \\};
		\arrows (11-) edge node[above]{$g$} (-12) (11) edge node[left]{$q$} (21) (21-) edge node[below]{$f$} (-22) (12) edge node[right]{$p$} (22);
	\end{tikzmath}
	be a cartesian square of $G$-schemes where $p$ is proper. Then the exchange transformation
	\[
	\Ex^*_*\colon f^*p_* \to q_*g^* \colon \SH^G(Y) \to \SH^G(X')
	\]
	is an equivalence.
\end{corollary}

\begin{proof}
	Any proper morphism in $\Sch_B^G$ is $G$-projective Nisnevich-locally on $B$. By Nisnevich separation and smooth base change, we can assume that $p$ is a closed $G$-immersion or that $p$ is smooth and $G$-projective. In the former case, the result holds by closed base change. In the latter case, we note that the square
	\begin{tikzmath}
		\diagram{
		f^*p_\sharp\Sigma^{-\Omega_p} & f^*p_* \\
		q_\sharp \Sigma^{-\Omega_q}g^* & q_*g^* \\
		};
		\arrows (11-) edge (-12) (21-) edge (-22) (11) edge[<-] node[left]{$\Ex_\sharp^*$} (21) (12) edge node[right]{$\Ex^*_*$} (22);
	\end{tikzmath}
	is commutative for formal reasons (see the proof of \cite[Lemma 2.4.23 (1)]{CD}).
	 Hence, the desired result follows from Theorem~\ref{thm:stableduality} and smooth base change.
\end{proof}

\begin{corollary}[Proper projection formula]
	\label{cor:properproj}
	Let $p\colon Y\to X$ be a proper $G$-morphism. For every $A\in\SH^G_\pt(X)$ and $B\in\SH^G_\pt(Y)$, the canonical map
	\[
		A\tens p_*B\to p_*(p^*A\tens B)
	\]
	is an equivalence.
\end{corollary}

\begin{proof}
	This follows from the smooth and closed projection formulas as in the proof of Corollary~\ref{cor:properBC}.
\end{proof}

\begin{corollary}[Smooth–proper base change]
	\label{cor:deligne}
	Let
	\begin{tikzmath}
			\diagram{Y' & Y \\ X' & X \\};
			\arrows (11-) edge node[above]{$g$} (-12) (11) edge node[left]{$q$} (21) (21-) edge node[below]{$f$} (-22) (12) edge node[right]{$p$} (22);
		\end{tikzmath}
	be a cartesian square of $G$-schemes where $p$ is proper and $f$ is smooth. Then the exchange transformation
	\[
	\Ex_{\sharp *}\colon f_\sharp q_*\to p_* g_\sharp\colon \SH^G(Y') \to \SH^G(X)
	\]
	is an equivalence.
\end{corollary}

\begin{proof}
	This follows from smooth–closed base change as in the proof of Corollary~\ref{cor:properBC}.
\end{proof}

\begin{corollary}[Atiyah duality]
	\label{cor:atiyah}
	Let $f\colon X\to S$ be a smooth proper $G$-morphism. Then $\Sigma^\infty X_+$ is strongly dual to $f_\sharp\s^{-\Omega_f}$ in $\SH^G(S)$.
\end{corollary}

\begin{proof}
	Recall that $A\in\SH^G(S)$ is strongly dualizable if and only if, for every $B\in \SH^G(S)$, the canonical map $\Hom(A,\1_S)\tens B\to \Hom(A,B)$ is an equivalence. By the smooth projection formula, we have $\Hom(f_\sharp\1_X,B)\simeq f_*f^*B$, and the canonical map $f_*\1_X\tens B\to f_*f^*B$ is an equivalence by the proper projection formula. Thus, $f_\sharp\1_X$ is strongly dualizable, with dual $f_*\1_X$. By Theorem~\ref{thm:stableduality}, $f_*\1_X\simeq f_\sharp\s^{-\Omega_f}$. 
\end{proof}

\subsection{The exceptional functors}

We now show how the results of the previous sections give rise to the formalism of six operations for $\SH^G(\ph)$. 
It is not difficult to define a functor
\begin{equation}\label{eqn:functor}
\SH^G\colon (\Sch_B^G)^\op \to \CAlg(\Pr^{\mathrm L,\otimes})
\end{equation}
that sends a $G$-scheme $S$ to $\SH^G(S)$ and a $G$-morphism $f\colon T\to S$ to $f^*\colon \SH^G(S)\to \SH^G(T)$; we refer to \cite[\S9.1]{RobaloThesis} for details.

To define the pushforward with compact support $f_!$ at the level of $\infty$-categories, we will use the technology developed in \cite[Chapter V.1]{GR}.
A $G$-morphism $f\colon Y\to X$ in $\Sch_B^G$ is called \emph{compactifiable} if there exists a proper $G$-scheme $P$ in $\Sch_B^G$ such that $f$ is the composition of a $G$-immersion $Y\into P\times_BX$ and the projection $P\times_BX\to X$. It is easy to show that the composition of two compactifiable $G$-morphisms is compactifiable, so that compactifiable $G$-morphisms define a wide subcategory $(\Sch_B^G)_\cpt$ of $\Sch_B^G$. Note also that if $g\circ f$ is compactifiable, then $f$ is compactifiable; in particular, any $G$-morphism between $G$-quasi-projective $B$-schemes is compactifiable.
By Lemma~\ref{lem:Gimmersion} (2), any compactifiable $G$-morphism $f$ can be written as $p\circ j$ where $j$ is an open $G$-immersion and $p$ is a proper compactifiable $G$-morphism; such a factorization is called a \emph{compactification} of $f$. 
It is then clear that the category $(\Sch_B^G)_\cpt$, equipped with its wide subcategories of open immersions and proper morphisms, satisfies the assumptions of \cite[Exposé XVII, 3.2.4]{SGA4-3}. By \cite[Exposé XVII, Proposition 3.2.6 (ii)]{SGA4-3}, the category of compactifications of any compactifiable $G$-morphism is cofiltered, and in particular weakly contractible.

Given $S\in\Sch_B^G$, let $\Corr(\Sch_S^G)_{\all,\cpt}^\prop$ denote the $2$-category whose $1$-morphisms are spans $X\from Y\to Z$ in $\Sch_S^G$ with $Y\to Z$ compactifiable and whose $2$-morphisms are proper $G$-morphisms between spans \cite[Chapter V.1, \S1]{GR}.
Restricting~\eqref{eqn:functor} to $\Sch_S^G$ gives a functor 
\begin{equation}\label{eqn:functor2}
	(\Sch_S^G)^\op\to \Mod_{\SH^G(S)}.
\end{equation}
We claim that this functor extends uniquely\footnote{Uniqueness means that the $\infty$-groupoid of such extensions is contractible.} to an $(\infty,2)$-functor
\begin{equation}\label{eqn:functor3}
\Corr(\Sch_S^G)_{\all,\cpt}^{\prop,\operatorname{2-op}} \to \Mod_{\SH^G(S)}
\end{equation}
satisfying the following condition:
	\begin{itemize}
		\item[(\textasteriskcentered)] Let $U\into X$ be an open $G$-immersion in $\Sch_S^G$. Then the canonical $2$-isomorphism between $\id_U$ and the composition of the two spans $U= U\into X$ and $X\hookleftarrow U= U$ becomes the unit of an adjunction in $\Mod_{\SH^G(S)}$.
	\end{itemize}
Here, ``$\operatorname{2-op}$'' means that we reverse the direction of the $2$-morphisms.
We first note that~\eqref{eqn:functor2} sends smooth (\resp{} proper) $G$-morphisms to left (\resp{} right) adjointable morphisms in $\Mod_{\SH^G(S)}$, by the smooth projection formula (\resp{} by the proper projection formula and Proposition~\ref{prop:SHgeneration} (4)).
Smooth base change for open immersions allows us to apply \cite[Chapter V.1, Theorem 3.2.2]{GR}: the functor~\eqref{eqn:functor2} admits a unique extension to the $2$-category $\Corr(\Sch_S^G)_{\all,\open}^\open$, whose $1$-morphisms are spans $X\from Y\to Z$ with $Y\to Z$ an open $G$-immersion.
Its restriction to the $(2,1)$-category $\Corr(\Sch_S^G)_{\all,\open}$ satisfies condition (\textasteriskcentered), and by \cite[Chapter V.1, Theorem 4.1.3]{GR} it is the unique extension of~\eqref{eqn:functor2} with this property. We now apply \cite[Chapter V.1, Theorem 5.2.4]{GR} with proper compactifiable morphisms as admissible morphisms and open immersions as co-admissible morphisms: this is justified by the weak contractibility of the categories of compactifications, proper base change (Corollary~\ref{cor:properBC}), and smooth–proper base change for open immersions (Corollary~\ref{cor:deligne}).
As a result, there is a unique further extension of~\eqref{eqn:functor2} from $\Corr(\Sch_S^G)_{\all,\open}$ to $\Corr(\Sch_S^G)_{\all,\cpt}^{\prop,\operatorname{2-op}}$, as claimed.

Let us unpack some of the data encoded by~\eqref{eqn:functor3}.
Given a compactifiable $G$-morphism $f\colon Y\to X$ in $\Sch_B^G$, we denote by
\[
f_!\colon \SH^G(Y)\to \SH^G(X)
\]
the image of the span $Y\stackrel\id\from Y\stackrel f\to X$ by~\eqref{eqn:functor3}, with $S=B$. Being a morphism in $\Pr^{\mathrm L}$, $f_!$ admits a right adjoint $f^!$.
The functors $f_!$ and $f^!$ are called the \emph{exceptional functors}.
If $f$ is proper, there is an adjunction of spans
\[
(X\stackrel f\from Y\stackrel\id\to Y)\dashv (Y\stackrel\id\from Y\stackrel f\to X)
\]
 in $\Corr(\Sch_B^G)_{\all,\cpt}^{\prop,\operatorname{2-op}}$, so that $f_!\simeq f_*$. On the other hand, if $f$ is an open immersion, condition (\textasteriskcentered) implies that $f_!\simeq f_\sharp$. Thus, in general, we have an equivalence $f_!\simeq p_*j_\sharp$ for any factorization $f=pj$ with $j$ an open $G$-immersion and $p$ a proper compactifiable $G$-morphism.

A compactifiable $G$-morphism $f\colon Y\to X$ may also be viewed as a morphism in $\Sch_X^G$. By uniqueness of the extensions~\eqref{eqn:functor3}, there is a commutative square
\begin{tikzmath}
	\diagram{
	\Corr(\Sch_X^G)_{\all,\cpt} & \Mod_{\SH^G(X)} \\
	\Corr(\Sch_B^G)_{\all,\cpt} & \Mod_{\SH^G(B)}\rlap, \\
	};
	\arrows (11-) edge (-12) (11) edge (21) (12) edge (22) (21-) edge (-22);
\end{tikzmath}
which shows that $f_!$ can be promoted to an $\SH^G(X)$-module functor. In particular, there is a canonical equivalence
\[
f_!(\ph\tens f^*(\ph))\simeq f_!(\ph)\tens \ph.
\]
By construction, this equivalence is the smooth projection formula if $f$ is an open immersion and the proper projection formula if $f$ is proper.

Given a cartesian square of $G$-schemes
\begin{tikzequation}\label{eqn:BCsquare}
	\diagram{Y' & Y \\ X' & X \\};
	\arrows (11-) edge node[above]{$g$} (-12) (11) edge node[left]{$q$} (21) (21-) edge node[below]{$f$} (-22) (12) edge node[right]{$p$} (22);
\end{tikzequation}
with $p$ compactifiable, the functor~\eqref{eqn:functor3} gives a canonical equivalence
\[
\Ex^*_!\colon f^*p_! \simeq q_!g^*\colon \SH^G(Y) \to \SH^G(X').
\]
By construction, $\Ex^*_!$ can be identified with the exchange equivalence $\Ex^*_\sharp$ if $p$ is an open immersion and with the exchange equivalence $\Ex^*_*$ if $p$ is proper. By adjunction, there is also a natural equivalence
\[
\Ex^!_*\colon p^!f_* \simeq g_*q^! \colon \SH^G(X') \to \SH^G(Y).
\]

Given the cartesian square \eqref{eqn:BCsquare} with $p$ compactifiable, we can define a natural transformation
\[
\Ex_{!*}\colon p_!g_* \to f_*q_!\colon \SH^G(Y') \to \SH^G(X)
\]
by the composition
\[
p_!g_*\stackrel{\eta}\to f_*f^*p_!g_* \stackrel{\Ex^*_!}\simeq f_*q_!g^*g_* \stackrel\epsilon\to f_*q_!.
\]
It is an equivalence if $f$ is proper: this is obvious if $p$ is proper, and it follows from Corollary~\ref{cor:deligne} if $p$ is an open immersion. Similarly, we can define a natural transformation
\[
\Ex^{*!}\colon g^*p^!\to q^!f^*\colon \SH^G(X)\to\SH^G(Y')
\]
by the composition
\[
g^*p^! \stackrel{\eta}\to g^*p^! f_*f^* \stackrel{\Ex^!_*}\simeq g^*g_* q^!f^* \stackrel\epsilon\to q^!f^*.
\]
It is an equivalence if $f$ is smooth: this is obvious if $p$ is an open immersion, and it follows from Corollary~\ref{cor:deligne} if $p$ is proper.

Finally, given $f\colon X\to S$ compactifiable, we will define an endofunctor $\Tw_f\colon \SH^G(X)\to \SH^G(X)$ and natural transformations
\begin{gather*}
	f_!\to f_*,\\
	\Tw_f\circ f^! \to f^*.
\end{gather*}
Consider the cartesian square of $G$-schemes
\begin{tikzmath}
	\diagram{X\times_SX & X \\ X & S\rlap, \\};
	\arrows (11-) edge node[above]{$\pi_2$} (-12) (11) edge node[left]{$\pi_1$} (21) (21-) edge node[below]{$f$} (-22) (12) edge node[right]{$f$} (22);
\end{tikzmath}
and let $\delta\colon X\to X\times_SX$ be the diagonal. Since compactifiable morphisms are separated, $\delta$ is proper and hence $\delta_*=\delta_!$.
Then the transformation $f_!\to f_*$ is the composition
\[
f_! \simeq f_! \pi_{2*}\delta_* \xrightarrow{\Ex_{!*}} f_*\pi_{1!}\delta_*\simeq f_*.
\]
It is an equivalence if $f$ is proper, since $\Ex_{!*}$ is.
We set $\Tw_f=\delta^!\pi_2^*$. The natural transformation $\Tw_f\circ f^!\to f^*$ is then the composition
\[
\delta^!\pi_2^*f^! \xrightarrow{\Ex^{*!}} \delta^!\pi_1^!f^*\simeq f^*.
\]
It is an equivalence if $f$ is smooth, since $\Ex^{*!}$ is. In that case, $\Tw_f$ itself is an equivalence of $\infty$-categories, by Proposition~\ref{prop:thomequiv}, and in fact it is canonically equivalent to $\Sigma^{-\Omega_f}$, by purity (Proposition~\ref{prop:purity}).

The following theorem summarizes the properties of the six operations established so far:

\begin{theorem}\label{thm:main}
The six operations 
\[(\ph)^*,\; (\ph)_*,\; (\ph)_!,\; (\ph)^!,\; \tens,\; \Hom\]
satisfy the following properties, whenever the exceptional functors are defined.
	\begin{enumerate}
		\item \textnormal{(Proper pushforward)} If $f$ is a proper $G$-morphism, there is a canonical equivalence
		\[
		f_!\simeq f_*.
		\]
		\item \textnormal{(Smooth pullback)} If $f$ is a smooth $G$-morphism, there is a canonical equivalence
		\begin{gather*}
		\Sigma^{-\Omega_f}\circ f^!\simeq f^*.
		\end{gather*}
		\item \textnormal{(Base change)} If
		\begin{tikzmath}
				\diagram{\bullet & \bullet \\ \bullet & \bullet \\};
				\arrows (11-) edge node[above]{$g$} (-12) (11) edge node[left]{$q$} (21) (21-) edge node[below]{$f$} (-22) (12) edge node[right]{$p$} (22);
			\end{tikzmath}
			is a cartesian square of $G$-schemes, there are canonical equivalences
			\begin{gather*}
			f^*p_!\simeq q_!g^*,\\
			f^!p_*\simeq q_*g^!.
			\end{gather*}
		\item \textnormal{(Gluing)} If $i$ is a closed $G$-immersion with complementary open $G$-immersion $j$, there are cofiber sequences
		\begin{gather*}
		j_!j^! \to \id \to i_*i^*,\\
		i_!i^! \to \id \to j_*j^*.
		\end{gather*}
		\item \textnormal{(Immersive pushforward)} If $i$ is a $G$-immersion, the functors $i_*$ and $i_!$ are fully faithful.
		\item \textnormal{(Monoidality)} If $f$ is any $G$-morphism, there is a canonical equivalence
		\[
		f^*(\ph\tens\ph)\simeq f^*(\ph)\tens f^*(\ph).
		\]
		\item \textnormal{(Projection formulas)} If $f$ is any $G$-morphism, there are canonical equivalences
		\begin{gather*}
		f_!(\ph\tens f^*(\ph))\simeq f_!(\ph)\tens \ph ,\\
		\Hom(f_!(\ph),\ph)\simeq f_*\Hom(\ph, f^!(\ph)),\\
		f^!\Hom(\ph,\ph)\simeq \Hom(f^*(\ph), f^!(\ph)).
		\end{gather*}
		\item \textnormal{(Homotopy invariance)} If $f$ is a $G$-affine bundle, the functors $f^*$ and $f^!$ are fully faithful.
	\end{enumerate}
\end{theorem}

\begin{corollary}
	Let $f\colon Y\to X$ be a compactifiable $G$-morphism. Then the functor $f^!\colon \SH^G(X)\to\SH^G(Y)$ preserves colimits.
\end{corollary}

\begin{proof}
	By Nisnevich separation and Theorem~\ref{thm:main} (2), the question is Nisnevich-local on $B$. Hence, we can assume that $f$ is $G$-quasi-projective. By Lemma~\ref{lem:Gimmersion} (1), we can further assume that $f$ is smooth or a closed $G$-immersion. The result then follows from Theorem~\ref{thm:main} (2) and (4), respectively.
\end{proof}

\begin{remark}\label{rmk:stacks}
	As we explained in \S\ref{sub:equivariant}, our goal was to construct an interesting extension of the functor $\SH\colon \{\text{schemes}\}^\op\to\{\text{symmetric monoidal $\infty$-categories}\}$ to algebraic stacks.
	It is easy to see from the definitions (and Corollary~\ref{cor:thomequiv}) that the symmetric monoidal $\infty$-categories $\H^G(S)$, $\H^G_\pt(S)$, and $\SH^G(S)$ are indeed intrinsic invariants of the quotient stack $[S/G]$. 
	Moreover, it is clear that the adjunction $f^*\dashv f_*$, with $f^*$ symmetric monoidal, exists for an arbitrary morphism of stacks $f$. For example, if $f\colon\mathbf BG\to B$ is the unique map, $f_*\colon \SH(\mathbf BG)\to\SH(B)$ is the motivic analog of the ``genuine $G$-fixed points'' functor in stable equivariant homotopy theory. From this perspective, the base change property of Theorem~\ref{thm:main} admits the following generalization (with the same proof). Given any cartesian square of stacks
	\begin{tikzmath}
		\diagram{\mathfrak Y' & \mathfrak Y \\ \mathfrak X' & \mathfrak X \\};
		\arrows (11-) edge node[above]{$g$} (-12) (11) edge node[left]{$q$} (21) (21-) edge node[below]{$f$} (-22) (12) edge node[right]{$p$} (22);
	\end{tikzmath}
	for which $\SH(\ph)$, $p_!$, and $q_!$ are defined, there is a canonical equivalence $f^*p_!\simeq q_!g^*$.
\end{remark}

\begin{remark}\label{rmk:nagatachow}
	Suppose that the tame group $G$ is finite and discrete. In that case, $\SH^G(S)$ can be defined for arbitrary qcqs $G$-schemes $S$, using for $\Sm_S^G$ the category of all finitely presented smooth $G$-schemes over $S$ (see Remark~\ref{rmk:nissite}) and inverting the regular representation sphere (see Remark~\ref{rmk:regular}).
	The proper base change theorem and related results then hold for arbitrary proper $G$-morphisms. Following the proof of \cite[Proposition C.13]{HoyoisGLV}, we only need a suitable equivariant version of Chow's lemma \cite[XII, \S7]{SGA4-3}, but such a statement follows easily from its nonequivariant version. Indeed, suppose given a separated $G$-morphism of finite type $f\colon X\to S$ and a projective morphism $\pi\colon X'\to X$ such that $f\circ\pi$ is quasi-projective and $\pi^{-1}(U)\simeq U$ for some nonempty open subset $U\subset X$.
	Replacing $U$ by a maximal nonempty intersection of its $G$-translates, we can assume that
	$V=\coprod_{g\in G/H} gU$ is an open subscheme of $X$, where $H$ is the stabilizer of $U$.
	Replacing $X'$ by the fiber product of the $\lvert H\rvert$ maps $h\circ \pi\colon X'\to X$, $h\in H$, we can assume that $\pi$ is $H$-equivariant.\footnote{In terms of stacks, this fiber product is the Weil restriction of $X'$ along $X\to [X/H]$.}
	Let $Z\subset X$ be an $H$-invariant closed complement to $V\minus U$. Then one can replace $X'$ by the induced $G$-scheme $(G\times X'_Z)/H$ and $U$ by $V$.
	
	 Similarly, the functor $f_!$ can be defined for any $G$-morphism $f\colon X\to S$ that is separated and of finite type, since Nagata's compactification theorem immediately generalizes to the $G$-equivariant setting. Indeed, given a nonequivariant compactification $f=p\circ j$ where $j\colon X\into \bar X$, one obtains an equivariant one by replacing $\bar X$ by the schematic closure of $X$ in the fiber product of the $\lvert G\rvert$ maps $g\circ p\colon \bar X\to S$, $g\in G$ (the schematic closure will be $G$-invariant by Lemma~\ref{lem:Gimmersion} (2)).
	 
	 In summary, if $G$ is a finite discrete group of order $d$, Theorem~\ref{thm:main} applies to all qcqs $G$-schemes over $\Spec \Z[1/d]$, the exceptional adjunction $f_!\dashv f^!$ being defined for $f$ separated and of finite type.
\end{remark}

\subsection{Descent properties}
\label{sub:descent}

We investigate the descent properties of the functor $S\mapsto \SH^G(S)$. We begin by introducing the equivariant analogs of the cdh topology and the constructible topology.

Let $X$ be a $G$-scheme. An \emph{abstract blowup square} over $X$ is a cartesian square
\begin{tikzequation}\label{eqn:blowupsquare}
	\diagram{W & Y \\ Z & X \\};
	\arrows (11-) edge[c->] (-12) (11) edge (21)
	(21-) edge[c->] node[above]{$i$} (-22) (12) edge node[right]{$p$} (22);
\end{tikzequation}
of $G$-schemes where $i$ is a closed $G$-immersion, $p$ is proper, and $p$ induces an isomorphism $Y\times_X(X\minus Z)\simeq X\minus Z$. 
The \emph{cdh topology} on $\Sch_B^G$ is the coarsest topology finer than the Nisnevich topology and such that, for every abstract blowup square~\eqref{eqn:blowupsquare}, $\{i,p\}$ generates a covering sieve of $X$.
The same argument as in the proof of Proposition~\ref{prop:topos} shows that a presheaf $F$ on $\Sch_B^G$ is a sheaf for the cdh topology if and only if $F(\emptyset)$ is contractible and $F$ sends Nisnevich squares and abstract blowup squares to cartesian squares.

The \emph{constructible topology} on $\Sch_B^G$ is the coarsest topology such that:
\begin{itemize}
	\item the empty sieve covers the empty scheme;
	\item if $Z\into X$ is a closed $G$-immersion with open complement $U\into X$, $\{U\into X,Z\into X\}$ generates a covering sieve.
\end{itemize}
Note that the constructible topology is finer than the cdh topology.

\begin{proposition}
	\label{prop:constructible}
	Let $\{f_i\colon U_i\to S\}$ be a constructible cover of a $G$-scheme $S$. Then the families of functors
\begin{gather*}
\{f_i^*\colon \SH^G(S)\to \SH^G(U_i)\},\\
\{f_i^!\colon \SH^G(S) \to \SH^G(U_i)\}
\end{gather*}
are conservative (assuming each $f_i$ compactifiable in the latter case).
\end{proposition}

\begin{proof}
	This is an immediate consequence of gluing.
\end{proof}

\begin{proposition}\label{prop:cdh}
	The functors
	\begin{gather*}
		(\Sch_B^G)^\op\to\Cat_\infty,\quad S\mapsto \SH^G(S),\quad f\mapsto f^*,\\
		(\Sch_B^G)_\cpt^\op\to\Cat_\infty,\quad S\mapsto \SH^G(S),\quad f\mapsto f^!,
	\end{gather*}
	are sheaves for the cdh topology.
\end{proposition}

\begin{proof}
	We give the proof for the first functor.
	We must show that:
	\begin{enumerate}
		\item[(a)] $\SH^G(\emptyset)\simeq *$;
		\item[(b)] $\SH^G(\ph)$ takes Nisnevich squares to cartesian squares;
		\item[(c)] $\SH^G(\ph)$ takes abstract blowup squares to cartesian squares.
	\end{enumerate}
	Assertion (a) is obvious. 
	We will prove (c) and omit the proof of (b) which is entirely similar.\footnote{Assertion (b) also follows from Proposition~\ref{prop:Nislocal} (2), but this alternative argument is not available for the second functor.}
	Let $Q$ be the abstract blowup square
	\begin{tikzmath}
		\diagram{W & Y \\ Z & X\rlap. \\};
		\arrows (11-) edge[c->] node[above]{$k$} (-12) (11) edge node[left]{$q$} (21)
		(21-) edge[c->] node[above]{$i$} (-22) (12) edge node[right]{$p$} (22);
	\end{tikzmath}
	By \cite[Lemma 5.17]{DAG7}, $\SH^G(Q)$ is cartesian if and only if:
	\begin{enumerate}
		\item[(d)] the pair $(i^*,p^*)$ is conservative;
		\item[(e)] given $E_Z\in\SH^G(Z)$, $E_Y\in\SH^G(Y)$, $E_W\in\SH^G(W)$, and $q^*E_Z\simeq E_W\simeq k^*E_Y$, if \[E_X=i_*E_Z\times_{(pk)_*E_W}p_*E_Y,\] then the maps
		\[
		i^*E_X \to E_Z\quad\text{and}\quad p^*E_X\to E_Y
		\]
		induced by the canonical projections are equivalences.
	\end{enumerate}
	Assertion (d) follows from Proposition~\ref{prop:constructible}, since $\{i,p\}$ is a constructible cover of $X$.
	Let us prove (e).
	Proper base change and the full faithfulness of $i_*$ immediately imply that $i^*E_X\to E_Z$ is an equivalence. To show that $p^*E_X\to E_Y$ is an equivalence, it suffices to show that $k^*p^*E_X\to k^*E_Y$ and $j^*p^*E_X\to j^*E_Y$ are equivalences, where $j$ is the open $G$-immersion complementary to $k$. 
	The former is $q^*$ of the map just shown to be an equivalence, and the latter is an equivalence by smooth base change.
\end{proof}

\begin{corollary}\label{cor:descent}
	Let $S$ be a $G$-scheme and let $E\in\SH^G(S)$. Then the functor
	\[
	(\Sch_S^G)^\op \to \SH^G(S),\quad (f\colon X\to S)\mapsto f_*f^*E,
	\]
	is a sheaf for the cdh topology, and the functor
	\[
	(\Sch_B^G)_{\cpt/S} \to \SH^G(S),\quad (f\colon X\to S)\mapsto f_!f^!E,
	\]
	is a cosheaf for the cdh topology.
\end{corollary}

In particular, taking mapping spaces or mapping spectra, we deduce that any space-valued or spectrum-valued invariant of $G$-schemes represented by a cocartesian section of $\SH^G(\ph)$ satisfies cdh descent.

\begin{remark}\label{rmk:generalcdh}
	If $G$ is finite and discrete, the results of this subsection hold if we replace $\Sch_B^G$ by the category of all qcqs $G$-schemes over $B$ and $(\Sch_B^G)_\cpt$ by the wide subcategory of separated $G$-morphisms of finite type (see Remark~\ref{rmk:nagatachow}).
\end{remark}

\providecommand{\bysame}{\leavevmode\hbox to3em{\hrulefill}\thinspace}

\end{document}